\newif\ifCRM

\ifCRM
  \documentclass[a4paper,reqno]{amsart}
  \usepackage{a4}
\else
  \documentclass[a4paper,twoside,12pt]{article}
  \usepackage[a4paper,vmargin={3cm,4cm},hmargin={4cm,3cm}]{geometry}
\fi

\usepackage[latin1]{inputenc}
\usepackage[T1]{fontenc}
\usepackage{textcomp}

\usepackage[english]{babel}

\usepackage{amsmath,amssymb,amsfonts,amsthm,dsfont}

\usepackage{graphicx,array}

\title{Condensation of polyhedric structures onto soap films}
\author{Vincent Feuvrier}
\ifCRM
  \address{CRM\\Apartar 50\\08193 Bellaterra\\Barcelona, Spain}
  \email{vincent.feuvrier@normalesup.org}
  \subjclass[2000]{49Q20}
  \keywords{Almgren Minimal sets, Plateau problem, Euclidean polyhedrons and complexes}
\fi

\newtheorem{definition}{Definition}
\newtheorem{proposal}{Proposal}
\newtheorem{theorem}{Theorem}
\newtheorem{lemma}{Lemma}

\newtheorem{corollary}{Corollary}

\def\open#1{\overset{\circ}{#1}}
\def\dist{\mathbf{d}}
\def\Hconverge#1{\stackrel{#1}{\displaystyle\relbar\joinrel\rightharpoondown}}
\DeclareMathOperator\diam{Diam}
\DeclareMathOperator\vect{Vect}
\DeclareMathOperator\affine{Affine}

\DeclareMathOperator\identity{Id}
\DeclareMathOperator\support{Spt}
\def\subsubset{\subset\kern-0.2em\subset}
\def\A{\mathcal{A}}
\def\B{\mathcal{B}}
\def\C{\mathcal{C}}
\def\F{\mathcal{F}}

\def\H{\mathcal{H}}

\def\R{\mathcal{R}}
\def\S{\mathcal{S}}

\def\U{\mathcal{U}}

\def\Int{\mathop{\int}}

\newlength\flsplitlabelwidth
\newlength\flsplitwidth
\newenvironment{flsplit}{%
  \settowidth\flsplitlabelwidth{\strut\hspace{2em}(\theequation)}%
  \setlength\flsplitwidth{\columnwidth}%
  \addtolength\flsplitwidth{-\flsplitlabelwidth}%
  \begin{equation}%
    \begin{array}{@{}>{\displaystyle}r@{\;}>{\displaystyle}l@{}}%
      \multicolumn{2}{@{}c@{}}{\strut\hspace{\flsplitwidth}\strut}\\[-2.5ex]%
}{%
    \end{array}%
  \end{equation}%
}

\def\shovebreak#1{
  \\&
  \multicolumn{1}{@{}>{\displaystyle}r@{}}{#1}\\[5ex]%
}

\allowdisplaybreaks[1]

\begin{document}

\ifCRM\else
  \maketitle
  \vfill
\fi

\begin{abstract}
We study the existence of solutions to general measure-minimization problems over topological classes that are stable under localized Lipschitz homotopy, including the standard Plateau problem without the need for restrictive assumptions such as orientability or even rectifiability of surfaces. In case of problems over an open and bounded domain we establish the existence of a ``minimal candidate'', obtained as the limit for the local Hausdorff convergence of a minimizing sequence for which the measure is lower-semicontinuous. Although we do not give a way to control the topological constraint when taking limit yet --- except for some examples of topological classes preserving local separation or for periodic two-dimensional sets --- we prove that this candidate is an Almgren-minimal set. Thus, using regularity results such as Jean Taylor's theorem, this could be a way to find solutions to the above minimization problems under a generic setup in arbitrary dimension and codimension.
\end{abstract}

\ifCRM
  \maketitle
\else
  \vfill
  \strut
  \newpage
  \strut
  \vfill
\fi

\ifCRM
  \enlargethispage{4ex}
\fi

\tableofcontents

\ifCRM\else
  \vfill
  \strut
  \newpage
\fi

\section*{Introduction}
\ifCRM\else
  \addcontentsline{toc}{section}{Introduction}
\fi

We consider a class $\mathfrak{F}$ of relatively closed subsets of a given domain $U$ in $\mathbb{R}^n$ --- that will be our competitors, and we also suppose that $\mathfrak{F}$ is stable under some class of admissible deformations (see definition~\ref{definitiondeformation}).

We then consider the following problem: find $E\in\mathfrak{F}$ such that
\begin{equation}
\mu(E)=\inf_{F\in\mathfrak{F}}\mu(F),
\end{equation}
where $\mu$ stands for a given $d$-dimensional measure functional with $0\leq d<n$ --- for instance the $d$-dimensional Hausdorff measure $\H^d$, but more general cases are also possible. The Plateau problem can be rewritten in these terms, by taking a class $\mathfrak{F}$ stable under Lipschitz deformations that only move a relatively compact subset of points of $U$. In that case, the boundary of $U$ acts as a topological constraint.

In case of a problem over an open bounded domain $U$ of $\mathbb{R}^n$ and in arbitrary dimension and codimension we prove the following theorem of existence of a minimal candidate (see theorem~\ref{theoremexistence} for a more precise statement):

\begin{quotation}\itshape
There is a relatively closed subset $E$ of $U$, Almgren almost-minimal and with no greater measure than any element of $\mathfrak{F}$, that is obtained as a local Hausdorff limit over all compact subsets of $U$ of a measure-minimizing sequence of elements of $\mathfrak{F}$.
\end{quotation}

Notice that we do not prove that $E\in\mathfrak{F}$ --- in fact it can be false, see section~\ref{subsectionapplication} where we also give two examples of usage of this result. However, we hope that in some cases at least for $2$-dimensional sets, by using regularity-related results about $E$ such as Jean Taylor theorem (see~\cite{taylor,david2,david3}) we may be able to build a Lipschitz retraction sending a neighborhood of $E$ onto $E$, which would be enough to control the topological constraint in $\mathfrak{F}$ when taking limit in our minimizing sequence.

One of the technical difficulties that arise in this approach is that the Hausdorff measure is generally not lower semicontinuous --- although the case of one-dimensional sets can be handled using Go\l\k ab's theorem --- which usually prevents directly taking limit in arbitrary minimizing sequences to study the existence of solutions to this kind of general, measure-related minimization problems.

In fact, we give a way to convert any measure-minimizing sequence into another minimizing sequence of ``regularized'' sets (i.e. quasiminimal with uniform constants) that verify an uniform concentration property initially introduced by Dal Maso, Morel and Solemini in~\cite{DMS}, and for which the Hausdorff measure is lower semicontinuous (see theorem~\ref{theoremdavid}, which is borrowed from~\cite{david1}).

The first step of this process is to find a way to build generalized Euclidean dyadic grids with several imposed orientation and uniform bounds on the flatness of their polyhedrons. Their construction is explained in \cite{feuvrier2} (see theorem~\ref{theoremfusion}): provided that they are far enough from each others, it is possible to glue several dyadic grids (with different orientations) together into a larger grid of convex polyhedrons that ``connect well'' (see definition~\ref{definitioncomplex} for a topological definition) and such that every polyhedron of the new grid (including its faces in all lower dimensions) is not too flat. In fact, we give an implicit uniform lower bound that depends only on dimension $n$ on the minimal angle of two faces of any dimension that meet at a given vertex (see our definition~\ref{definitionrotondity} of ``rotondity'').

\begin{figure}
\includegraphics[width=0.45\textwidth]{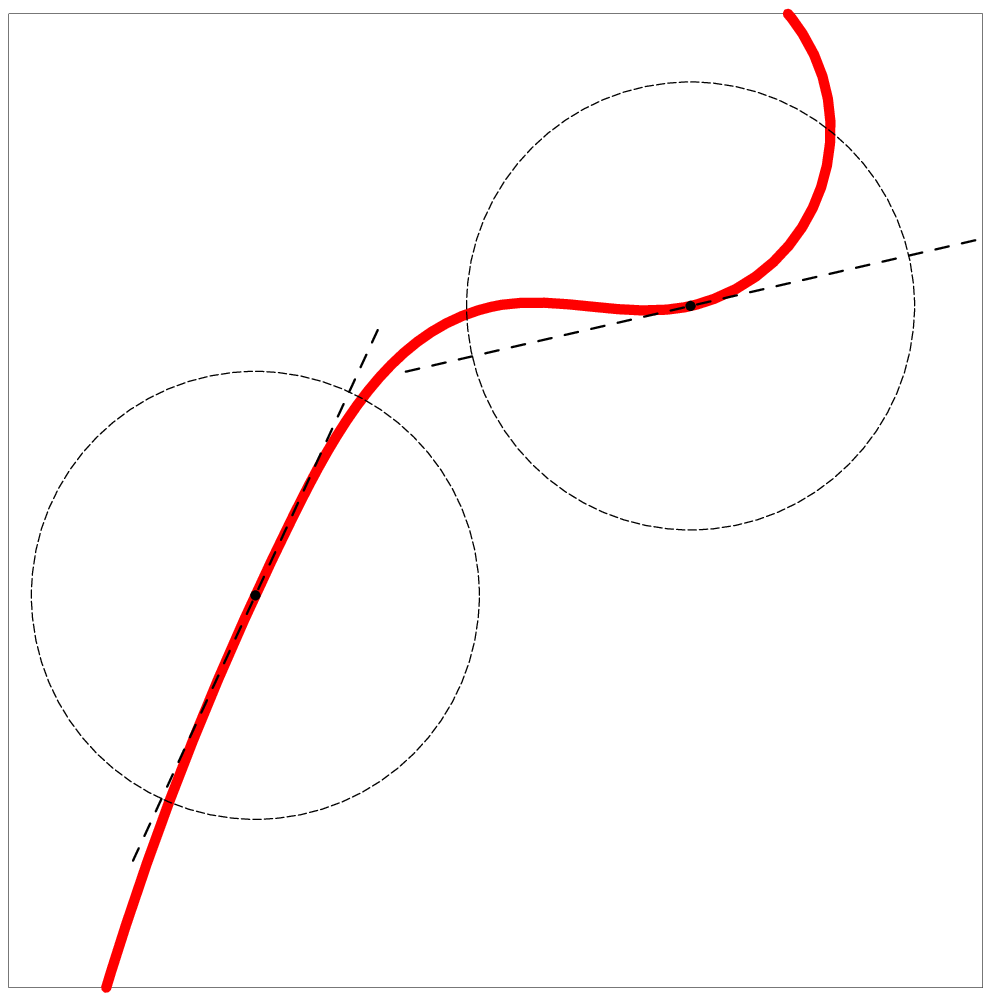}\hfill\includegraphics[width=0.45\textwidth]{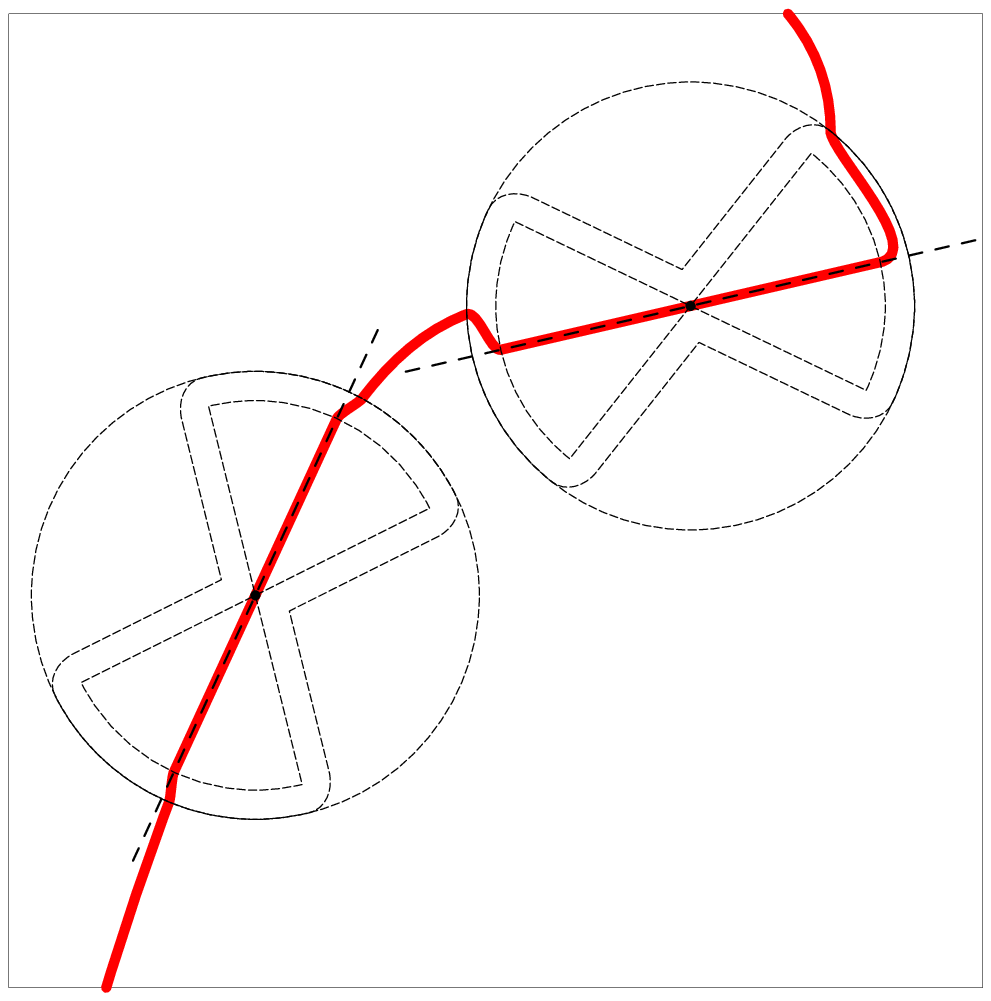}
\caption{\label{figureintroductionA}On the left, a competitor in $\mathfrak{F}$ and an almost-covering by disjoint balls centered on its rectifiable part. On the right, we project it onto some of its tangent planes inside a cone and a ball, while keeping the measure of the patches that connect the flat part to the remaining one arbitrary small.}
\end{figure}

\begin{figure}
\includegraphics[width=0.45\textwidth]{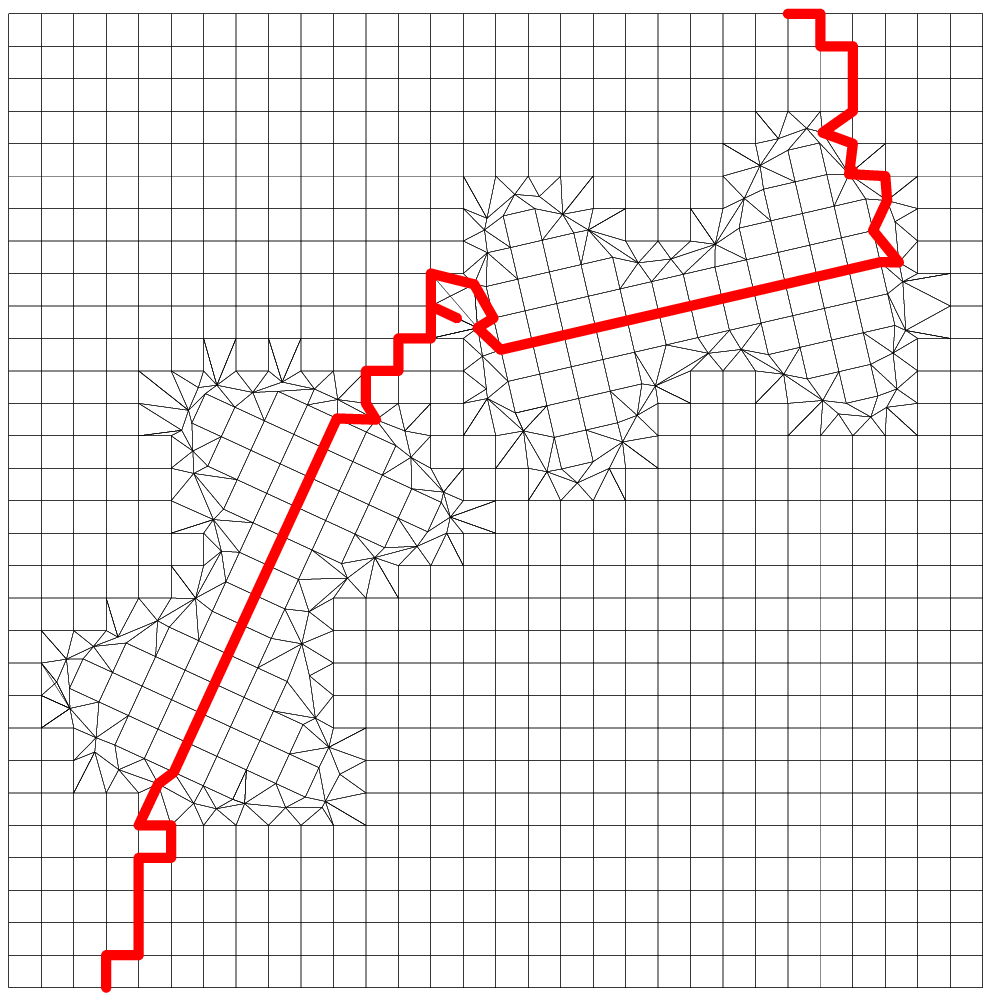}\hfill\includegraphics[width=0.45\textwidth]{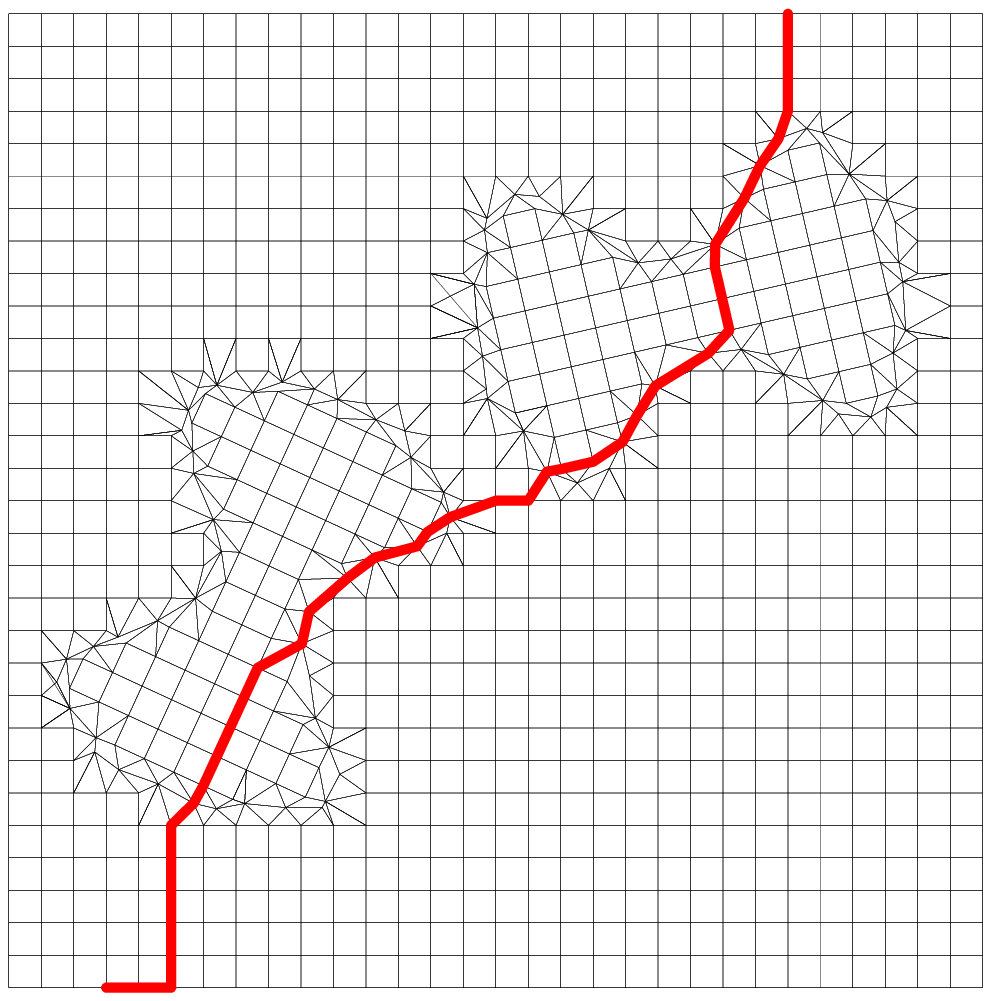}
\caption{\label{figureintroductionB}On the left, we use Federer-Fleming radial projections inside a polyhedric grid designed to keep the measure increase as small as needed with respect to the initial competitor on figure~\ref{figureintroductionA}. Notice that the measure did not increase in the cubes that are parallel to the tangent planes. On the right, we do a finite measure-minimization amongst the polyhedric competitors in $\mathfrak{F}$. The set we obtain is quasiminimal with constants depending on the flatness of the polyhedrons of the grid, and is even better than the polyhedric competitor on the left.}
\end{figure}

The second step is to carefully design polyhedric grids to approximate a given compact $d$-dimensional set while keeping control on the measure increase introduced by the approximation (see theorem~\ref{theorempolyhedralapproximation}). For this purpose, we use an almost covering of the rectifiable part of the set by dyadic grids that roughly follow the direction of its tangent planes and then use the above method to merge these grids together (see figure~\ref{figureintroductionA}). The uniform lower bound obtained on the flatness of the polyhedrons is useful when approximating our sets using successive Lipschitz Federer-Fleming-like projections (see~\cite{federer-flemming}) onto decreasing dimensional polyhedrons of the grid till dimension $d$, to obtain additional measure-related regularity constants (in fact, quasiminimality constants, as introduced earlier by Almgren) that depend only on dimensions $d$ and $n$ (see figure~\ref{figureintroductionB}).

This polyhedral approximation theorem is the key result of this paper. One can see it as a version for non-orientable surfaces of the classic polyhedral approximation theorem for integral currents. It may also be used to generalize to higher dimension and codimension a result of T. De Pauw in \cite{depauw} for two-dimensional surfaces in $\mathbb{R}^3$.

The plan of the paper is as follows.

Section~1 is devoted to summarize the basic definitions and notations we will be using through the next sections. We start with Euclidean polyhedrons, complexes and dyadic cubes. We also give an Almgren-like formalism (see \cite{almgren,david1}) for quasiminimal, almost-minimal and minimal sets.

In section~2 we give some technical lemmas that are to be used later in the polyhedral approximation process. First we give some Lipschitz extension lemmas before studying basic measure-related properties of orthogonal and radial projection extensions.

In section~3 we give an optimization lemma which allows converting any competitor into another one that is quasiminimal with constants depending only on the dimension, without increasing its measure too much. Then, we proceed in proving the main theorem, before giving some examples of setup under which the topological constraint behaves well when taking limit.

The proposed research of solutions is actually quite close in spirit to that of Reifenberg (see \cite{reifenberg}), although based on Almgren's initial formalism. It is not as ``elementary'' and flexible as any of the classic distributional approaches, but fits problems that cannot be handled by currents and finite perimeter sets. Compared to Reifenberg theory, it might end up to be simpler and more flexible because it heavily relies on technical geometric tools which involve long proofs and complicated constructions but hopefully will be turned into ready-to-use results.

\ifCRM\else
  The author would like to express his thanks to Guy David for his many advices and suggestions. He also gratefully acknowledges partial support from the Centre de Recerca Matemàtica at the Universitat Autònoma de Barcelona.
\fi

\section{Preliminaries}

We begin with some notations and basic definitions.

\subsection{Euclidean polyhedrons}

We place ourself in $\mathbb{R}^n$ with its usual Euclidean structure. We say that a set $A$ is an affine half-space if one can find an affine hyperplane $H$ and a non-parallel vector $u$ such that
\begin{equation}
A=\left\{x+ru\colon x\in H\text{ and }r\geq 0\right\}.
\end{equation}
We will say that a non-empty intersection of affine half-spaces is a polyhedron according to the following definition.

\begin{definition}[Polyhedrons]\label{definitionpolyhedron}
A polyhedron $\delta$ of dimension $n$ is a compact with non-empty interior intersection of finitely many affine half-spaces.
\end{definition}

By keeping only affine half-spaces whose boundary intersects $\delta$ over a set of $n-1$ Hausdorff dimension it is easy to check that amongst all half-spaces families that are suitable for this definition one can find one that is minimal for inclusion. We will denote it by $\A(\delta)$.

By allowing non-empty compact sets with empty interior we generalize this definition to $k$-dimensional polyhedrons (with $k\leq n$) by placing ourselves in the smallest affine subspace $\affine(\delta)$ of dimension $k$ that contains them. In that case, the usual topological operators (closure, interior and boundary) will be taken relatively to $\affine(\delta)$, as well as the affine half-subspaces in $\A(\delta)$. By convention we consider singletons as polyhedrons of dimension zero, equal to their interior and with empty boundary.

Polyhedrons as we defined them are convex. With a simple convexity argument it is easy to check that the affine dimension of $\affine(\delta)$ is the same as the Hausdorff dimension of $\delta$. We will denote both by $\dim(\delta)$.

In fact, it is possible to show (but we will not do it here) that our definition is equivalent to the one of usual convex polytopes, as the convex hull of a finite family of points --- typically the ``vertices'', that we will introduce shortly. Indeed, the previous notations allow for an easy writing of the definition of polyhedric faces. For convenience we will call them ``subfaces'' in the general case and keep the word ``face'' to specifically designate a subface of dimension one less than the relative polyhedron.

\begin{definition}[Subfaces]
Let $\delta$ be a $n$-dimensional polyhedron such that $\A(\delta)=\left\{A_1,\ldots,A_p\right\}$ and $\left\{A'_1,\ldots,A'_p\right\}$ a family of subsets of $\mathbb{R}^n$ such that $A'_i=A_i$ or $A'_i=\partial A_i$ for $1\leq i\leq p$. By putting $\alpha=\bigcap_i A'_i$, if $\alpha\neq\emptyset$ we say that $\alpha$ is a subface of $\delta$ and more precisely:
\begin{itemize}
\item if $\dim\alpha<\dim\delta$ then $\alpha$ is a strict subface;
\item if $\dim\alpha=\dim\delta-1$ then $\alpha$ is a face;
\item if $\dim\alpha=0$ (i.e. if $\alpha$ is a singleton) then $\alpha$ is a vertex and we will mistake it for the point it contains for convenience.
\end{itemize}
We will denote by $\F(\delta)$ the set of all subfaces of $\delta$, and for $0\leq k\leq\dim\delta$:
\begin{equation}
\F_k(\delta)=\left\{\alpha\in\F(\delta)\colon\dim\alpha=k\right\}.
\end{equation}
\end{definition}

Again, we naturally generalize this definition to $k$-dimensional polyhedrons with $k\leq n$. It is not difficult to check that subfaces are also polyhedrons, that the faces are of disjoint interior and that their union is the boundary of the polyhedron. For any polyhedron $\delta$ we can even write that
\begin{equation}
\delta=\bigsqcup_{\alpha\in\F(\delta)}\open{\alpha}
\end{equation}
where $\sqcup$ stands for a disjoint union and the interior $\open{\alpha}$ of all subfaces is taken relatively to the corresponding generated affine subspace $\affine(\alpha)$.

We now give ourselves some way to control the flatness of polyhedrons, which will be used later to control the measure increase when approximating rectifiable sets using radial projections onto them.

\begin{definition}[Shape control]\label{definitionrotondity}
For any non-empty compact set $A$ we define the following quantities:
\begin{itemize}
\item the outer radius, by taking the infimum of radii of balls containing $A$ (with the convention $\inf\emptyset=0$)
\begin{equation}
\overline{R}(A)=\inf\left\{r>0\colon\exists x\in\mathbb{R}^n,A\subset B(x,r)\right\};
\end{equation}
\item the inner radius, by taking the supremum of radii of included balls (with the convention $\sup\emptyset=0$)
\begin{equation}
\underline{R}(A)=\sup\left\{r>0\colon\exists x\in\mathbb{R}^n,B(x,r)\cap\affine(A)\subset A\right\};
\end{equation}
\item the rotondity, by taking the ratio of the two (with the convention $R(A)=1$ when $\overline{R}(A)=0$)
\begin{equation}
R(A)=\frac{\underline{R}(A)}{\overline{R}(A)}\in[0,1].
\end{equation}
\end{itemize}
\end{definition}

Of course, the more $R(A)$ is close to $1$, the more $A$ look like a ball and the less it is flat. By a compacity argument, it is easy to show that the supremum in the calculus of $\underline{R}(A)$ is reached for some ball $B$ such that $B\cap\affine A\subset A$. We will call it an inscribed ball inside $A$.

\subsection{Polyhedric complexes and dyadic cubes}

We now consider a finite set $S$ of $k$-dimensional polyhedrons. We introduce the following notations:
\begin{itemize}
\item the union of the polyhedrons
\begin{equation}
\U(S)=\bigcup_{\delta\in S}\delta;
\end{equation}
\item the set of the subfaces
\begin{equation}
\F(S)=\bigcup_{\delta\in S}\F(\delta).
\end{equation}
\end{itemize}
Additionally, when all the polyhedrons in $S$ have the same dimension $k$ we will also use:
\begin{itemize}
\item the set of $k'$-dimensional subfaces (for $0\leq k'\leq k$)
\begin{equation}
\F_{k'}(S)=\bigcup_{\delta\in S}\F_{k'}(\delta);
\end{equation}
\item the set of boundary faces
\begin{equation}
\F_\partial(S)=\left\{\alpha\in\F_{k-1}(S)\colon\forall(\beta,\gamma)\in S^2,\alpha\neq\beta\cap\gamma\right\}.
\end{equation}
\end{itemize}

To formalize the idea of polyhedric meshes made of polyhedrons that connect well we give the following definition.

\begin{definition}[Complexes]\label{definitioncomplex}
We say that a set $S$ of $k$-dimensional polyhedrons is a $k$-dimensional complex if all its subfaces have disjoint interiors (again, relative to the generated affine subspace):
\begin{equation}
\forall(\alpha,\beta)\in\F(S)^2\colon\alpha\neq\beta\Rightarrow\open{\alpha}\cap\open{\beta}=\emptyset.
\end{equation}
\end{definition}

For instance, it is easy to check that for any polyhedron $\delta$ and $0\leq k\leq\dim\delta$, the set $\F_k(\delta)$ is a complex. So is $\F_{k'}(S)$ for $0\leq k'\leq k$ when $S$ is a $k$-dimensional complex. Furthermore, when $k=n$ we also have $\partial(\U(S))=\U(\F_\partial(S))$. When $S$ is a complex, we call any subset of $\F(S)$ made of subfaces of dimension at most $k$ a ``$k$-dimensional skeleton'' of $S$.

To control the shape of all polyhedrons --- including their subfaces --- within a complex we also generalize our notations for inner or outer radii and rotondity to complexes as well:
\begin{equation}\label{equationcomplexrotondity}
\overline{\R}(S)=\max_{\delta\in\F(S)}\overline{R}(\delta)\qquad\underline{\R}(S)=\min_{\delta\in\F(S)}\underline{R}(\delta)\qquad\R(S)=\min_{\delta\in\F(S)}R(\delta).
\end{equation}

Generic and easy-to-use examples of complexes are those made of dyadic cubes. For $r>0$ a dyadic cube is a polyhedron that can be written as $[0,r]^n$ in some orthonormal basis of $\mathbb{R}^n$, and an unit dyadic cube when $r=1$. Such cubes can be naturally placed on a grid to form a complex.

\begin{definition}[Dyadic complexes]\label{definitiondyadiccomplex}
We call dyadic complex of stride $r$ any set of dyadic cubes that can be written as
\begin{equation}
S=\left\{rz+[0,r]^n\colon z\in Z\right\}
\end{equation}
in an orthonormal basis, where $Z$ is a finite subset of $\mathbb{Z}^n$.
\end{definition}

Dyadic cubes are very convenient to locally approximate rectifiable sets of arbitrary dimension using their subfaces, because we can always choose their orientations to locally match those of the set's tangent planes while taking them as small as needed. To closely match an arbitrary set we would end up with many disjoint dyadic complexes with different orientations. Then, to complete the polyhedral approximation process these complexes should be merged into a larger one that covers the entire set to be approximated. However, although anyone would believe that such polyhedrons can be built it is not obvious that the non-dyadic polyhedrons needed to fill the gaps between all the dyadic grids can always be designed so they are never too flat.

\begin{figure}
\includegraphics[width=\textwidth]{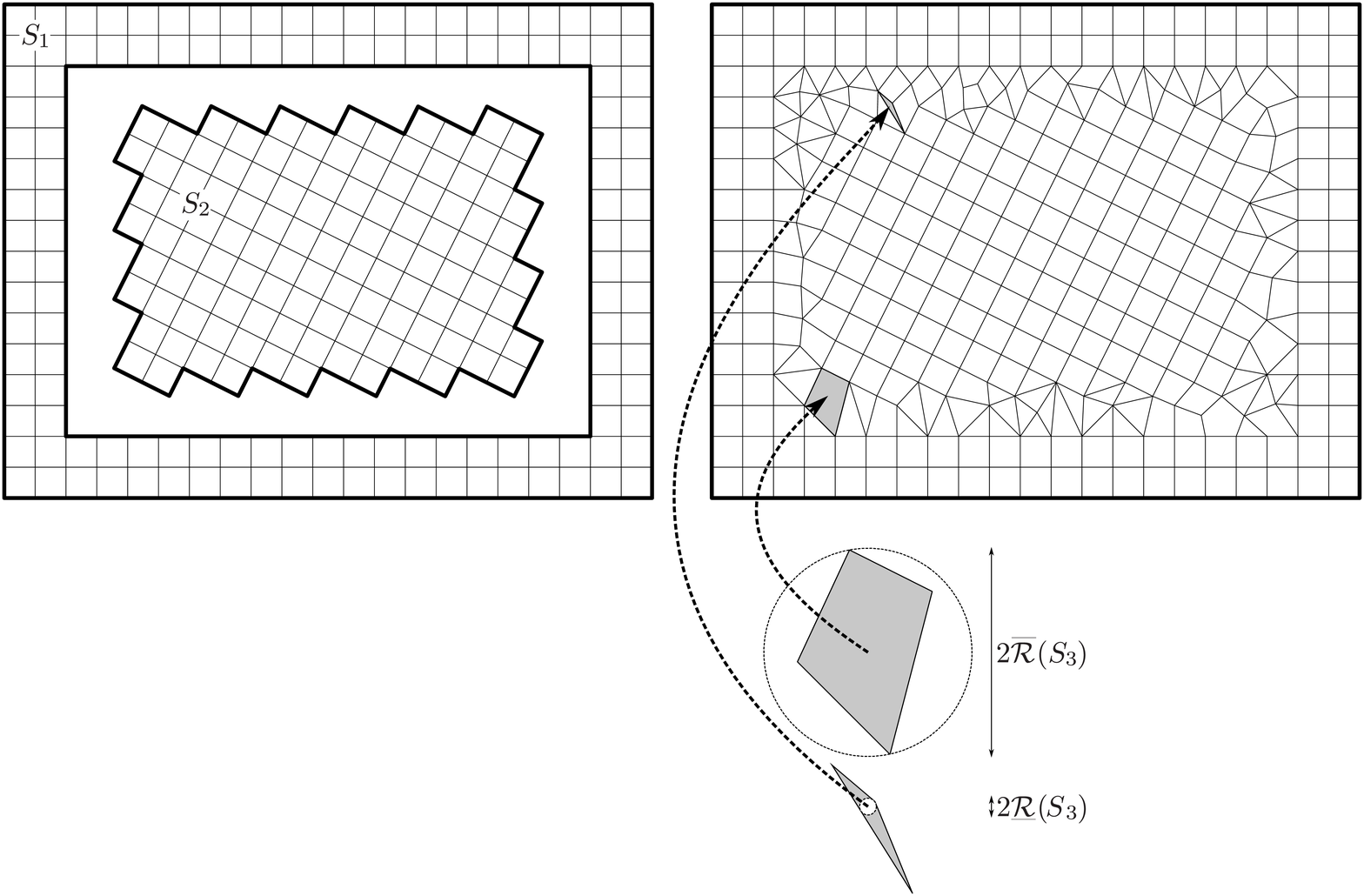}
\caption{\label{figurefusionA}Merging of two dyadic complexes with different orientations, and the associated shape constants.}
\end{figure}

In~\cite{feuvrier2} we proved the following result that can be used to merge two dyadic complexes together while keeping uniform bounds on the rotondity of all added polyhedrons and their subfaces (see figure~\ref{figurefusionA}).

\begin{theorem}[Merging of dyadic complexes with uniform rotondity]\label{theoremfusion}
One can find three positive constants $\rho$, $c_1$ and $c_2$ depending only on $n$ such that for all compact set $K$, for all open set $O\subset K$ and for all unit dyadic complexes $S_1$ and $S_2$ such that
\begin{equation}
\U(S_1)=K\setminus O\qquad\U(S_2)\subset O\qquad\min_{(x,y)\in\U(S_1)\times\U(S_2)}\Vert x-y\Vert\geq\rho
\end{equation}
then one can build $S_3$ such that $S'=S_1\sqcup S_2\sqcup S_3$ is a $n$-dimensional complex verifying
\begin{equation}
\U(S')=K\qquad\overline{\R}(S')\leq c_1\overline{\R}(S_1\cup S_2)\qquad\underline{\R}(S')\geq c_2\underline{\R}(S_1\cup S_2).
\end{equation}
\end{theorem}

Later, we will use this theorem to merge a large number of disjoint dyadic grids of arbitrary orientations together --- assuming their stride is small enough to build it --- by considering a global dyadic grid with ``holes'' separately enclosing each one.

\subsection{Quasiminimal and (almost-)minimal sets}

Let $U$ be a nonempty domain of $\mathbb{R}^n$. For a map $f\colon U\rightarrow U$ we denote by $\xi_f$ the set of points that are actually moved by $f$:
\begin{equation}
\xi_f=\left\{x\in U\colon x\neq f(x)\right\}.
\end{equation}
We also call support of $f$ the set of these points and their images:
\begin{equation}
\support f=\xi_f\cup f(\xi_f).
\end{equation}

Suppose that $M\geq 1$. In what follows, we assume that we are given a measurable function $h$ over $U$, with values in $[1,M]$. For $0\leq d<n$ we will consider the following $d$-dimensional set functional, for any measurable set $E\subset U$:
\begin{equation}\label{equationfunctionnal}
J_h^d(E)=\int_{x\in E}h(x)d\H^d(x)
\end{equation}
where $\H^d$ stands for the $d$-dimensional Hausdorff measure (see for instance Mattila's book~\cite{mattila}).

The following definition will be useful to describe our so-called ``topological classes stable under local Lipschitz homotopy''.

\begin{definition}[Admissible deformations]\label{definitiondeformation}
For $\delta>0$ we say that a one-parameter family $(\phi_t)_{t\in[0,1]}$ of maps from $U$ into itself is a $\delta$-deformation over $U$ if the following requirements are met:
\begin{itemize}
\item$\phi_0=\identity_U$ and $\phi_1$ is Lipschitz;
\item$(t,x)\mapsto\phi_t(x)$ is continuous over $[0,1]\times U$;
\item by putting
\begin{equation}\support\phi=\bigcup_{t\in[0,1]}\support\phi_t\end{equation}
then $\support\phi$ is compact relatively in $\open U$ (i.e. $\overline{\support\phi}$ is compact and included in $U$, which we will denote by $\support\phi\subsubset\open U$) and $\diam(\support\phi)\leq\delta$.
\end{itemize}
\end{definition}

When $(\phi_t)$ is a deformation over $U$ and $E\subset U$ we say that $\phi_1(E)$ is an Almgren competitor of $E$.

For $X\subset\mathbb{R}^n$ and $\rho>0$ we denote by $X_\rho$ the $\rho$-neighborhood of $X$:
\begin{equation}
X_\rho=\bigcup_{x\in X}B(x,\rho)=\left\{x\in\mathbb{R}^n\colon\dist(x,X)<\rho\right\}.
\end{equation}
For convenience we give the following statement which will be used later to easily build a deformation from a Lipschitz map whose support is small enough.

\begin{proposal}[Automatic building of deformation]\label{proposalautomaticdeformation}
Suppose that $U\subset\mathbb{R}^n$, that $f$ is a Lipschitz map over $U$ and that $(\phi_t)$ is a $\diam(U)$-deformation over $U$. If there is $\rho>0$ such that
\begin{equation}
\Vert\phi_1-f\Vert_\infty<\rho\qquad\text{and}\qquad(\xi_{\phi_1}\cup\xi_f)_\rho\subsubset U
\end{equation}
then the one-parameter family $(\psi_t)$ of maps on $U$ defined for $0\leq t\leq 1$ by
\begin{equation}
\psi_t(x)=\begin{cases}
\phi_{2t}(x)&\text{if $t\leq\frac{1}{2}$}\\
(2-2t)\phi_{1}(x)+(2t-1)f(x)&\text{if $t>\frac{1}{2}$}
\end{cases}
\end{equation}
is also a $\diam(U)$-deformation over $U$ such that $\psi_1=f$.
\end{proposal}

The proof is really easy, and consists only in proving that $\overline{\support\psi}$ is relatively compact in $\open U$.

\begin{proof}
Suppose that $x\in U$ and consider the three possible cases:
\begin{itemize}
\item if $x\notin\xi_{\phi_1}\cup\xi_f$ then
\begin{equation}\label{equationautomaticdeformationA}
\left\{\psi_t(x)\colon t\in[0,1]\right\}=\left\{\phi_t(x)\colon t\in[0,1]\right\};
\end{equation}
\item if $x\in\xi_{\phi_1}$ then for all $t\in[0,1/2]$:
\begin{equation}\label{equationautomaticdeformationB}
\psi_t(x)=\phi_{2t}(x)\in\support(\phi).
\end{equation}
For $t\geq1/2$, $\psi_t(x)$ is on the line segment $[\phi_1(x),f(x)]$ which is included in the closed ball $\overline{B}\left(\phi_1(x),\Vert\phi_1(x)-f(x)\Vert\right)$. Since $\Vert\phi_1(x)-f(x)\Vert<\rho$ we get
\begin{equation}\label{equationautomaticdeformationC}
\psi_t(x)\in B(\phi_1(x),\rho)\subset(\xi_{\phi_1})_\rho;
\end{equation}
\item if $x\in\xi_f\setminus\xi_{\phi_1}$ then $\psi_t(x)=x$ for $t\leq1/2$. Using the same argument as above, for $t\geq1/2$ we have
\begin{equation}\label{equationautomaticdeformationD}
\psi_t(x)\in B(f(x),\rho)\subset(\xi_f)_\rho;
\end{equation}
\end{itemize}
Notice that
\begin{equation}
\bigcup_{t\in[0,1]}\left\{x\colon\psi_t(x)\neq x\right\}\subset\support(\phi)\cup\xi_f.
\end{equation}
By~\eqref{equationautomaticdeformationA}, \eqref{equationautomaticdeformationB}, \eqref{equationautomaticdeformationC} and~\eqref{equationautomaticdeformationD} we also get
\begin{equation}
\bigcup_{t\in[0,1]}\psi_t\left(\xi_{\psi_t}\right)\subset\support(\phi)\cup\left(\xi_f\cup\xi_{\phi_1}\right)_\rho,
\end{equation}
which in turn gives
\begin{equation}
\support(\psi)\subset\support{\phi}\cup\left(\xi_f\cup\xi_{\phi_1}\right)_\rho\subsubset U.
\end{equation}
\end{proof}

Let us now define quasiminimal sets, which were introduced by Almgren in~\cite{almgren}. These sets are such that their measure can decrease when deformed, but only in a controlled manner in regards of the size of the points being affected.

\begin{definition}[Quasiminimal sets]\label{definitionquasiminimality}
Let $M\geq 1$ and $\delta>0$. We say that $E$ is a $(M,\delta)$-quasiminimal set over $U$ if $E$ is a relatively closed subset of $U$ with locally finite measure (i.e. $\H^d(E\cap K)<\infty$ for all compact set $K$) such that for all $\delta'$-deformation $(\phi_t)$ over $U$ (with $0<\delta'\leq\delta$) we have
\begin{equation}
\H^d(E\cap\xi_{\phi_1})\leq M\H^d(\phi_1(E\cap\xi_{\phi_1})).
\end{equation}
\end{definition}

In the special case when $M=1$ and $\delta=\diam U$ we say that $E$ is minimal. Now suppose that we are given a function $h\colon]0,\delta]\rightarrow[0,+\infty]$ such that $\lim_{t\rightarrow0}h(t)=0$ and for all $\delta'\leq\delta$, $E$ is $(1+h(\delta'),\delta')$-quasiminimal. We will call such sets --- that look more and more like minimal sets when looked at closely --- almost-minimal sets with gauge function $h$.

To make future statements easier to write, we will also call ``$d$-set'' any $\H^d$-measurable set with locally finite measure, and ``null $d$-set'' any set with null measure.

Since our proofs will involve delicate hair-cutting and measure control tools in varying dimensions, we define the $d$-dimensional core of a set $E$ (which is usually denoted as $E^*$) as follows:
\begin{equation}
\ker^d(E)=\left\{x\in E\colon\forall r>0,\H^d(E\cap B(x,r))>0\right\}.
\end{equation}
We will also use the following notations, for $0\leq l\leq d$:
\begin{equation}
\begin{cases}
\ker_d^d(E)=\ker^d(E)\\
\displaystyle\ker_d^l(E)=\ker^l\left(E\setminus\bigcup_{d\geq l'>l}\ker_d^{l'}(E)\right),
\end{cases}
\end{equation}
and it is easy to check that the $\ker_d^l(E)$ (for $0\leq l\leq d$) are pairwise disjoint and form a partition of $E$. Also, $E\setminus\ker^d(E)$ is a null $d$-set and $\ker^l(E)$ is a relatively closed subset of $E$. Furthermore, if $E$ is $(M,\delta)$-quasiminimal, so is $\ker^d(E)$, and if $E=\ker^d(E)$ we say that $E$ is reduced.

We denote by $\dist_\H$ the Hausdorff distance, which is defined as follows for two non-empty sets $A$ and $B$:
\begin{equation}
\dist_\H(A,B)=\max\left(\sup_{x\in A}\dist(x,B),\sup_{x\in B}\dist(x,A)\right),
\end{equation}
with the conventions $\dist_\H(\emptyset,B)=\dist_\H(A,\emptyset)=\infty$ and $\dist_\H(\emptyset,\emptyset)=0$. For any compact set $K\subset\mathbb{R}^n$ we define the local Hausdorff distance $\dist_K$ over $K$ by:
\begin{equation}
\dist_K(A,B)=\dist_\H(K\cap A,K\cap B).
\end{equation}

We say that a sequence $(E_k)_{k\in\mathbb{N}}$ of sets converges towards $E$ locally on every compact of $U$ if $E$ is a relatively closed subset of $U$ and for all compact subset $K\subset U$:
\begin{equation}
\lim_{k\rightarrow\infty}\dist_K(E_k,E)=0.
\end{equation}
We will denote it by $E_k\Hconverge{U}E$. One can check that this defines an unique limit, and that any domain $U\subset\mathbb{R}^n$ is compact for this convergence in the sense that every sequence has a convergent subsequence.

Finally, in order to prove our main result we need the following theorem, which can be found in~\cite{david1}.
\begin{theorem}\label{theoremdavid}
Suppose that $U\subset\mathbb{R}^n$, $0\leq d<n$, $\delta>0$, $M\geq1$ and $(E_k)_{k\geq 0}$ is a sequence of $(M,\delta)$-quasiminimal sets over $U$ such that $\ker^d(E_k)\Hconverge{U}E$. Then the following holds:
\begin{itemize}
\item $E$ is reduced and $(M,\delta)$-quasiminimal over $U$;
\item for all open subset $W\subset U$,
\begin{equation}
\H^d(E\cap W)\leq\liminf_{k\rightarrow\infty}\H^d(E_k\cap W);
\end{equation}
\item there is $C>0$ such that for all open subset $W\subsubset U$,
\begin{equation}
\H^d(E\cap\overline{W})\geq C^{-1}\limsup_{k\rightarrow\infty}\H^d(E_k\cap\overline{W});
\end{equation}
\item for all $\delta$-deformation $(f_t)_{0\leq t\leq 1}$ over $U$ and $\epsilon>0$, one can build a Lipschitz map $g$ over $U$ such that
\begin{equation}
\Vert f_1-g\Vert_\infty<\epsilon\qquad\text{and}\qquad\xi_g\subsubset\xi_{f_1},
\end{equation}
and for $k$ large enough:
\begin{equation}\label{equationdavidA}
\begin{split}
\H^d(g(E_k\cap\xi_g))&\leq\H^d(f_1(E_k\cap\xi_{f_1}))+\epsilon\\
\H^d(E\cap\xi_{f_1})&\leq\H^d(E_k\cap\xi_g)+\epsilon.
\end{split}
\end{equation}
\end{itemize}
\end{theorem}

In fact, although the first three points gathered in theorem~\ref{theoremdavid} are given as independent statements in~\cite{david1}, the last point is adapted from the proof of the second one (which is called ``Theorem~4.1'' in~\cite{david1}). More precisely, we borrowed equations~[4.93], [4.108] and [4.109] from~\cite{david1}. Starting with $f_1$, a new map $g$ is built such that $\xi_g\subset\xi_{f_1}$ and to which we apply the measure inequalities for $E'_k$. As emphasized by the author, the reason for this process is that we cannot actually use the argument with $f_1$, since it could be injective on $E_k$ and at the same time glue large pieces of $E$ together onto the same image. For this reason we use a small variation of $f_1$ that mimics its behavior and send distinct pieces of $E'_k$ onto the same image when $f_1$ do the same with $E$. Combined with proposal~\ref{proposalautomaticdeformation}, $g$ can also be turned into a $\delta$-deformation over $U$ in order to stay in our topological class $\mathfrak{F}$, and will be used in the proof of theorem~\ref{theoremexistence}.

\section{Orthogonal and radial projections onto polyhedrons}

Our first step is to establish some properties of deformations that will be used later to approximate any given set with polyhedrons. Basically, we will use two kind of deformations: ``magnetic projections'' (see proposal~\ref{proposalmagneticprojection}) that are used to locally flatten a given rectifiable set onto a tangent plane, and radial projections (see definition~\ref{definitionradialprojection}) that send the inside of a polyhedron onto its faces.

\subsection{Fine-tuned Lipschitz extensions}

Before we start building our projections onto polyhedrons, we give some Lipschitz extension lemmas. Although Kirszbraun's theorem (originally in~\cite{kirszbraun}) would be sufficient to get the expected Lipschitz constants, in some cases we also need additional control on the size of the support of the extensions. For this reason we prefer building them explicitly ``by hand''.

\begin{lemma}[Ring-like Lipschitz extensions around a compact]\label{lemmaringextension}
Let $K$ be a nonempty compact set of $\mathbb{R}^n$ and $f$ a $k$-Lipschitz map over $K$ with $k\geq 1$. Suppose that there exists a map $\Pi\colon\mathbb{R}^n\rightarrow K$ such that $f\circ\Pi$ is also $k$-Lipschitz and $\Pi\vert_K=\identity_K$ and put
\begin{equation}K_\rho=\{x\in\mathbb{R}^n\colon\dist(x,K)\leq\rho\}.\end{equation}

Then, for all $\rho>0$ one can find a Lipschitz map $g\colon K_\rho\rightarrow K_\rho$ with constant at most $k+1+\frac{\overline{\dist}(f(K),K)}{\rho}$ such that $g\vert_K=f$, $g\vert_{\partial(K_\rho)}=\identity_{\partial(K_\rho)}$.
\end{lemma}

For instance, if $K$ is convex one can take the convex projector onto $K$ as $\Pi$. Later, we will use this lemma in proposal~\ref{proposalmagneticprojection} when $K$ is the intersection of a cone with a ball to build ``magnetic projections'' that coincide with an affine projector inside $K$ and the identity map outside $K_\rho$.

\begin{proof}
Take $\rho>0$ and suppose that $f$ and $\Pi$ are as above. We define the following map $g$ on $K_\rho$:
\begin{equation}
g(x)=\left(1-\frac{\dist(x,K)}{\rho}\right)f\circ\Pi(x)+\frac{\dist(x,K)}{\rho}x.
\end{equation}
It is easy to check that $g$ is continuous, that $g\vert_K=f$ and $g\vert_{\partial(K_\rho)}=\identity_{\partial(K_\rho)}$. Now all we have to do is to get the required Lipschitz constants for $g$. For that purpose, take $(x,y)\in K_\rho$ and consider the three possible cases:
\begin{itemize}
\item when $(x,y)\in K^2$, since $f$ is $k$-Lipschitz we easily get
\begin{equation}
\Vert g(x)-g(y)\Vert=\Vert f(x)-f(y)\Vert\leq k\Vert x-y\Vert;
\end{equation}
\item when $(x,y)\in(K_\rho\setminus K)^2$, put $x'=\Pi(x)$ and $y'=\Pi(y)$. We now get:
\begin{equation}
\begin{split}
\Vert g(x)-g(y)\Vert&=\tfrac{\Vert\rho f(x')-\rho f(y')+\dist(x,K)(x-f(x'))-\dist(y,K)(y-f(y'))\Vert}{\rho}\\
&\leq\tfrac{\Vert\rho f(x')-\rho f(y')+\dist(x,K)((x-f(x'))-(y-f(y')))\Vert}{\rho}\\
&\qquad+\frac{\Vert(\dist(x,K)-\dist(y,K))(y-f(y'))\Vert}{\rho}\\
&\leq\frac{\rho-\dist(x,K)}{\rho}\Vert f(x')-f(y')\Vert+\frac{\dist(x,K)}{\rho}\Vert x-y\Vert\\
&\qquad+\left\vert\frac{\dist(x,K)-\dist(y,K)}{\rho}\right\vert\Vert y-f(y')\Vert.
\end{split}
\end{equation}
Since we also know that $k\geq 1$, $\dist(x,K)\leq\rho$, $\Vert f(x')-f(y')\Vert\leq k\Vert x-y\Vert$ and $\Vert y-f(y')\Vert\leq\rho+\dist_\H(f(K),K)$ we finally get:
\begin{equation}
\begin{split}
\Vert g(x)-g(y)\Vert&\leq\left(\frac{k\rho-(k-1)\dist(x,K)}{\rho}+\frac{\rho+\dist_\H(f(K),K)}{\rho}\right)\Vert x-y\Vert\\
&\leq\left(k+1+\frac{\dist_\H(f(K),K)}{\rho}\right)\Vert x-y\Vert;
\end{split}
\end{equation}
\item when $x\in K$ and $y\in K_\rho\setminus K$, we put as above $y'=\Pi(y)$ and get:
\begin{equation}
\begin{split}
\Vert g(x)-g(y)\Vert&=\frac{\Vert\rho f(x)-\dist(y,K)y-(\rho-\dist(y,K))f(y')\Vert}{\rho}\\
&=\frac{\Vert\dist(y,K)(f(x)-y)-(\rho-\dist(y,K))(f(x)-f(y'))\Vert}{\rho}\\
&\leq\frac{\dist(y,K)}{\rho}\Vert f(x)-y\Vert+\frac{\rho-\dist(y,K)}{\rho}\Vert f\circ\Pi(x)-f\circ\Pi(y)\Vert\\
&\leq\frac{\Vert x-y\Vert}{\rho}(\rho+\dist_\H(f(K),K))+k\Vert x-y\Vert\\
&\leq\left(k+1+\frac{\dist_\H(f(K),K)}{\rho}\right)\Vert x-y\Vert.
\end{split}
\end{equation}
\end{itemize}
In all cases, we have shown that $g$ is $k'$-Lipschitz with $k'=1+k+\frac{\dist_\H(f(K),K)}{\rho}$.
\end{proof}

Conversely, the following lemma is used to build a Lipschitz extension inside a ball that have been subtracted from a compact.

\begin{lemma}[Lipschitz extension inside a ball]\label{lemmaholeextension}
Suppose that $K$ is a star compact with respect to $x$ that contains an open ball $B$ centered at $x$ with radius $r$ and put $K'=K\setminus B$. For $\rho>0$ we denote by $\rho B$ the ball centered at $x$ with radius $\rho r$.

For all $k$-Lipschitz map $f\colon K'\rightarrow K'$ and $\rho\in]0,1[$ one can build a $k'$-Lipschitz map $g\colon K\rightarrow K$ such that $g\vert_{K'}=f\vert_{K'}$, $g\vert_{\rho B}=\identity_{\rho B}$ and $k'$ depends only on $\rho$, $\diam(K)$ and $r$.
\end{lemma}

\begin{proof}
For all $y\in B\setminus\rho B$ there is only one point in $[x,y)\cap\partial B$ which we call $\Pi(y)$. We can notice already that $\Pi$ is $\frac{1}{\rho}$-Lipschitz. When $y\in B\setminus\rho B$ we put
\begin{equation}
u(y)=\frac{\Vert\Pi(y)-y\Vert}{r(1-\rho)}\in[0,1[
\end{equation}
and we define $h\colon K\setminus\rho B\rightarrow K\setminus\rho B$ as
\begin{equation}
g(y)=\begin{cases}
f(y)&\text{if $x\in K\setminus B$}\\
u(y)y+(1-u(y))f\circ\Pi(y)&\text{if $y\in B\setminus\rho B$.}
\end{cases}
\end{equation}
It is easy to check that $h$ is continuous, and that $h\vert_{\partial\rho B}=\identity_{\partial\rho B}$. Now suppose that $(y,z)\in(K\setminus\rho B)^2$ and consider the three following cases:
\begin{itemize}
\item if $(y,z)\in(K\setminus B)^2$ then
\begin{equation}
\Vert h(y)-h(z)\Vert=\Vert f(y)-f(z)\Vert\leq k\Vert y-z\Vert;
\end{equation}
\item if $(y,z)\in B^2$ we get
\begin{equation}
\begin{split}
\Vert h(y)-h(z)\Vert&=\Vert f\circ\Pi(y)-f(z)+u(y)(y-f\circ\Pi(y))-u(z)(z-f\circ\Pi(z))\Vert\\
&\leq\Vert f\circ\Pi(y)-f\circ\Pi(z)\Vert+\Vert u(y)(y-f\circ\Pi(y)-z+f\circ\Pi(z))\Vert\\
&\qquad+\Vert (u(y)-u(z))(z-f\circ\Pi(z))\Vert\\
&\leq\frac{2k}{\rho}\Vert y-z\Vert+\Vert y-z\Vert+\frac{k}{\rho}\Vert y-z\Vert+\diam(K)\vert u(y)-u(z)\vert\\
&=\left(2\frac{k}{\rho}+1\right)\Vert y-z\Vert+\diam(K)\left\vert\frac{\Vert y-\Pi(y)\Vert-\Vert z-\Pi(z)\Vert}{r(1-\rho)}\right\vert\\
&\leq\left(\frac{2k}{\rho}+1\right)\Vert y-z\Vert+\diam(K)\frac{\Vert y-z-\Pi(y)+\Pi(z)\Vert}{r(1-\rho)}\\
&\leq\left(\frac{2k}{\rho}+1\right)\Vert y-z\Vert+\diam(K)\frac{\Vert y-z\Vert+\Vert\Pi(y)-\Pi(z)\Vert}{r(1-\rho)}\\
&\leq\left(\frac{2k}{\rho}+1\right)\Vert y-z\Vert+\diam(K)\frac{\Vert y-z\Vert+\rho^{-1}\Vert y-z\Vert}{r(1-\rho)}\\
&=\left(1+\frac{2k}{\rho}+\frac{2\diam(K)}{\rho(1-\rho)r}\right)\Vert y-z\Vert;
\end{split}
\end{equation}
\item finally, if $y\in K\setminus B$ and $z\in B$ we have
\begin{equation}
\begin{split}
\Vert h(y)-h(z)\Vert&=\Vert f(y)-u(z)z-(1-u(z))f\circ\Pi(z)\Vert\\
&\leq u(z)\Vert y-z\Pi(z)\Vert+(1-u(z))\Vert f(y)-f\circ\Pi(z)\Vert\\
&\leq\Vert y-z\Vert+k\Vert y-\Pi(z)\Vert\\
&\leq\left(1+\frac{k}{\rho}\right)\Vert y-z\Vert.
\end{split}
\end{equation}
\end{itemize}

We have just shown that $h$ is Lipschitz. Now all we have to do is to apply lemma~\ref{lemmaringextension} to extend $h$ inside $\rho B$ and lemma~\ref{lemmaholeextension} will be proven.
\end{proof}

The last extension theorem we provide is used to extend a Lipschitz map defined on the subfaces of a complex to the whole Euclidean space, while keeping its support as small as prescribed.

\begin{lemma}[Lipschitz extension around a complex]\label{lemmacomplexextension}
Let $k\in\left\{0,\ldots,n\right\}$, $S$ a $k$-dimensional complex and $U$ an open bounded set such that $\U(S)\subset U$. Suppose that for each $\delta\in S$ we are being given a Lipschitz map $\phi_\delta\colon\delta\rightarrow\delta$ such that $\phi\vert_{\partial\delta}=\identity_{\partial\delta}$.

Then we can find a Lipschitz map $\phi\colon\mathbb{R}^n\rightarrow\mathbb{R}^n$ such that:
\begin{equation}
\forall\delta\in S\colon\phi\vert_\delta=\phi_\delta\quad\text{and}\quad\phi\vert_{\mathbb{R}^n\setminus U}=\identity_{\mathbb{R}^n\setminus U}.
\end{equation}
\end{lemma}

Notice that we do not really care about the Lipschitz constant of the final map, although we could give an estimate based upon the largest one of those of the $\phi_\delta$ and the rotondity of $S$.

\begin{proof}
All we really have to do is to prove that the map $\psi_S$ defined on the closed set $F=\U(S)\cup(\mathbb{R}^n\setminus U)$ as
\begin{equation}
\psi_S(x)=\begin{cases}
x&\text{if $x\notin U$}\\
\phi_\delta(x)&\text{if $x\in\delta\in S$}
\end{cases}
\end{equation}
is Lipschitz and to apply Kirszbraun theorem to it.

To begin with, one can check that the definition of $\psi_S$ is consistent. Firstly, notice that any polyhedron inside $S$ is disjoint of $\mathbb{R}^n\setminus U$. Additionally, if one can find $x\in\delta_1\cap\delta_2$ such that $(\delta_1,\delta_2)\in S^2$ and $\delta_1\neq\delta_2$ then --- by definition~\ref{definitioncomplex} of a complex --- $\delta_1\cap\delta_2\subset\partial\delta_1\cup\partial\delta_2$ and we have $\phi_{\delta_1}(x)=\phi_{\delta_2}(x)=x$.

We will now prove that $\psi_S$ is Lipschitz by induction over the number of polyhedrons in $S$. In what follows, for each $\delta\in S$ we suppose that $\phi_\delta$ is $k_\delta$-Lipschitz.

If $S$ is made of only one polyhedron $\delta$, put
\begin{equation}
a=\min_{x\in\delta}\left(\dist(x,\mathbb{R}^n\setminus U)\right)>0\quad\text{and}\quad b=\max_{x\in\delta}\left(\dist(x,\mathbb{R}^n\setminus U)\right)<\infty.
\end{equation}
We already know that $\psi_S$ is $1$-Lipschitz over $\mathbb{R}^n\setminus U$ and $k_\delta$-Lipschitz over $\delta$. If $x\in\mathbb{R}^n\setminus U$ and $y\in\delta$ we have $\psi_S(y)=\phi_\delta(y)\in\delta$ and we get
\begin{equation}
\Vert\psi_S(x)-\psi_S(y)\Vert=\Vert x-\psi_S(y)\Vert\leq b\leq\frac{b}{a}\Vert x-y\Vert.
\end{equation}

Now suppose that $S=S'\sqcup\left\{\delta\right\}$ (with $S\neq\emptyset$) and that $\psi_{S'}$ is $k$-lipschitz. Let $x\in\U(S')$ and $y\in\delta$ and consider the line segment $[\psi_{S'}(x),\phi_\delta(y)]$; since $\phi_\delta(y)\in\delta$ and $\psi_{S'}(x)\notin\partial\delta$ we know that this line segment meets $\partial\delta$ at at least one point $y'$.

Let us verify that one can always find a subface $F\in\F(\delta)$ such that $y'\in F$. First, since $\U(\F_{k-1}(\delta))=\partial\delta$ then one can find a face $F_1$ such that $y'\in F_1$. If $y'\in\open{F_1}$ we have finished. Otherwise $y'\in\partial F_1$ and again, one can find $F_2\in\F_{k-2}(F_1)$ such that $y'\in F_2$. By iterating this argument while $y'\notin\open{F_{i-1}}$ one can find a subface $F_i\in\F_{k-i}(\delta)$ such that $y'\in\open{F_i}$ or $y'\in\partial F_i$. Since subfaces of dimension zero are singletons --- equal to their interior, following our conventions --- this building process will stop eventually with at most $i=k$ (in such case $y'$ is a vertex of $\delta$) and in all cases we can find $F\in\F(\delta)$ such that $y'\in\open{F}$.

Denote by $S''$ the subset of $S$ made of the polyhedrons that do not intersect $F$ and by $\delta'$ a polyhedron in $S'$ such that $\psi_{S'}(x)\in\delta'$. There are three possible cases:
\begin{itemize}
\item if $F$ is a common subface of both $\delta$ and at least one polyhedron of $S'$ we have $\psi_S(y')=\psi_{S'}(y')=\phi_\delta(y')=y'$ and we get
\begin{equation}
\begin{split}
\Vert\psi_S(x)-\psi_S(y)\Vert&=\Vert\psi_{S'}(x)-\phi_\delta(y)\Vert\\
&=\Vert\psi_{S'}(x)-\psi_{S'}(y')\Vert+\Vert\phi_\delta(y')-\phi_\delta(y)\Vert\\
&\leq(k+k_\delta)(\Vert x-y'\Vert+\Vert y'-y\Vert)\\
&=(k+k_\delta)\Vert x-y\Vert;
\end{split}
\end{equation}
\item if $\delta'\in S''$ we put
\begin{equation}
a(F)=\min_{x\in\U(S'')}\dist(x,F)>0\qquad\text{and}\qquad b(F)=\max_{x\in\U(S'')}\dist(x,F)\in]0,1]
\end{equation}
and we get
\begin{equation}
\begin{split}
\Vert\psi_S(x)-\psi_S(y)\Vert&=\Vert\psi_{S'}(x)-\phi_\delta(y)\Vert\\
&=\Vert\psi_{S'}(x)-y'\Vert+\Vert\phi_\delta(y')-\phi_\delta(y)\Vert\\
&\leq b(F)+k_\delta\Vert y'-y\Vert\\
&\leq\left(\frac{b(F)}{a(F)}+k_\delta\right)\Vert x-y\Vert;
\end{split}
\end{equation}
\item lastly, if $\delta'\notin S''$ we put $H=\affine(F)$ and $G=\F_0(F)\cap\F_0(\delta')$ (i.e. $G$ is the set of vertices common to both $F$ and $\delta'$). We consider the minimal ratio of the distance to $H$ by the distance to $G$ of vertices of $\delta'$ that are not in $G$:
\begin{equation}
a(F)=\min\left\{\frac{\dist(c,H)}{\dist(c,G)}\colon c\in\F_0(\delta')\text{ et }c\notin G\right\}>0.
\end{equation}
By a convexity argument it is easy to check that for all $t\in\delta'$ --- and in particular for $t=\psi_{S'}(x)$ --- we have
\begin{equation}
\dist(t,H)\geq a(F)\dist(t,G).
\end{equation}
By denoting by $c$ a vertex common to both $F$ and $\delta'$ whose distance to $\psi_{S'}(x)$ is minimal we also get:
\begin{equation}\label{equationcomplexextensionA}
\begin{split}
\Vert\psi_S(x)-\psi_S(y)\Vert&=\Vert\psi_{S'}(x)-\phi_\delta(y)\Vert\\
&=\Vert\psi_{S'}(x)-\psi_{S'}(c)\Vert+\Vert\phi_\delta(c)-\phi_\delta(y)\Vert\\
&\leq(k+k_\delta)\left(\Vert x-c\Vert+\Vert c-y\Vert\right).
\end{split}
\end{equation}
Consider triangle $xcy$ and denote by $\hat x$, $\hat c$ and $\hat y$ the non-oriented angles respectively at vertices $x$, $c$ and $y$. A simple planar geometry identity gives us that
\begin{equation}
\frac{\Vert x-y\Vert}{\sin\hat c}=\frac{\Vert x-c\Vert}{\sin\hat y}=\frac{\Vert c-y\Vert}{\sin\hat x}.
\end{equation}
To conclude, notice that the sinus of the non-oriented angle between the lines $(x,c)$ and $(y,c)$ is between $a(F)$ and $1$. It follows that
\begin{equation}
\begin{split}
\Vert x-c\Vert+\Vert c-y\Vert&=\frac{\sin\hat y}{\sin\hat c}\Vert x-y\Vert+\frac{\sin\hat x}{\sin\hat c}\Vert x-y\Vert\\
&\leq\frac{2}{a(F)}\Vert x-y\Vert
\end{split}
\end{equation}
and by inequality~\eqref{equationcomplexextensionA} we finally get
\begin{equation}
\Vert\psi_S(x)-\psi_S(y)\Vert\leq\frac{2(k+k_\delta)}{a(F)}\Vert x-y\Vert.
\end{equation}
\end{itemize}
In all three cases we could give a constant $c(F,\delta')$ such that $\Vert\psi_S(x)-\psi_S(y)\Vert\leq c(F,\delta')\Vert x-y\Vert$. By taking the maximum of $c(F,\delta')$ for all possible subfaces $F$ of $\delta$ and polyhedrons $\delta'\in S'$ --- which are in finite number --- we get a global constant $c$.

This achieve proving that $\psi_S\vert_{\U(S)}$ is $c$-Lipschitz. To prove that it is also Lipschitz over $\U(S)\cup(\mathbb{R}^n\setminus U)$ one can easily adapt the argument at the beginning of the proof by induction when $S$ contains only one polyhedron.

By induction over the number of polyhedrons in $S$, this achieves proving that $\psi_S$ is Lipschitz over $\U(S)\cup(\mathbb{R}^n\setminus U)$. By using Kirszbraun's theorem it is also possible to build a Lipschitz extension of $\psi_S$ over the whole space $\mathbb{R}^n$ that meets the announced requirements.
\end{proof}

\subsection{Measure-optimal projections}

We now introduce the two basic tools that will allow us later to build a deformation of a given rectifiable set onto a polyhedric mesh without increasing its measure too much. We begin with Lipschitz maps with compact support used to locally ``flatten'' the set onto its approximate tangent planes.

\begin{proposal}[Magnetic projection]\label{proposalmagneticprojection}
Let $K$ be a nonempty compact set of $\mathbb{R}^n$ and $H$ an affine subspace. Let $p$ be the orthogonal projector on $H$, $\vec H$ the linear subspace $H-p(0)$ and suppose that $p(K)\subset K$, $H\cap K$ is convex and for all $x\in H\cap K$, the compact set $K(x)=K\cap(x+\vec H^\bot)$ is convex.

Then for all $\rho>0$, one can find a so-called ``$\rho$-magnetic projection onto $H$ inside $K$'' map $\Pi_{H,\rho,K}\colon\mathbb{R}^n\rightarrow H$ verifying the following properties:
\begin{itemize}
\item$\Pi_{H,\rho,K}(K_\rho)\subset K_\rho$, where $K_\rho=\{x\in\mathbb{R}^n\colon\dist(x,K)\leq\rho\}$;
\item$\Pi_{H,\rho,K}\vert_{K}=p\vert_K$ and $\Pi_{H,\rho,K}\vert_{\mathbb{R}^n\setminus K_\rho}=\identity_{\mathbb{R}^n\setminus K_\rho}$;
\item$\Pi_{H,\rho,K}$ is Lipschitz with constant at most $2+\frac{\dist_\H(H\cap K,K)}{\rho}$.
\end{itemize}
\end{proposal}

\begin{proof}
Suppose that $A$ is a nonempty convex compact set of $\mathbb{R}^n$. By compacity, for all $x\in\mathbb{R}^n$ one can find $y\in A$ such that $\Vert x-y\Vert=\dist(x,A)$ and by convexity, $y$ is unique; we call it the projection of $x$ onto the convex set $A$ and denote it by $\pi_A(x)$. Let us rapidly verify that $\pi_A$ is $1$-Lipschitz. When $A$ is a singleton or a line segment it is very easy to check. Otherwise, take $(x,y)\in\mathbb{R}^n$ and put $u=\pi_A(x)$ and $v=\pi_A(y)$. Since $[u,v]\subset A$ we have $\Vert\pi_A(x)-\pi_A(y)\Vert\leq\Vert\pi_{[u,v]}(x)-\pi_{[u,v]}(y)\Vert$ and the Lipschitz constant of $\pi_A$ follows immediately from the one of $\pi_{[u,v]}$.

Now, fix $\rho>0$ and for $x\in\mathbb{R}^n$ consider its projection $\pi_{H\cap K}(x)$ onto the nonempty convex set $H\cap K$. Since $\pi_{H\cap K}(x)\in H\cap K$ then the compact set $K(\pi_{H\cap K}(x))=K\cap(\pi_{H\cap K}(x)+\vec H^\bot)$ is nonempty and by hypothesis, convex. We will denote by $\Pi(x)$ the projection of $x$ onto this new convex:
\begin{equation}
\forall x\in\mathbb{R}^n\colon\Pi(x)=\pi_{K(\pi_{H\cap K}(x))}(x).
\end{equation}
By construction $x\in K(\pi_{H\cap K}(x))(x)$, therefore
\begin{equation}
p\circ\Pi(x)\in p(K(\pi_{H\cap K}(x)))\subset H\cap(\pi_{H\cap K}(x)+\vec H^\bot)=\{\pi_{H\cap K}(x)\}.
\end{equation}
It follows that $p\circ\Pi=\pi_{H\cap K}$ is $1$-Lipschitz, and it is easy to check that $\Pi\vert_K=\identity_K$. 
\end{proof}

In what follows we suppose that $0\leq d<n$ and consider a $d$-set $E$. For $x\in\mathbb{R}^n$, we define the lower and upper radial $d$-dimensional densities of $E$ at $x$ respectively by putting
\begin{equation}
\underline{\nu}^d_E(x)=\liminf_{r\rightarrow 0}\frac{\H^d(E\cap B(x,r))}{c_dr^d}\qquad\overline{\nu}^d_E(x)=\limsup_{r\rightarrow 0}\frac{\H^d(E\cap B(x,r))}{c_dr^d}
\end{equation}
where $c_d$ stands for the measure of the $d$-dimensional unit ball.

Also, we say that $H$ is an approximate tangent plane for $E$ at $x$ if $H$ is a $d$-plane containing $x$, $\underline{\nu}^d_E(x)>0$ and
\begin{equation}\label{equationA}
\forall u>0\colon\limsup_{r\rightarrow 0}\frac{\H^d(E\setminus\C(x,r,u))}{r^d}=0
\end{equation}
where $\C(x,r,u)$ stands for the following intersection between a cone that ``follows'' $H$ and a closed ball centered at $x$:
\begin{equation}\label{equationcone}
\C(x,r,u)=\left\{y\in\overline{B}(x,r)\colon\dist(y,H)\leq u\Vert x-y\Vert\right\}.
\end{equation}

If $E$ has such an approximate tangent plane at $\H^d$ almost every point we say that $E$ is $d$-rectifiable. Conversely, if $E$ has no approximate tangent plane at almost every point we say that it is $d$-irregular, which is the same as saying that any rectifiable set intersects $E$ only on a null set. It is well-known (again, see for instance Mattila's book~\cite{mattila}) that $E$ is rectifiable if and only if $\underline{\nu}^d_E$ and $\overline{\nu}^d_E$ are equal to the characteristic set function of $E$, $\H^d$ almost everywhere. Conversely, $E$ is irregular if and only if $\overline{\nu}^d_E$ is less than $1$ almost everywhere. As a consequence, any $d$-set $E$ can be written as
\begin{equation}
E=E_R\sqcup E_I
\end{equation}
with $E_R$ rectifiable and $E_I$ irregular. We will refer to $E_R$ and $E_I$ respectively as the rectifiable and irregular parts of $E$ --- which are defined up to a null set.

The next lemma makes use of the previous proposal introducing magnetic projections in the following context. At almost every point of $E$ where there is an approximate tangent plane, one can find a ball such that the magnetic projection onto the tangent plane inside a small neighborhood of the ball does not increase the measure of the set too much.

\begin{lemma}[Magnetic projection inside a high density cone]\label{lemmamagneticprojection}
Let $E$ be a $d$-set. For all $\epsilon>0$ and at $\H^d$ almost every point $x$ of the rectifiable part of $E$ one can find $r_{\max}>0$, $\rho\in]0,1[$, $u>0$ and an approximate tangent plane $H$ at $x$ such that for all $r\in]0,r_{\max}[$:
\begin{equation}\label{equationmagneticprojection}
\H^d(\Pi_{H,\rho r,\C(x,r,u)}(E\cap B(x,r+r\rho)\setminus\C(x,r,u)))\leq\epsilon\H^d(E\cap B(x,r+r\rho)).
\end{equation}
\end{lemma}

\begin{proof}
First, notice that the above $\C(x,r,u)$ is suitable to be used as ``$K$'' in proposal~\ref{proposalmagneticprojection}. Fix $\epsilon'>0$, $u>0$ and $\rho\in]0,1[$.

Suppose that the lower and upper radial densities of $E$ at $x$ are equal to $1$. We can find $r_1>0$ such that for all $t\leq r_1$:
\begin{equation}\label{equationmagneticprojectionA}
c_d(2t)^d(1+\epsilon')^{-1}\leq\H^d(E\cap B(x,t))\leq c_d(2t)^d(1+\epsilon').
\end{equation}
By taking $t=r$ and $t=r+r\rho$ in~\eqref{equationmagneticprojectionA} it follows that for all $r<\frac{r_1}{2}$:
\begin{multline}\label{equationmagneticprojectionB}
\H^d(E\cap B(x,r+r\rho)\setminus B(x,r))\\
\begin{aligned}
&\leq 2^dc_d(1+\epsilon')(r+r\rho)^d-2^dc_d(1+\epsilon')^{-1}r^d\\
&\leq 2^dc_d(r+r\rho)^d(1+\epsilon')^{-1}\left((1+\epsilon')^2-\frac{r^d}{(r+r\rho)^d}\right)\\
&\leq\left((1+\epsilon')^2-\frac{1}{(1+\rho)^d}\right)\H^d(E\cap B(x,r+r\rho))
\end{aligned}
\end{multline}
Suppose that $\rho$ is small enough so $(1+\rho)^d<\frac{1}{(1-\epsilon')^2}$. By replacing in~\eqref{equationmagneticprojectionB} we obtain
\begin{multline}
\H^d(E\cap B(x,r+r\rho)\setminus B(x,r))\\
\begin{aligned}
&\leq((1+\epsilon')^2-(1-\epsilon')^2)\H^d(E\cap B(x,r+r\rho))\\
&=2\epsilon'\H^d(E\cap B(x,r+r\rho)).
\end{aligned}
\end{multline}

Also, suppose that $H$ is an approximate tangent plane at $x$. By~\eqref{equationA} we can find $r_2>0$ such that for all $r<r_2$:
\begin{equation}
\begin{split}
\H^d(E\cap B(x,r)\setminus\C(x,r,u))&\leq\epsilon'r^dc_d\\
&\leq\epsilon'(1+\epsilon')\H^d(E\cap B(x,r))\\
&\leq\epsilon'(1+\epsilon')\H^d(E\cap B(x,r+r\rho)).
\end{split}
\end{equation}

On the other hand, we can write that
\begin{multline}
E\cap B(x,r+r\rho)\setminus\C(x,r,u)
\\=\left(E\cap B(x,r+r\rho)\setminus B(x,r)\right)\sqcup\left(E\cap B(x,r)\setminus\C(x,r,u)\right)
\end{multline}
and since $\Pi_{H,\rho r,\C(x,r,u)}$ is $2+\frac{u}{\rho}$-Lipschitz by proposal~\ref{proposalmagneticprojection} we get
\begin{multline}
\H^d(\Pi_{H,\rho r,\C(x,r,u)}(E\cap B(x,r+r\rho)\setminus\C(x,r,u)))\\
\leq\left(2+\frac{u}{\rho}\right)^d(2\epsilon'+\epsilon'(1+\epsilon'))\H^d(E\cap B(x,r+r\rho)).
\end{multline}

To conclude, all we have to do is taking $u>0$ small enough such that $\left(2+\frac{u}{\rho}\right)^d<2^d+\epsilon'$ and we get
\begin{multline}
\H^d(\Pi_{H,\rho r,\C(x,r,u)}(E\cap B(x,r+r\rho)\setminus\C(x,r,u)))\leq\\
\epsilon'(2^d+\epsilon')(3+\epsilon')\H^d(E\cap B(x,r+r\rho)).
\end{multline}
Put $r_{\max}=\min\left(\frac{r_1}{2},r_2\right)$ and recall that at $\H^d$ almost every point of the rectifiable part of $E$, the radial densities are equal to $1$ and $E$ has an approximate tangent plane. Being given $\epsilon>0$, by taking $\epsilon'$ small enough to get $\epsilon'(2^d+\epsilon')(3+\epsilon')<\epsilon$ this achieves proving the lemma.
\end{proof}

Following definition~\ref{definitionpolyhedron}, our polyhedrons are nonempty, convex and compact. Inside the generated affine subspace, any half-line starting in the interior of a polyhedron will intersect its boundary at one unique point, which legitimates the following definition.

\begin{definition}[Radial projection]\label{definitionradialprojection}
Suppose that $\delta$ is a $k$-dimensional polyhedron (with $1\leq k\leq n$) and that $x\in\open{\delta}$. We define the radial projection $\Pi_{\delta,x}$ onto the faces of $\delta$ by
\begin{equation}
\Pi_{\delta,x}\colon\begin{cases}
\delta\setminus\{x\}\rightarrow\partial\delta\\
y\mapsto z\in[x,y)\cap\partial\delta.
\end{cases}
\end{equation}
\end{definition}

It is easy to check that $\Pi_{\delta,x}\vert_{\partial\delta}=\identity_{\partial\delta}$ and that $\Pi_{\delta,x}\vert_{\delta\setminus U}$ is Lipschitz for all open set $U$ containing $x$. The following lemma will allow us to control the measure increase of the radial projection of a given $d$-set with constants depending on the polyhedron's rotondity.

\begin{lemma}[Optimal radial projection]\label{lemmaoptimalradialprojection}
Suppose that $0\leq d<k\leq n$. There exists a constant $K>0$ depending only on $d$, $k$ and $n$ such that for all $k$-dimensional polyhedron $\delta$ and closed $d$-set $E$ contained in $\delta$, one can find $X\subset\open{\delta}$ with positive $\H^k$-measure such that:
\begin{equation}
\forall x\in X\colon\H^d(\Pi_{\delta,x}(E))\leq KR(\delta)^{-2d}\H^d(E).
\end{equation}
\end{lemma}

The proof will use a mean value argument and Fubini's theorem. Although it would have been more convenient to use the Jacobian determinant of $\phi_{\delta,x}$ and a change of variables when computing the mean value of $\H^d(\Pi_{\delta,x}(E))$, this approach would have required additional assumptions on the regularity of $E$. For this reason we will slice $\delta$ in thin pieces parallel to its faces and approximate the integral by summing the measure in each piece.

\begin{proof}
Suppose that $B$ is an inscribed ball inside $\delta$, put $B'=\frac{1}{2}B$ and $H=\affine(\delta)$. For $z\in\partial\delta$ we denote by $n(z)$ an unit vector parallel to $H$ which is $\H^{k-1}$ almost everywhere normal to $\partial\delta$ at $z$, and by $\partial^*\delta$ the subset of $\partial\delta$ where $n(z)$ is affectively normal to $\partial\delta$. We also define
\begin{equation}
\tau_x(z)=\frac{\Vert z-x\Vert}{\vert\left<n(z),z-x\right>\vert}\qquad\text{and}\qquad A=\sup_{x\in B',z\in\partial^*\delta}\tau_x(z),
\end{equation}
where $\left<\cdot,\cdot\right>$ stands for the usual Euclidean dot product in $\mathbb{R}^n$.

For all $z\in\partial^*\delta$ one can find a face $F\in\F_{k-1}(\delta)$ containing $z$. Put $H'=\affine(F)$: by construction, $n(z)$ is normal to $H'$ and
\begin{equation}\label{equationoptimalradialprojectionA}
\tau_x(z)=\frac{\Vert z-x\Vert}{\dist(x,H')}.
\end{equation}
Since we supposed that $x\in B'$ and since $\delta$ is contained in a ball with the same center as $B'$ with radius $2\overline{R}(\delta)$ we also get:
\begin{equation}\label{equationoptimalradialprojectionB}
\dist(x,H)\geq\frac{\underline{R}(\delta)}{2}\qquad\text{and}\qquad\dist(x,z)\leq 2\overline{R}(\delta).
\end{equation}
Using~\eqref{equationoptimalradialprojectionA} and~\eqref{equationoptimalradialprojectionB} we deduce that
\begin{equation}
A\leq4\frac{\overline{R}(\delta)}{\underline{R}(\delta)}=\frac{4}{R(\delta)}.
\end{equation}

Fix an integer $p>0$ and a point $x\in B'$. Consider the set $\left\{F_1,\ldots,F_m\right\}=\F_{k-1}(\delta)$ of faces of $\delta$ and put $H_i=\affine(F_i)$ for $1\leq i\leq m$ (each $H_i$ is an affine hyperplane of $H$ of dimension $k-1$). For $r>0$, denote by $h_r$ the homothecy centered at $x$ with dilatation factor $r$ and consider the following sets:
\begin{gather}
\C^i_l(x)=\bigcup_{\frac{l}{p}<r\leq\frac{l+1}{p}}h_r(F_i),\\
\C^i(x)=\bigcup_{0\leq r\leq 1}h_r(F_i)=\bigcup_{0\leq l<p}\C^i_l(x),\\
\delta_l=\bigcup_{\frac{l}{p}<r\leq\frac{l+1}{p}}h_r(\delta)=\bigcup_{1\leq i\leq m}\C^i_l(x).
\end{gather}
Since $x\in B'\subset\open{\delta}$ and by convexity we have the following identities:
\begin{equation}\label{equationoptimalradialprojectionC}
\delta\setminus\{x\}=\bigcup_{i,l}\C^i_l(x)=\bigcup_i\C^i(x)=\bigcup_l\delta_l(x).
\end{equation}
Furthermore, the sets $\delta_l(x)$ are disjoint for $0\leq l<p$.

Suppose that $l>0$ and notice that the restriction of $\Pi_{\delta,x}$ to $\C^i(x)$ is the radial projection centered at $x$ on $H_i$. Then it is Lipschitz with constant at most
\begin{equation}\label{equationoptimalradialprojectionD}
\frac{p}{l}\sup_{z\in F_i\cap\partial^*\delta}\tau_x(z)\leq\frac{pA}{l}.
\end{equation}

Following~\eqref{equationoptimalradialprojectionC}, the measure of the radial projection of $E$ can be rewritten as
\begin{equation}
\begin{split}
\H^d(\Pi_{\delta,x}(E))&=\sum_{0\leq l<p}\H^d(\Pi_{\delta,x}(E\cap\delta_l))\\
&=\H^d(\Pi_{\delta,x}(E\cap\delta_0))+\sum_{1<l<p}\H^d(\Pi_{\delta,x}(E\cap\delta_l)).
\end{split}
\end{equation}
Since $x\in B'\setminus E$ and we supposed that $E$ is closed then for $p$ large enough we have $E\cap\delta_0=\emptyset$ and using~\eqref{equationoptimalradialprojectionD} we get:
\begin{equation}\label{equationoptimalradialprojectionE}
\H^d(\Pi_{\delta,x}(E))=\sum_{1<l<p}\H^d(\Pi_{\delta,x}(E\cap\delta_l))\leq A^d\sum_{1<l<p}\left(\frac{p}{l}\right)^d\H^d(E\cap\delta_l).
\end{equation}
When $y\in\delta_l$ we have $\Vert y-x\Vert<\frac{l+1}{p}\overline{R}(\delta)<\frac{2l}{p}\overline{R}(\delta)$. It follows that
\begin{equation}
\H^d(E\cap\delta_l)=\int_{y\in E\cap\delta_l}d\H^d(y)\leq\left(\frac{2l}{p}\overline{R}(\delta)\right)^d\int_{y\in E\cap\delta_l}\frac{d\H^d(y)}{\Vert y-x\Vert^d}
\end{equation}
and by replacing in~\eqref{equationoptimalradialprojectionE}:
\begin{multline}\label{equationoptimalradialprojectionF}
\H^d(\Pi_{\delta,x}(E))\leq(2A\overline{R}(\delta))^d\sum_{1<l<p}\int_{y\in E\cap\delta_l}\frac{d\H^d(y)}{\Vert y-x\Vert^d}\\
=(2A\overline{R}(\delta))^d\int_{y\in E}\frac{d\H^d(y)}{\Vert y-x\Vert^d}.
\end{multline}

Let us now compute the mean value of $\H^d(\Pi_{\delta,x}(E))$ when $x\in B'\setminus E$. Using~\eqref{equationoptimalradialprojectionF} we already have
\begin{equation}\label{equationoptimalradialprojectionG}
\int_{x\in B'\setminus E}\H^d(\Pi_{\delta,x}(E))d\H^k(x)\leq(2A\overline{R}(\delta))^d\int_{x\in B'\setminus E}\int_{y\in E}\frac{d\H^d(y)d\H^k(x)}{\Vert y-x\Vert^d}
\end{equation}
and since $B'$ is a $k$-dimensional ball with radius $\frac{\underline{R}(\delta)}{2}$ and $1\leq d\leq k$ we also get:
\begin{equation}
\int_{x\in B'\setminus E}\frac{d\H^k(x)}{\Vert y-x\Vert^d}=\int_{x\in B'}\frac{d\H^k(x)}{\Vert y-x\Vert^d}=C\underline{R}(\delta)^{k-d}<\infty,
\end{equation}
where $C$ is a positive constant depending only on $d$ and $k$. Also, we supposed that $E$ is a $d$-set included in $\delta$ and since $\delta$ is compact we can write that
\begin{equation}
\int_{y\in E}C\underline{R}(\delta)^{k-d}d\H^d(y)=C\underline{R}(\delta)^{k-d}\H^d(E)<\infty
\end{equation}
which allows using Fubini's theorem in~\eqref{equationoptimalradialprojectionG}:
\begin{equation}\label{equationoptimalradialprojectionH}
\int_{x\in B'\setminus E}\H^d(\Pi_{\delta,x}(E))d\H^k(x)\leq (2A)^dC\underline{R}(\delta)^{k-d}\overline{R}(\delta)^d\H^d(E).
\end{equation}

On the other hand, one can find $D>0$ depending only on $k$ such that
\begin{equation}
\H^k(B'\setminus E)=\H^k(B')=D\underline{R}(\delta)^k.
\end{equation}
Along with~\eqref{equationoptimalradialprojectionH} this proves that it is possible to find a subset $X\subset B'\setminus E$ of positive measure such that, for instance:
\begin{equation}
\begin{split}
\forall x\in X\colon\H^d(\Pi_{\delta,x}(E))&\leq2\frac{\displaystyle\int_{x\in B'}\H^d(\Pi_{\delta,x}(E))d\H^k(x)}{\H^k(B')}\\
&\leq\frac{2(2A)^dC\underline{R}(\delta)^{k-d}\overline{R}(\delta)^d}{D\underline{R}(\delta)^k}\H^d(E)\\
&\leq\frac{8^{d+1}C}{DR(\delta)^{2d}}\H^d(E).
\end{split}
\end{equation}
Since $C$ and $D$ depend only on $d$ and $k$, this achieves proving lemma~\ref{lemmaoptimalradialprojection}.
\end{proof}

In the special case when $E$ is irregular we also provide the following statement. It will be useful later to make the irregular part's measure vanish when approximating a given $d$-set with polyhedrons --- and thus allow us to give the main statement without restricting ourself to rectifiable sets only.

\begin{lemma}[Radial projection and irregular sets]\label{lemmairregularradialprojection}
Suppose that $0\leq d<k\leq n$, that $\delta$ is a $k$-dimensional polyhedron and that $E$ is a closed irregular $d$-set contained in $\delta$. Then, for $\H^k$ almost all $x\in\open{\delta}$, $\Pi_{\delta,x}(E)$ is also irregular.
\end{lemma}

Recall that an irregular set intersects a regular one only on a null set and that $\Pi_{\delta,x}(E)$ is contained in $\partial\delta$ --- which is $k-1$-rectifiable. As a consequence, $\H^d(\Pi_{\delta,x}(E))=0$ for $\H^k$ almost every $x$ as soon as $d=k-1$.

\begin{proof}
The first step of the proof is to show that for $\H^n$ almost any center, the radial projection of a given $d$-irregular set onto a given affine hyperplane is also $d$-irregular. Although this may not be the most natural way to prove it, we will rely on the well-known result about orthogonal projections of irregular sets onto linear subspaces --- again, see for instance Mattila's book~\cite{mattila}: for almost every linear $d$-plane $H$, the orthogonal projection of $E$ onto $H$ is a null $d$-set --- and conversely, any set verifying this property is $d$-irregular.

To define what we mean by ``almost every linear $d$-plane'' we will denote by $G(n,d)$ the Grassmannian manifold of all $d$-dimensional linear subspaces of $\mathbb{R}^n$ and consider the following Radon measure $\gamma_{n,d}$ on $G(n,d)$:
\begin{multline}\label{equationirregularradialprojectionA}
\forall X\subset G(n,d)\colon\gamma_{n,d}(X)=\underbrace{\H^n\times\ldots\times\H^n}_{\text{$d$ times}}\left(\left\{(v_1,\ldots,v_d)\in(\mathbb{R}^n)^d\colon\right.\right.\\\left.\left.\Vert v_i\Vert\leq1\text{ and }\vect(v_1,\ldots,v_d)\in X\right\}\right).
\end{multline}
By ``for almost every linear $d$-plane'' we are referring to a subset $Y\subset G(n,d)$ such that $\gamma_{n,d}(G(n,d)\setminus Y)=0$.

Suppose that $x\in\mathbb{R}$ and $y=(y_2,\ldots,y_n)\in\mathbb{R}^{n-1}$. For convenience, in what follows we will denote by $(x,y)$ the element $(x,y_2,\ldots,y_n)\in\mathbb{R}^n$. We will also use the following notations and variables:
\begin{itemize}
\item$a\in\mathbb{R}^{n-1}$, $0<\alpha<1$ and $\beta>0$;
\item$P$ is the affine hyperplane $\{1\}\times\mathbb{R}^{n-1}$ (identified with $\mathbb{R}^{n-1}$) and $p$ is the orthogonal projector onto $P$;
\item$\Pi_a$ is the radial projection onto $P$ centered at $(0,a)\in\mathbb{R}^n$;
\item$F$ is an irregular $d$-set (with $d\geq2$) contained in
\begin{equation}
D=[\alpha,1]\times[-\beta,\beta]^{n-1}.
\end{equation}
\end{itemize}

Firstly, we want to show that for $\H^{n-1}$ almost every $a$, $\Pi_a(F)$ is $d$-irregular. For that purpose, define
\begin{equation}
\phi_a\colon\begin{cases}
D\longrightarrow\mathbb{R}^n\\
(x,y)\longmapsto\left(\frac{1}{x},a+\frac{y-a}{x}\right),
\end{cases}
\end{equation}
and notice that $\Pi_a=p\circ\phi_a$. By putting $(x',y')=\phi_0(x,y)=\left(\frac{1}{x},\frac{y}{x}\right)$ we get
\begin{equation}
\phi_a(x,y)=\left(\frac{1}{x},a+\frac{y}{x}-\frac{a}{x}\right)=(x',a+y'-x'a).
\end{equation}
Besides, put $b=(1,a)$ and consider the three following affine maps onto $\mathbb{R}^n$:
\begin{alignat}{2}
p_a\colon z=(x,y)&\longmapsto z-\frac{\left<z,b\right>}{\Vert b\Vert^2}b=\left(x-\frac{x+\left<y,a\right>}{1+\Vert a\Vert^2},y-\frac{x+\left<y,a\right>}{1+\Vert a\Vert^2}a\right),\label{equationirregularradialprojectionB}\\
f_a\colon z=(x,y)&\longmapsto\left(x,y-xa\right)
\intertext{and}
\tau_a\colon z=(x,y)&\longmapsto z+b=\left(x+1,y+a\right).
\end{alignat}
Notice that
\begin{equation}
f_a\circ p_a(x',y')=\left(x'-\frac{x'+\left<y',a\right>}{1+\Vert a\Vert^2},y'-x'a\right)
\end{equation}
which in turn gives
\begin{equation}\label{equationirregularradialprojectionC}
p\circ\tau_a\circ f_a\circ p_a\circ\phi_0=p\circ\phi_a=\Pi_a.
\end{equation}
For convenience, let us identify $P$ with $\mathbb{R}^{n-1}$ and for $H\in G(n-1,d)$, suppose that $H$ (in fact, $\{1\}\times H$) is a $d$-dimensional linear subspace of $P$. Also, put $H'=\mathbb{R}\times H$ and denote by $p_H$ and $p_{H'}$ respectively the orthogonal projections onto $H$ and $H'$. Since $p_H\circ p=p\circ p_{H'}$ and $p_{H'}\circ f_a=f_{p_H(a)}\circ p_{H'}$ we deduce from~\eqref{equationirregularradialprojectionC} that
\begin{equation}
p_H\circ\Pi_a=p\circ\tau_a\circ f_{p_H(a)}\circ p_{H'}\circ p_a\circ\phi_0.
\end{equation}

Since $f_{p_H(a)}$ is $1+\Vert p_H(a)\Vert$-Lipschitz we get
\begin{equation}\label{equationirregularradialprojectionD}
\begin{split}
\H^d(p_H\circ\Pi_a(F))&=\H^d(p\circ\tau_a\circ f_{p_H(a)}\circ p_{H'}\circ p_a\circ\phi_0(F))\\
&\leq\H^d(f_{p_H(a)}\circ p_{H'}\circ p_a\circ\phi_0(F))\\
&\leq(1+\Vert a\Vert)^d\H^d(p_{H'}\circ p_a\circ\phi_0(F)).
\end{split}
\end{equation}
Also, recall that $p_a$ is defined in~\eqref{equationirregularradialprojectionB} as the orthogonal projector onto the linear hyperplane $H_a$ perpendicular to $b=(1,a)$. By putting $V(a,H)=H_a\cap H'$, $p_{H'}\circ p_a$ is the linear projection onto the linear $d$-plane $V(a,H)$. Suppose that $(v_1,\ldots,v_d)\in(\mathbb{R}^{n-1})^d$ are such that $\mathbb{R}\times\vect(u_1,\ldots,u_d)=H'$ and $\Vert u_i\Vert\leq 1$. Then
\begin{equation}
V(a,H)=\vect\left((-\left<u_1,a\right>,u_1),\ldots,(-\left<u_d,a\right>,u_d)\right),
\end{equation}
with $\left\Vert(-\left<u_1,a\right>,u_1)\right\Vert\leq 1+\Vert a\Vert$.

Take $r>0$, suppose that $X\subset(\mathbb{R}^{n-1}\cap B(0,r))\times G(n-1,d)$ and put
\begin{equation}
Y=\left\{V(a,H)\colon(a,H)\in X\right\}\subset G(n,d).
\end{equation}
In what follows, for convenience we will denote by $(\H^a)^b$ the product measure $\H^a\times\ldots\times\H^a$. Using inequalities of Hausdorff measure of Lipschitz images, we get the following, where $C$ and $C'$ depend only on $d$ and $n$:
\begin{equation}\label{equationirregularradialprojectionE}
\begin{split}
&(\H^{n-1}\times\gamma_{n-1,d})(X)\\
=&(\H^{n-1})^{d+1}\left(\left\{\left(a,u_1,\ldots,u_d\right)\colon\Vert u_i\Vert\leq1\text{ and }(a,\vect(u_1,\ldots,u_d))\in X\right\}\right)\\
\leq&\sum_{j\geq1}2^{j(n-1)}\H^{n-d-1}\times(\H^n)^d\left(\left\{\left((\left<u_{d+1},a\right>,\ldots,\left<u_{n-1},a\right>),\right.\right.\right.\\
&\qquad\left.(-\left<u_1,a\right>,u_1),\ldots,(-\left<u_d,a\right>,u_d))\right)\colon a\in X,2^{-j}<\Vert u_i\Vert\leq2^{-j+1},\\
&\qquad\left.\left.\vect(u_1,\ldots,u_n)=\{0\}\times\mathbb{R}^{n-1}\text{ and }(a,\vect(u_1,\ldots,u_d))\in X\right\}\right)\\
\leq&\sum_{j\geq1}2^{-jd(n-1)}\H^{n-d-1}\times(\H^n)^d\left(\left\{\left((\left<u_{d+1},a\right>,\ldots,\left<u_{n-1},a\right>),\right.\right.\right.\\
&\qquad\left.(-\left<u_1,a\right>,u_1),\ldots,(-\left<u_d,a\right>,u_d))\right)\colon a\in X,1/2<\Vert u_i\Vert\leq1,\\
&\qquad\left.\left.\vect(u_1,\ldots,u_n)=\{0\}\times\mathbb{R}^{n-1}\text{ and }(a,\vect(u_1,\ldots,u_d))\in X\right\}\right)\\
\leq&C\H^{n-d-1}\times(\H^n)^d\left(\left\{\left((\left<u_{d+1},a\right>,\ldots,\left<u_{n-1},a\right>),\right.\right.\right.\\
&\qquad\left.(-\left<u_1,a\right>,u_1),\ldots,(-\left<u_d,a\right>,u_d))\right)\colon\Vert u_i\Vert\leq1,\\
&\qquad\left.\left.\vect(u_1,\ldots,u_n)=\{0\}\times\mathbb{R}^{n-1}\text{ and }(a,\vect(u_1,\ldots,u_d))\in X\right\}\right)\\
\leq&C\H^{n-d-1}\times(\H^n)^d\left(\left\{\left((\left<u_{d+1},a\right>,\ldots,\left<u_{n-1},a\right>),\right.\right.\right.\\
&\qquad\left.(-\left<u_1,a\right>,u_1),\ldots,(-\left<u_d,a\right>,u_d))\right)\colon\Vert u_i\Vert\leq1,\\
&\qquad\left.\left.\vect(u_1,\ldots,u_n)=\{0\}\times\mathbb{R}^{n-1}\text{ and }(a,\vect(u_1,\ldots,u_d))\in X\right\}\right)\\
\leq&C(1+r)^{nd}\H^{n-d-1}\times(\H^n)^d\left(\left\{(b,v_1,\ldots,v_d)\colon b\in\mathbb{R}^{n-d-1},\Vert b\Vert\leq r,\right.\right.\\
&\left.\left.\qquad v_i\in\mathbb{R}^n,\Vert v_i\Vert\leq1\text{ and }\vect(v_1,\ldots,v_d)\in Y\right\}\right)\\
\leq&CC'r^{n-d-1}(1+r)^{nd}(\H^n)^d\left(\left\{(v_1,\ldots,v_d)\in(\mathbb{R}^n)^d\colon\Vert v_i\Vert\leq1\text{ and}\right.\right.\\
&\left.\left.\qquad\vect(v_1,\ldots,v_d)\in Y\right\}\right)\\
\leq&CC'(1+r)^{(n-1)(d-1)}\gamma_{n,d}(Y).
\end{split}
\end{equation}

Since $\phi_0$ is biLipschitz on $D$ and $F\subset D$ is $d$-irregular, $F'=\phi_0(F)$ is also $d$-irregular, which can be expressed as
\begin{equation}
\gamma_{n,d}\left(\left\{H\in G(n,d)\colon\H^d(p_H(F'))>0\right\}\right)=0.
\end{equation}
We are now ready to show that for $H^{n-1}$ almost all $a\in B(0,r)$, the radial projection $\Pi_a(F)$ is $d$-irregular. For that purpose, suppose that
\begin{equation}
X=\left\{(a,H)\in\mathbb{R}^{n-1}\times G(n-1,d)\colon\Vert a\Vert\leq r\text{ and }\H^d\left(p_H(F')\right)>0\right\},
\end{equation}
and let us compute the following quantity $M(r)$, using~\eqref{equationirregularradialprojectionD} and~\eqref{equationirregularradialprojectionE}:
\begin{equation}\label{equationirregularradialprojectionF}
\begin{split}
M(r)&=\Int_{\Vert a\Vert\leq r}\left(\Int_{H\in G(n-1,d)}\H^d(p_H\circ\Pi_a(F))d\gamma(n-1,d)(H)\right)d\H^{n-1}(a)\\
&=\Int_{\substack{\Vert a\Vert\leq r\\H\in G(n-1,d)}}\H^d(p_H\circ\Pi_a(F))d(\H^{n-1}\times\gamma(n-1,d))(a,H)\\
&\leq(1+r)^d\Int_{(a,H)\in X}\H^d(p_{H'}\circ p_a(F'))d\left(\H^{n-1}\times\gamma(n-1,d)\right)(a,H)\\
&\leq(1+r)^d\H^d(F')\left(\H^{n-1}\times\gamma(n-1,d)\right)(X)\\
&\leq CC'(1+r)^{(n-1)(d-1)}\H^d(F')\gamma_{n,d}(Y)\\
&\leq CC'(1+r)^{(n-1)(d-1)}\H^d(F')\gamma_{n,d}\left(\left\{H\in G(n,d)\colon\H^d(p_H(F'))>0\right\}\right)\\
&=0.
\end{split}
\end{equation}
Equation~\eqref{equationirregularradialprojectionF} is valid for any $r>0$, which is enough to prove that $\Pi_a(F)$ is $d$-irregular for $\H^{n-1}$ almost all $a\in\mathbb{R}^{n-1}$. It is also clear that all the above calculations could have been done with any radial projection centered at $(x,a)$, with $x<0$. As a consequence, for $\H^n$ almost all center $(a,x)$ (with $x\leq 0$), the radial projection of $F$ onto $P$ is $d$-irregular.

Let us resume proof of lemma~\ref{lemmairregularradialprojection}. Without loss in generality, by working in the affine subspace $\affine(\delta)$ we can assume that $k=n$. Fix $x\in\open{\delta}\setminus E$. Since $E$ is closed, one can find a ball $B(x)$ centered at $x$ such that $B\subset\open{\delta}\setminus E$. If we consider a face $F_i\in\F_{k-1}(\delta)$, included in the affine hyperplane $H_i$, by using the same notations as those in lemma's~\ref{lemmaoptimalradialprojection} proof we have
\begin{equation}
\forall x\in\open{\delta}\setminus E\colon\inf_{y\in\C^i(x)\cap E}\dist(x,H_i)>0.
\end{equation}
Using the above part of the proof, one can find a ball $B_i(x)\subset B(x)$ such that $\Pi_{\delta,y}(E\cap\C^i(y))$ is $d$-irregular for $\H^k$ almost all $y\in B_i(x)$. By iterating this argument over all faces of $\delta$, one can find a ball $B'(x)=\bigcap_i B_i(x)$, centered at $x$, such that for $\H^k$ almost all $y\in B'(x)$:
\begin{equation}
\Pi_{\delta,y}(E)\text{ is $d$-irregular.}
\end{equation}
Since $E$ is a null $k$-set (recall that $d<k$), by repeating over all $x\in\open{\delta}\setminus E$ this achieves proving the lemma.
\end{proof}

\section{Existence of a minimal candidate}

Before we start with the main result, we give ourself two handy tools that will allow us either to build a polyhedric mesh and a Lipschitz map that send a given $d$-set onto its $d$-dimensional subfaces, or to build a Lipschitz map that sends a given $d$-set onto the subfaces of an existing grid, each time with some kind of optimal control over the potential $d$-dimensional measure increase.

\subsection{Polyhedral approximation}

We will now proceed into proving the following analogous for compact $d$-sets of the classical polyhedral approximation theorem for integral currents. Notice that the requirements on $E$ are very minimalist: we do not even suppose that $E$ is rectifiable.

\begin{theorem}[Polyhedral approximation]\label{theorempolyhedralapproximation}
Suppose that $0<d<n$ and that $h\colon\mathbb{R}^n\rightarrow[1,+\infty[$ is continuous.

There is a positive constant $J>0$ such that for all open bounded domain $U\subset\mathbb{R}^n$, for all closed $d$-set $E\subset U$ and for all $\epsilon>0$, $R>0$, one can build a $n$-dimensional complex $S$ and a Lipschitz map $\phi\colon\mathbb{R}^n\rightarrow\mathbb{R}^n$ satisfying the following properties:
\begin{itemize}
\item $\phi\vert_{\mathbb{R}^n\setminus U}=\identity_{\mathbb{R}^n\setminus U}$ and $\Vert\phi-\identity_{\mathbb{R}^n}\Vert_\infty\leq\epsilon$;
\item $\R(S)\geq M$, $\overline{\R}(S)\leq J$ and the boundary faces $\F_\partial(S)$ of $S$ are the same as the ones of a dyadic complex;
\item $\phi(E)\subset\U(\F_d(S))$ and $\U(S)\subset U$;
\item $J_h^d(\phi(E))\leq (1+\epsilon)J_h^d(E)$.
\end{itemize}
\end{theorem}

\begin{proof}
To begin with, suppose that $E=E_R\cup E_I$, where $E_R$ is $d$-rectifiable, $E_I$ is $d$-irregular and $E_R\cap E_I=\emptyset$. Let us fix $\epsilon>0$, $\epsilon'>0$, $R>0$ and apply lemma~\ref{lemmamagneticprojection} to $E$: at $\H^d$ almost every point of $E_R$, one can find $r_{\max}(x)>0$, $\rho$ and $u$ such that for all $r<r_{\max}(x)$, inequality~\eqref{equationmagneticprojection} is true. Since $h$ is continuous over the compact set $\overline{U}$, one can find $A>0$ such that $1\leq h(x)\leq A$ for all $x\in U$, and for all $x\in U$ one can find $r'_{\max}(x)>0$ such that
\begin{equation}
\forall y\in B(x,r'_{\max}(x))\colon (1-\epsilon')h(x)\leq h(y)\leq (1+\epsilon')h(x).\label{equationpolyhedralapproximationA}
\end{equation}
Denote by $\B$ the collection of closed balls centered at a point $x$ of $E_R$ where $r_{\max}$ is defined, with radius at most $\min\left(\frac{r_{\max}(x)}{1+\rho},r'_{\max}(x),\frac{\epsilon}{2}\right)$. By a Vitali covering lemma, one can extract a countable subset $\widehat{\B}=\left\{B_i\colon i\in\mathbb{N}\right\}$ from $\B$ of pairwise disjoint balls such that
\begin{equation}
\H^d\left(E_R\setminus\bigcup_iB_i\right)=0.
\end{equation}

For each ball $B_i\in\widehat{\B}$ centered at $x_i$ with radius $r$, denote by $\rho_i$ and $u_i$ the constants given by lemma~\ref{lemmamagneticprojection} at $x_i$, put $r_i=\frac{r}{1+\rho_i}$ and consider the compact set
\begin{equation}
K_i=\C(x_i,r_i,u_i),
\end{equation}
as defined in~\eqref{equationcone}. Call $H_i$ the approximate tangent $d$-plane at $x_i$. Our upper bound on the radii of balls in $\B$ implies that
\begin{equation}\label{equationpolyhedralapproximationB}
\H^d(\Pi_{H_i,r_i\rho_i,K_i}(E\cap B(x_i,r_i+r_i\rho_i)\setminus K_i))\leq\epsilon'\H^d(E\cap B(x_i,r_i+r_i\rho_i)).
\end{equation}

Consider a finite subset $\overline{\B}$ from $\widehat{\B}$ such that
\begin{equation}\label{equationpolyhedralapproximationC}
\H^d\left(E_R\setminus\bigcup_{B\in\overline{\B}}B\right)\leq\epsilon'\H^d(E)
\end{equation}
and define the magnetic projections product (see proposal~\ref{proposalmagneticprojection})
\begin{equation}\label{equationpolyhedralapproximationD}
\psi_0=\prod_{B_i\in\overline{\B}}\Pi_{H_i,r_i\rho,K_i}.
\end{equation}
Notice that $\support\Pi_{H_i,r_i\rho,K_i}\subset(K_i)_{r_i\rho}\subset B_i$ --- which are pairwise disjoint balls of radii at most $\frac{\epsilon}{2}$ --- and $\Pi_{H_i,r_i\rho,K_i}$ is $\left(2+\frac{u_i}{\rho_i}\right)$-Lipschitz, so $\psi_0$ is $\gamma$-Lipschitz with
\begin{gather}
\gamma=2+\max_{B_i\in\overline{\B}}\frac{u_i}{\rho_i},\\
\Vert\psi_0-\identity_U\Vert_\infty\leq\frac{\epsilon}{2},\label{equationpolyhedralapproximationE}
\end{gather}
and the definition of $\psi_0$ does not depend upon the choice of the order of multiplication in~\eqref{equationpolyhedralapproximationF}.

Suppose that $\alpha>0$. If $\alpha$ is small enough, one can build in each $K_i$ a dyadic complex $S_i$ of stride $\alpha$ (see definition~\ref{definitiondyadiccomplex}) in an orthonormal basis centered at $x_i$ with $d$ vectors parallel to $H_i$. There is also a constant $\alpha_i$ depending on $u_i$ and $r_i$ such that, if $\alpha<\alpha_i$ and by taking in $S_i$ every possible dyadic cube included in $K_i$:
\begin{equation}\label{equationpolyhedralapproximationF}
\H^d(\psi_0(E)\cap K_i\setminus\U(S_i))\leq\epsilon'\H^d(\psi_0(E)\cap K_i)\leq\epsilon'\H^d(E\cap B_i).
\end{equation}
By putting $\alpha_{\max}=\min_i\alpha_i$ and by taking $\alpha<\alpha_{\max}$, one can build all these dyadic complexes $S_i$ of stride $\alpha$ such that $\Sigma_2=\bigcup_i S_i$ is a $n$-dimensional complex obtained as a finite union of dyadic complexes verifying~\eqref{equationpolyhedralapproximationF}. Let us define:
\begin{equation}
\begin{split}
E_1&=E\setminus\bigcup_{B_i\in\overline{\B}}B_i,\\
E_2&=E\cap\bigcup_{B_i\in\overline{\B}}B_i\setminus K_i,\\
E_3&=\{x\in E\cap\bigcup_{B_i\in\overline{\B}}K_i\colon\psi_0(x)\notin\U(S_i)\},\\
E_4&=\{x\in E\cap\bigcup_{B_i\in\overline{\B}}K_i\colon\psi_0(x)\in\U(S_i)\}.
\end{split}
\end{equation}
Notice that $E=E_1\sqcup E_2\sqcup E_3\sqcup E_4$, that $\psi_0\vert_{E_1}=\identity_{E_1}$ and by~\eqref{equationpolyhedralapproximationC}, \eqref{equationpolyhedralapproximationB} and~\eqref{equationpolyhedralapproximationF} we also have the following inequalities:
\begin{equation}
\begin{split}
\H^d(\psi_0(E_1\cap E_R))&=\H^d(E_1\cap E_R)\leq\epsilon'\H^d(E),\\
\H^d(\psi_0(E_2))&\leq\epsilon'\H^d(E_2)\leq\epsilon'\H^d(E),\\
\H^d(\psi_0(E_3))&\leq\H^d(E_3)\leq\epsilon'\H^d(E).
\end{split}
\end{equation}
By summing and putting $\epsilon''=3\epsilon'A$ we obtain
\begin{equation}\label{equationpolyhedralapproximationG}
J_h^d\left(\psi_0((E_1\cap E_R)\sqcup E_2\sqcup E_3)\right)\leq 3\epsilon'AJ_h^d(E)=\epsilon''J_h^d(E).
\end{equation}
On the other hand, $\psi_0\vert_{\U(S_i)}$ is the orthogonal projector onto $H_i$ with $\U(S_i)\subset B_i$. Since each $B_i$ has radius at most $r'_{\max}(x_i)$, by~\eqref{equationpolyhedralapproximationB} we have
\begin{equation}\label{equationpolyhedralapproximationH}
J_h^d(\psi_0(E_4))\leq(1+\epsilon')\H^d(\psi_0(E_4))\leq(1+\epsilon')\H^d(E_4)\leq(1+\epsilon')^2J_h^d(E),
\end{equation}
and we can notice already that $\psi_0(E_4)\subset\U(\F_d(S))$, since we chose the orientation of $S_i$ parallel to $H_i$.

By hypothesis, $E$ and $\partial U$ are compact, and since $E\cap\partial U=\emptyset$ we have
\begin{equation}
a=\inf_{(x,y)\in E\times\partial U}\dist(x,y)>0.
\end{equation}
Consider theorem~\ref{theoremfusion} (in what follows, $\rho$ is the minimal distance required to merge dyadic grids together, and $c_1$ the constant used to control the upper radii) and suppose that we took
\begin{equation}
\alpha<\min\left(\alpha_{\max},\frac{a}{4\sqrt{n}},\frac{\min_i{\rho_i}}{16\sqrt{n}},\frac{\min_i{\rho_i}}{2\rho},\frac{R}{2c_1\sqrt{n}},\frac{\epsilon}{2c_1\sqrt{n}}\right)
\end{equation}
when building our dyadic grids $S_i$. Fix an arbitrary orthonormal basis in $\mathbb{R}^n$. By taking all possible cubes of stride $\alpha$ in this basis that are included in $U$ and disjoint with all the $(K_i)_{r_i\rho_i/2}$, one can build a dyadic complex $\Sigma_1$ such that:
\begin{equation}
\begin{split}
\U(\Sigma_1)&\subset U\setminus\bigcup_i(K_i)_{r_i\rho_i/2},\\
\U(\Sigma_1)&\supset E_1,\\
\U(\Sigma_1)&\supset\bigcup_i\left((K_i)_{r_i\rho_i}\setminus(K_i)_{7r_i\rho_i/8}\right).
\end{split}
\end{equation}
By using theorem~\ref{theoremfusion} separately in each $(K_i)_{r_i\rho_i}$ (which are pairwise disjoint) we can build a complex $S$ such that $\Sigma_1\sqcup\Sigma_2\subset S$, $E\subset\U(S)\subset U$, $\R(S)\geq M$ (with $M$ depending only on $n$) and
\begin{equation}\label{equationpolyhedralapproximationI}
\overline{\R}(S)\leq c_1\overline{\R}(\Sigma_1\sqcup\Sigma_2)\leq\min\left(R,\frac{\epsilon}{2}\right).
\end{equation}

Put $F_0=\psi_0(E)$ and let us reason by induction. Suppose that at rank $k\in\{1,\ldots,n-d\}$ we have found a Lipschitz map $\psi_{k-1}$ which verifies, by putting $F_{k-1}=\psi_{k-1}(E)$:
\begin{itemize}
\item $\psi_{k-1}\vert_{E_4}=\psi_0\vert_{E_4}$;
\item $\psi_{k-1}(E_1\cap E_I)$ is $d$-irregular and there is a constant $C>0$ depending only on $d$ and $n$ such that $J_h^d(\psi_{k-1}((E_1\cap E_R)\sqcup E_2\sqcup E_3))\leq C^{k-1}\epsilon''J_h^d(E)$;
\item $F_{k-1}\subset\U(\F_{n-k+1}(S))$.
\end{itemize}
Notice that by construction and~\eqref{equationpolyhedralapproximationG}, $\psi_0$ verifies all three properties at rank $k=1$.

For all $\delta\in\F_{n-k+1}(S)$ we can apply lemma~\ref{lemmairregularradialprojection} and lemma~\ref{lemmaoptimalradialprojection} to respectively $\psi_{k-1}(E_1\cap E_I)\cap\delta$ and $\psi_{k-1}((E_1\cap E_R)\sqcup E_2\sqcup E_3)\cap\delta$, and find a center $x_\delta\in\open{\delta}\setminus F_{k-1}$ such that $\Pi_{\delta,x_\delta}\circ(\psi_{k-1}(E_1\cap E_I)\cap\delta)$ is also $d$-irregular and
\begin{multline}\label{equationpolyhedralapproximationJ}
\H^d(\Pi_{\delta,x_\delta}\circ\psi_{k-1}(((E_1\cap E_R)\sqcup E_2\sqcup E_3)\cap\delta))\leq\\
K_{d,k}\R(S)^{-2d}\H^d(\psi_{k-1}(((E_1\cap E_R)\sqcup E_2\sqcup E_3)\cap\delta))
\end{multline}
where $K_{d,k}$ depends only on $d$ and $k$. Notice that $\psi_0\vert_{E_4}$ is defined as the orthogonal projector onto $H_i$ inside $K_i$, and since we supposed that $\psi_{k-1}\vert_{E_4}=\psi_0\vert_{E_4}$ we have
\begin{equation}
\psi_{k-1}(E_4)\subset\U(\F_d(S))\subset\U(\F_{n-k}(S)).
\end{equation}
As a consequence, for all subface $\delta\in\F_{n-k+1}(\Sigma_2)$ we have $\psi_{k-1}(E_4)\cap\delta\subset\partial\delta$, and since $E_4\subset\U(\Sigma_2)$:
\begin{equation}\label{equationpolyhedralapproximationK}
\forall\delta\in\F_{n-k+1}(S)\colon\Pi_{\delta,x_\delta}\vert_{\psi_0(E_4)\cap\delta}=\identity_{\psi_0(E_4)\cap\delta}.
\end{equation}

Since $E$ is closed, for all $\delta\in\F_{n-k+1}(S)$ we can find some $n-k$-dimensional ball $B_\delta\subset\delta$ such that $B_\delta\cap F_{k-1}=\emptyset$. Since $\Pi_{\delta,x_\delta}\vert_{\delta\setminus B_\delta}$ is Lipschitz, by applying lemma~\ref{lemmaholeextension} we can extend it on $\delta$ as a Lipschitz map $\psi_\delta$. And since $\psi_\delta\vert_{\partial\delta}=\identity_{\partial\delta}$, by applying lemma~\eqref{lemmacomplexextension} to the $n-k+1$-dimensional complex $\F_{n-k+1}(S)$, we can build a Lipschitz extension $\psi$ on $U$.

Put $\psi_k=\psi\circ\psi_{k-1}$ and let us check that $\psi_k$ verifies all three induction hypothesis:
\begin{itemize}
\item $\psi_k\vert_{E_4}=\psi\circ\psi_{k-1}\vert_{E_4}=\psi\circ\psi_0\vert_{E_4}=\psi_0\vert_{E_4}$ by~\eqref{equationpolyhedralapproximationK};
\item we already know that $\psi_k(E_1\cap E_I)$ is $d$-irregular. Since $\R(S)\geq M$, by~\eqref{equationpolyhedralapproximationJ} and by putting $C=AM^{-2d}\max_k K_{d,k}$ we also obtain
\begin{equation}
\begin{split}
J_h^d(\psi_k((E_I\cap E_1)\sqcup E_2\sqcup E_3))&\leq CJ_h^d(\psi_{k-1}((E_I\cap E_1)\sqcup E_2\sqcup E_3))\\
&\leq C^k\epsilon''J_h^d(E);
\end{split}
\end{equation}
\item by construction, for all $\delta\in\F_{n-k+1}(S)$ we have $\psi_k(\delta)\subset\partial\delta\in\F_{n-k}(S)$. Since we supposed that $F_{k-1}\subset\U(\F_{n-k+1}(S))$, we also have $F_k=\psi_k(E)\subset\U(\F_{n-k}(S))$, which achieves proving the induction.
\end{itemize}

Take $k=n-d$, put $\phi=\psi_{n-d}$ and recall that we built $\phi$ as the product $\phi=f\circ\psi_0$ where $f$ is such that $f(\delta)\subset\delta$ for all $\delta\in\F(S)$. Using~\eqref{equationpolyhedralapproximationI} we get $\Vert f-\identity_{\mathbb{R}^n}\Vert_\infty\leq\frac{\epsilon}{2}$ and by~\eqref{equationpolyhedralapproximationE}:
\begin{equation}
\Vert\phi-\identity_{\mathbb{R}^n}\Vert_\infty\leq\Vert\psi_0-\identity_{\mathbb{R}^n}\Vert_\infty+\Vert f-\identity_{\mathbb{R}^n}\Vert_\infty\leq\epsilon.
\end{equation}
Notice that since $\phi(E_1\cap E_I)$ is $d$-irregular and included in $\U(\F_d(S))$ (which is $d$-rectifiable) then $\H^d(\phi(E_1\cap E_I))=0$. Using~\eqref{equationpolyhedralapproximationH} we finally get:
\begin{equation}
\begin{split}
J_h^d(\phi(E))&\leq J_h^d(\phi(E_1\cap E_I))+J_h^d(\phi((E_1\cap E_R)\sqcup E_2\sqcup E_3\sqcup E_4))\\
&\leq J_h^d(\phi((E_1\cap E_R)\sqcup E_2\sqcup E_3))+J_h^d(E_4)\\
&\leq C^{n-d}\epsilon''J_h^d(E)+J_h^d(\psi_0(E))\\
&\leq (C^{n-d}\epsilon''+(1+\epsilon')^2)J_h^d(E).
\end{split}
\end{equation}
By taking $\epsilon'$ small enough such that $C^{n-d}\epsilon''+(1+\epsilon')^2\leq 1+\epsilon$, this achieves proving theorem~\ref{theorempolyhedralapproximation}.
\end{proof}

The following lemma is very similar, except that the polyhedric mesh is fixed. The control over the potential measure increase is given by a multiplicative constant depending on the shape of the polyhedrons and subfaces of the mesh.

\begin{lemma}[Polyhedral deformation]\label{lemmapolyhedraldeformation}
Suppose that $0\leq d<n$, that $U\subset\mathbb{R}^n$ is an open bounded domain and that $S$ is a $n$-dimensional complex such that $\U(S)\subset U$.

There exists a constant $K>0$ depending only on $d$ and $n$ such that for all closed $d$-set $E\subset\U(S)$, one can build a Lipschitz map $\phi\colon\mathbb{R}^n\rightarrow\mathbb{R}^n$ satisfying the following properties:
\begin{itemize}
\item$\phi\vert_{\mathbb{R}^n\setminus U}=\identity_{\mathbb{R}^n\setminus U}$ and for all subface $\alpha\in\F(S)$: $\phi(\alpha)=\alpha$ and $\phi\vert_\alpha=\identity_\alpha$ if $\dim\alpha\leq d$;
\item$\phi(E)\subset\U(\F_d(S))$;
\item$\H^d(\phi(E))\leq K\R(S)^{-2d(n-d)}\H^d(E)$ and for all subface $\alpha\in\F(S)$: $\H^d(\phi(E\cap\open{\alpha}))\leq K\R(S)^{-2d(n-d)}\H^d(E\cap\open{\alpha})$.
\end{itemize}
\end{lemma}

The proof is pretty straightforward: we just have to use an induction reasoning like the one in the above proof of  theorem~\ref{theorempolyhedralapproximation}.

\begin{proof}
By building optimal radial projections in subfaces of dimension $n$, $n-1$, $\ldots$ till dimension $d+1$ and extend them on $\mathbb{R}^n$ using lemma~\ref{lemmacomplexextension} we build a map $\phi$ that verifies all the required topological constraints, and such that
\begin{equation}
\forall\alpha\in\F(S)\colon\H^d(\phi(E\cap\open{\alpha}))\leq K\R(S)^{-2d(n-d)}\H^d(E\cap\open{\alpha})
\end{equation}
where $K$ depends only on $d$ and $n$.
\end{proof}

Our two previous polyhedral approximation and deformation statements (theorem~\ref{theorempolyhedralapproximation} and lemma~\ref{lemmapolyhedraldeformation}) are not complete, in the sense that the set we obtain in the end may not be made of complete polyhedrons, but instead may contain ``holes''. In each polyhedron that is not completely covered, it is possible to continue our radial projections in the previous dimension till all remaining subfaces are completely covered. At the end, the set we obtain is a finite union of subfaces of dimension at most $d$ (i.e. a $d$-dimensional skeleton, as introduced in section~1).

\begin{lemma}[Polyhedral erosion]\label{lemmapolyhedralerosion}
Suppose that $0\leq d<n$, that $U\subset\mathbb{R}^n$ is an open bounded domain and that $S$ is a $n$-dimensional complex such that $\U(S)\subset U$.

For all closed set $E\subset\U(\F_d(S))$ one can build a Lipschitz map $\phi\colon\mathbb{R}^n\rightarrow\mathbb{R}^n$ satisfying the following properties:
\begin{itemize}
\item$\phi\vert_{\mathbb{R}^n\setminus U}=\identity_{\mathbb{R}^n\setminus U}$ and for all subface $\alpha\in\F(S)$: $\phi(\alpha)=\alpha$, and $\phi\vert_\alpha=\identity_\alpha$ or $\phi(\alpha\cap E)\subset\partial\alpha$;
\item there is a $d$-dimensional skeleton $S'$ of $S$ such that $\phi(E)=\U(S')$;
\item $\H^d(\phi(E))\leq\H^d(E)$.
\end{itemize}
\end{lemma}

Later, this lemma will be used in conjunction with theorem~\ref{theorempolyhedralapproximation} or lemma~\ref{lemmapolyhedraldeformation} to restrict ourselves to a finite subclass of competitors for which finding a minimal set is trivial.

\begin{proof}
For $\leq j\leq d$ and $F\subset\mathbb{R}^n$, put
\begin{align}
S_j(F)&=\bigcup_{\substack{\delta\in\F_j(S)\\F\cap\delta=\delta}}\delta,&S'_j(F)&=\bigcup_{\substack{\delta\in\F_d(S)\\F\cap\open{\scriptstyle\delta}\neq\open{\scriptstyle\delta}}}F\cap\delta.
\end{align}
Notice that when $F\subset\U(\F_j(S))$, $S_j(F)\cap S'_j(F)\subset\U(\F_{j-1}(S))$ and we can find $F'\subset\U(\F_{j-1}(S))$ such that $F=S_j(F)\cup S'_j(F)\cup F'$.

We will use again a similar argument as in lemma~\ref{lemmacomplexextension}. Put $\psi_0=\identity_{\mathbb{R}^n}$, $E_d=E$ and notice that since $S$ is a complex, $S'_j(E_d)=\emptyset$ for all $j>d$. Let us reason by decreasing induction over $j$, and suppose that at rank $j\in\{1,\ldots,d\}$ we have built a Lipschitz map $\psi_j$ over $\mathbb{R}^n$ such that, by putting $E_j=\psi_j(E)$ we have:
\begin{equation}\label{equationpolyhedralerosionA}
\H^d(E_j)\leq\H^d(E)\qquad\text{and}\qquad\forall k\in\{j+1,\ldots,n\}\colon S'_k(E_j)=0. 
\end{equation}

Put
\begin{equation}
T=\left\{\alpha\in\F_j(S)\colon E_j\cap\open{\alpha}\neq\emptyset\text{ and }E_j\cap\open{\alpha}\neq\open{\alpha}\right\}.
\end{equation}
If $T=\emptyset$ we have finished. If not, since $E_j$ is closed then for all $\alpha\in T$ we can find a $d$-dimensional open ball $B\subset\open{\alpha}\setminus E_j$ centered at $x_\alpha$. By using lemma~\ref{lemmaholeextension} we can extend $\Pi_{\alpha,x_\alpha}\vert_{\alpha\setminus B}$ over $\alpha$ and obtain a Lipschitz map $\psi_\alpha$ such that
\begin{equation}
\begin{split}
\psi_\alpha(E_j\cap\alpha)&\subset\U(\F_{j-1}(S))\\
\H^d(\psi_\alpha(E_j\cap\alpha))&=0\leq\H^d(E_j\cap\alpha).
\end{split}
\end{equation}
Suppose that $\alpha\in T$, $k>j$ and that $\beta\in\F_k(S)$ is such that $\open{\alpha}\cap\beta\neq\emptyset$. Since $S$ is a complex, this implies that $\alpha\subset\partial\beta\subset\beta$. By~\eqref{equationpolyhedralerosionA}, either $E_j\cap\beta=\beta$ or $E_j\cap\open{\beta}=\emptyset$ and since $E_j\cap\open{\alpha}\neq\open{\alpha}$ the second case is true. As we have previously done in lemma~\ref{lemmacomplexextension}, we can build Lipschitz extensions of all the $\psi_\alpha$ (for $\alpha\in T$) over $\mathbb{R}^n$ with pairwise disjoint supports and such that
\begin{equation}
\support\psi_\alpha\cap\left((\mathbb{R}^n\setminus U)\cup(F_j\setminus S'_j(F_j))\right)=\emptyset.
\end{equation}
Put
\begin{equation}
\begin{split}
\psi&=\prod_{\alpha\in T}\psi_\alpha\\
\psi_{j-1}&=\psi\circ\psi_j.
\end{split}
\end{equation}
Since $\psi\vert_{F_j\setminus S'_j(F_j)}=\identity_{F_j\setminus S'_j(F_j)}$ and $\psi\vert_{S'_j(F_j)}$ is a product of extensions of radial projections in $j$-dimensional subfaces of $S$, then for all $k>j-1$:
\begin{equation}
S'_k(S_j(E_{j-1}))=S'_k(S'_j(E_{j-1}))=\emptyset.
\end{equation}
Besides, since $E_{j-1}=S_j(E_{j-1})\cup S'_j(E_{j-1})\cup E'$ where $E'\subset\U(\F_{j-1}(S))$ then $S'_k(E')=\emptyset$ and we get:
\begin{equation}
S'_k(E_{j-1})=S'_k(S_j(E_{j-1}))\cup S'_k(S'_j(E_{j-1}))\cup S'_k(E')=\emptyset.
\end{equation}
Also, it is clear that $\H^d(E_{j-1})\leq\H^d(E)$ because $E_{j-1}\subset E$, which achieves proving the induction.

If we iterate the above process till rank $j=0$ and put $\phi=\psi_0$, for all $k>0$ we have $S'_k(\phi(E))=\emptyset$, which is enough to conclude.
\end{proof}

\subsection{Limits of uniformly concentrated minimizing sequences}

In what follows we give a way to convert any minimizing sequence of elements of $\mathfrak{E}$ into another minimizing sequence of polyhedric and quasiminimal competitors, with uniform constants (depending only on dimensions $d$ and $n$). Notice that the following lemma may prove to be more useful in ``real problems'' than theorem~\ref{theoremexistence}, because it gives more control over the topological constraint embedded in $\mathfrak{F}$, especially when involving the boundary of $U$.

\begin{lemma}[Polyhedral optimization]\label{theorempolyhedraloptimization}
Suppose that $0<d<n$ and that $U\subset\mathbb{R}^n$.

There is a positive constant $M'>0$ (depending only on $d$ and $n$) such that
\begin{itemize}
\item for all continuous function $h\colon U\rightarrow[1,M]$,
\item for all relatively closed $d$-subset $E\subset U$,
\item for all relatively compact subset $V\subsubset\open{U}$ and for all $\epsilon>0$,
\end{itemize}
one can find a $n$-dimensional complex $S$ and a subset $E''\subset U$ satisfying the following properties:
\begin{itemize}
\item $E''$ is a $\diam(U)$-deformation of $E$ over $U$ and by putting $W=\open{\U(S)}$ we have $V\subsubset W\subsubset U$ and there is a $d$-dimensional skeleton $S'$ of $S$ such that $E''\cap\overline{W}=\U(S')$;
\item $J_h^d(E'')\leq(1+\epsilon)J_h^d(E)$;
\item there are $d+1$ complexes $S^0,\ldots,S^d$ with $S^l\subset\F_l(S)$ such that, by putting
\begin{align}
\begin{cases}
E^d=\U(S^d)\cap W\\
E^l=\U(S^l)\cap W^l
\end{cases}
&&
\begin{cases}
W^d=W\\
W^{l-1}=W^l\setminus E^l,
\end{cases}
\end{align}
then $E''\cap W=E^d\sqcup E^{d-1}\sqcup\ldots\sqcup E^0$ and for all $l\in\{0,\ldots,d\}$, $E^l$ is $(MM',\diam(W))$-quasiminimal over $W^l$ for $\H^l$. Furthermore, $E^l$ is optimal in the sense that if all the $E^{l'}$ are fixed for $l'>l$, any deformation of $E$ over $W^{l'}$ verifying the same above properties cannot decrease $J_h^d(E^l)$.
\end{itemize}
\end{lemma}

\begin{proof}
To begin with, we can always suppose that $U$ is bounded. Otherwise, take an open bounded neighborhood $U'\subset U$ of $V$ such that $V\subsubset U'$ and replace $U$ by $U'$. That way, we can assume that $\H^d(E)<\infty$. Since $V\subsubset U$ we have
\begin{equation}\label{equationpolyhedraloptimizationA}
A=\inf_{(x,y)\in\partial U\times\partial V}\dist(x,y)>0,
\end{equation}
which means that in any orthonormal basis and for $R<\frac{A}{8\sqrt{N}}$ one can build a dyadic complex $T$ of stride $R$ such that $V\subsubset\U(T)\subsubset U$.

Fix $\epsilon>0$, put
\begin{equation}\label{equationpolyhedraloptimizationB}
R=\frac{A}{8\sqrt{n}},
\end{equation}
and apply theorem~\ref{theorempolyhedralapproximation} to the closed $d$-set $E\cap\overline{V}$ in the open domain $U$, with the above constant $R$: we build a dyadic complex $S$ such that $\overline{\R}(S)<R$, $\R(S)>J$ and $\U(S)\subsubset U$, and a Lipschitz map $\psi_1$ such that $\Vert\psi_1-\identity_{\mathbb{R}^n}\Vert_\infty<\epsilon$, $\psi_1(E\cap\overline{V})\subset\U(\F_d(S))$ and $J_h^d(\psi_1(E\cap\overline{V}))\leq(1+\epsilon)J_h^d(E\cap\overline{V})$. Using lemma~\ref{lemmapolyhedralerosion} with $\psi_1(E\cap\overline{V})$, we build a Lipschitz map $\psi_2$ such that $\psi_2\circ\psi_1(E\cap\overline{V})=\U(S')$ where $S'$ is a $d$-dimensional skeleton of $S$ and
\begin{equation}
J_h^d(\psi_2\circ\psi_1(E\cap\overline{V}))\leq(1+\epsilon)J_h^d(E\cap\overline{V}).
\end{equation}
If we build an additional layer of cubes around $S$, and by stopping the radial projections of theorem~\ref{theorempolyhedralapproximation} and lemma~\ref{lemmapolyhedralerosion} at dimension $n-1$ in the boundary faces of $S$ we can even assume that
\begin{equation}
J_h^d(\psi_2\circ\psi_1(E))\leq (1+\epsilon)J_h^d(E)
\end{equation}
and
\begin{equation}
\psi_2\circ\psi_1\vert_{\mathbb{R}^n\setminus\open{\U(S)}}=\identity_{\mathbb{R}^n\setminus\open{\U(S)}}.
\end{equation}
Later, we will implicitly make the same assumptions when using lemmas~\ref{lemmapolyhedraldeformation} and~\ref{lemmapolyhedralerosion}.

Since $\F_\partial(S)$ is the same as a dyadic complex, and since $\overline{\R}(S)<R$, by~\eqref{equationpolyhedraloptimizationA} and~\eqref{equationpolyhedraloptimizationB} we can add dyadic cubes around $S$ until
\begin{equation}
V\subsubset\U(S)\subsubset U.
\end{equation}
Put $W=\open{\U(S)}$ and $E'=\psi_2\circ\psi_1(E)$, and recall that by lemma~\ref{lemmapolyhedralerosion}:
\begin{equation}
\forall\delta\in S\colon\psi_2(\delta)\subset\delta.
\end{equation}
This implies that $\Vert\psi_2-\identity_{\mathbb{R}^n}\Vert_\infty\leq\overline{\R}(S)<R$, and we get
\begin{equation}
\Vert\psi_2\circ\psi_1-\identity_{\mathbb{R}^n}\Vert_\infty<2R<A.
\end{equation}
By~\eqref{equationpolyhedraloptimizationA} and using proposal~\ref{proposalautomaticdeformation} with $\psi_2\circ\psi_1$ and the identity deformation over $U$ we build a deformation $(\phi_t)$ over $U$ such that $\phi_1=\psi_2\circ\psi_1$, and $E'$ is an Almgren competitor of $E$ such that $E'\setminus W=E\setminus W$.

Consider the set $\mathfrak{S}$ of subsets of $U$ obtained as an union of $E\setminus W$ with a $d$-dimensional skeleton of $S$:
\begin{equation}
\mathfrak{S}=\{\U(T)\cup(E\setminus W)\colon T\subset\F_d(S)\cup\ldots\cup\F_0(S)\},
\end{equation}
and the set $\mathfrak{E}$ of competitors of $E$ obtained by a deformation with support in $W$:
\begin{equation}
\mathfrak{E}=\{\phi_1(E)\colon\text{$(\phi_t)$ is a $\diam(W)$-deformation over $U$ and $\support(\phi)\subset W$}\}.
\end{equation}
Notice that $\mathfrak{S}\cap\mathfrak{E}$ is finite since $\F(S)$ is finite, and non-empty since it contains $E'$. Then we can find $E''\in\mathfrak{S}\cap\mathfrak{E}$ such that
\begin{equation}
J_h^d(E'')=\min\{J_h^d(F)\colon F\in\mathfrak{S}\cap\mathfrak{C}\},
\end{equation}
and furthermore
\begin{equation}
J_h^d(E'')\leq J_h^d(E')\leq(1+\epsilon)J_h^d(E).
\end{equation}

Let us check that $E''$ meets all the announced quasiminimality requirements. Suppose that $F$ is an Almgren competitor of $E''$ obtained by a $\diam(W)$-deformation $(\phi_t)$ over $W$. Since $F$ is also an Almgren competitor of $E$ we have $F\in\mathfrak{E}$. By applying lemmas~\ref{lemmapolyhedraldeformation} and~\ref{lemmapolyhedralerosion} to $F$ and $S$, as we did previously with $E$ we can build an Almgren competitor $F'\in\mathfrak{E}\cap\mathfrak{S}$ of $F$ obtained by a deformation $(\psi_t)$ over $W$ such that for all subface $\alpha\in\F(S)$:
\begin{equation}
\H^d(\psi_1(F\cap\open{\alpha}))\leq K\R(S)^{-2d(n-d)}\H^d(F\cap\open{\alpha})\leq K'\H^d(F\cap\open{\alpha}),
\end{equation}
where $K'=KJ^{-2d(n-d)}$ depends only on $d$ and $n$. Recall that $E''\cap W$ is an union of subfaces of dimension at most $d$ of $S$. Then, for all subface $\alpha$ of dimension at least $d+1$, $F\cap E''\cap\open{\alpha}=\emptyset$ and as a consequence:
\begin{equation}
F\setminus E''=\left(\bigsqcup_{\substack{\alpha\in\F(S)\\\dim(\alpha)>d}}F\cap\open{\alpha}\right)\sqcup\left(\bigsqcup_{\substack{\alpha\in\F(S)\\\dim(\alpha)\leq d}}(F\setminus E'')\cap\open{\alpha}\right).
\end{equation}
Notice that the Lipschitz maps given by lemmas~\ref{lemmapolyhedraldeformation} and~\ref{lemmapolyhedralerosion} are such that for all subface $\alpha$ of dimension at most $d$, $\psi_1\vert_\alpha=\identity_\alpha$ or $\psi_1(\alpha\cap F)\subset\partial\alpha$, which gives:
\begin{flsplit}\label{equationpolyhedraloptimizationC}
\H^d(\psi_1(F\setminus E''))&=\left(\sum_{\substack{\alpha\in\F(S)\\\dim(\alpha)>d}}\H^d(\psi_1(F\cap\open{\alpha}))\right)\shovebreak{+d\left(\sum_{\substack{\alpha\in\F(S)\\\dim(\alpha)\leq d}}\H^d(\psi_1((F\setminus E'')\cap\open{\alpha}))\right)}
&\leq\left(\sum_{\substack{\alpha\in\F(S)\\\dim(\alpha)>d}}K'\H^d(F\cap\open{\alpha})\right)\shovebreak{+\left(\sum_{\substack{\alpha\in\F(S)\\\dim(\alpha)\leq d}}\H^d((F\setminus E'')\cap\open{\alpha})\right)}
&\leq\max(K',1)\sum_{\alpha\in\F(S)}\H^d((F\setminus E'')\cap\open{\alpha})\\
&=K'\H^d(F\setminus E'').
\end{flsplit}
Since $F'\in\mathfrak{E}\cap\mathfrak{S}$ we have $J_h^d(F')\geq J_h^d(E'')$, and more precisely, by removing $E''\cap F'$:
\begin{equation}
J_h^d(E''\setminus F')\leq J_h^d(F'\setminus E'').
\end{equation}
Besides, $F'\setminus E''=\psi_1(F)\setminus E''\subset\psi_1(F\setminus E'')$ because $\psi_1(E'')=E''$ (recall that $E''\cap\overline{W}$ is an union of subfaces of $S$, and that by lemma~\ref{lemmapolyhedraldeformation}, for all $\alpha\in\F(S)$, $\psi_1(\alpha)=\alpha$) and $J_h^d(E''\setminus F')\leq J_h^d(\psi_1(F\setminus E''))$. Using our bounds on $h$ and~\eqref{equationpolyhedraloptimizationC} we get
\begin{equation}\label{equationpolyhedraloptimizationD}
\H^d(E''\setminus F')\leq M\H^d(\psi_1(F\setminus E''))\leq KM\H^d(F\setminus E'').
\end{equation}

Suppose that $\delta$ is a subface of $S$ of dimension at least $d+1$. Notice that $\psi_1(F\cap\open{\delta})$ is included in $\U(\F_d(\delta))$, and that by lemma~\eqref{lemmapolyhedraldeformation}:
\begin{equation}
\H^d(\psi_1(F\cap\open{\delta}))\leq K'\H^d(F\cap\open{\delta}).
\end{equation}
Conversely, if $\alpha\in\F_d(S)$ then either $\alpha\in E''\cap F'$ or $\open{\alpha}\cap E''\cap F'=\emptyset$ since $E''\cap\overline{W}$ and $F'\cap\overline{W}$ are both unions of subfaces of $S$. In the first case, the topological properties of the Lipschitz map given by lemmas~\ref{lemmapolyhedraldeformation} and~\ref{lemmapolyhedralerosion} imply that
\begin{equation}
\alpha\setminus F\subset\psi_1\left(\bigcup_{\delta\in S(\alpha)}\open{\delta}\right),
\end{equation}
where
\begin{equation}
S(\alpha)=\left\{\beta\in\F(S)\colon\beta\neq\alpha\text{ and }\alpha\in\F(\beta)\right\}.
\end{equation}
Consequently, for all $\alpha\in\F_d(S)$ such that $\alpha\subset E''\cap F'$:
\begin{equation}
\H^d(\alpha\setminus F)\leq\H^d\left(\alpha\cap\bigcup_{\delta\in S(\alpha)}\psi_1(F\cap\open{\delta})\right).
\end{equation}
By summing over all $d$-dimensional faces of $S$ that are included in $E''\cap F'\cap\overline{W}$ and by~\eqref{equationpolyhedraloptimizationD} we get:
\begin{align}\notag
\H^d((E''\cap F')\setminus F)&=\sum_{\alpha\subset E''\cap F'}\H^d\left(\alpha\setminus F\right)\\\notag
&\leq\sum_{\alpha\subset E''\cap F'}\H^d\left(\alpha\cap\bigcup_{\delta\in S(\alpha)}\psi_1(F\cap\open{\delta})\right)\\\notag
&=\H^d\left(\bigcup_{\alpha\subset E''\cap F'}\bigcup_{\delta\in S(\alpha)}\alpha\cap\psi_1(F\cap\open{\delta})\right)\\\notag
&\leq\H^d\left(\bigcup_{\alpha\subset F'}\bigcup_{\delta\in\F(S),\dim(\delta)>d}\alpha\cap\psi_1(F\cap\open{\delta})\right)\\\notag
&\leq\H^d\left(\left(\bigcup_{\alpha\subset F'}\alpha\right)\cap\left(\bigcup_{\delta\in\F(S),\dim(\delta)>d}\psi_1(F\cap\open{\delta})\right)\right)\\\notag
&=\H^d\left(F'\cap\bigcup_{\delta\in\F(S),\dim(\delta)>d}\psi_1(F\cap\open{\delta})\right)\\\notag
&=\H^d\left(\bigcup_{\delta\in\F(S),\dim(\delta)>d}\psi_1(F\cap\open{\delta})\right)\\\notag
&\leq\sum_{\delta\in\F(S),\dim(\delta)>d}\H^d\left(\psi_1(F\cap\open{\delta})\right)\\\notag
&\leq\sum_{\delta\in\F(S),\dim(\delta)>d}K'\H^d\left(F\cap\open{\delta}\right)\\\notag
&=K'\H^d\left(\bigcup_{\delta\in\F(S),\dim(\delta)>d}F\cap\open{\delta}\right)\\\notag
&=K'\H^d(F\setminus\mathcal{U}(\F_d(S)))\\\label{equationpolyhedraloptimizationE}
&\leq K'\H^d(F\setminus E'').
\end{align}

To achieve proving that $E''$ is quasiminimal, let us split $E''\setminus F$:
\begin{equation}\label{equationpolyhedraloptimizationF}
E''\setminus F=(E''\setminus(F'\cup F))\sqcup((E''\cap F')\setminus F)\subset(E''\setminus F')\cup((E''\cap F')\setminus F).
\end{equation}
Using~\eqref{equationpolyhedraloptimizationF}, \eqref{equationpolyhedraloptimizationD} and~\eqref{equationpolyhedraloptimizationE} we obtain
\begin{equation}\label{equationpolyhedraloptimizationG}
\begin{split}
\H^d(E''\setminus F)&\leq\H^d(E''\setminus F')+\H^d((E''\cap F')\setminus F)\\
&\leq K'(M+1)\H^d(F\setminus E'')\\
&\leq MM'\H^d(F\setminus E''),
\end{split}
\end{equation}
where $M'=2K'$ depends only on $d$ and $n$. Using the fact that $E''\setminus F\subset\xi_{\phi_1}$ and $F\setminus E''\subset\phi_1(\xi_{\phi_1})$ we also have the following set equalities:
\begin{equation}\label{equationpolyhedraloptimizationH}
\begin{split}
E''\cap\xi_{\phi_1}&=((E''\setminus F)\cap\xi_{\phi_1})\sqcup(E''\cap F\cap\xi_{\phi_1})\\
&=(E''\setminus F)\sqcup(E''\cap F\cap\xi_{\phi_1}),
\end{split}
\end{equation}
and
\begin{equation}\label{equationpolyhedraloptimizationI}
\begin{split}
\phi_1(E''\cap\xi_{\phi_1})&=\phi_1(E'')\cap\left((\phi_1(\xi_{\phi_1})\setminus E'')\sqcup(\xi_{\phi_1}\cap E'')\right)\\
&=F\cap\left((\phi_1(\xi_{\phi_1})\setminus E'')\sqcup(E''\cap\xi_{\phi_1})\right)\\
&=((F\setminus E'')\cap\phi_1(\xi_{\phi_1}))\sqcup(F\cap E''\cap\xi_{\phi_1})\\
&=(F\setminus E'')\sqcup(F\cap E''\cap\xi_{\phi_1}).
\end{split}
\end{equation}
Using~\eqref{equationpolyhedraloptimizationG}, \eqref{equationpolyhedraloptimizationH} and~\eqref{equationpolyhedraloptimizationI} we finally get
\begin{equation}
\begin{split}
\H^d(E''\cap\xi_{\phi_1})&=\H^d(E''\setminus F)+\mathcal{H}^d(E''\cap F\cap\xi_{\phi_1})\\
&\leq MM'\H^d(F\setminus E'')+\H^d(E''\cap F\cap\xi_{\phi_1})\\
&\leq\max(MM',1)\left(\H^d(F\setminus E'')+\H^d(F\cap E''\cap\xi_{\phi_1})\right)\\
&=MM'\H^d(\phi_1(E''\cap\xi_{\phi_1})),
\end{split}
\end{equation}
which achieves proving that $E''$ is $(MM',\diam(W))$-quasiminimal over $W$ for $\H^d$.

Let us now verify the last point of the lemma. Suppose that $S'\subset\F(S)$ is a $d$-dimensional skeleton of $S$ and consider the following definition for $0<l\leq d$:
\begin{equation}
\begin{cases}
\F^{*}_d(S')=\F_d(S)\cap S'\\
\F^{*}_{l-1}(S')=\left\{\alpha\in\F_{l-1}(S)\cap S'\colon\forall l'\geq l,\forall\beta\in\F^*_{l'}(S'),\alpha\not\subset\beta\right\}.
\end{cases}
\end{equation}
The complexes $\F^{*}_l(S')$ are in fact the $l$-dimensional polyhedrons of $S'$ that are not a subface of any polyhedron of $S'$ with higher dimension. Let us also define
\begin{equation}
\begin{cases}
\displaystyle\S^d=\left\{T\subset\F(S)\colon(E\setminus\overline{W})\cup\U(T)\in\mathfrak{C}\text{ and }J_h^d(\U(T))=\min_{G\in\mathfrak{C}\cap\mathfrak{S}}J_h^d(G\cap\overline{W})\right\}\\
\displaystyle\S^{l-1}=\left\{T\in\S^l\colon J_h^{l-1}(\U(\F^{*}_{l-1}(T)))=\min_{T'\in\S^l}J_h^{l-1}(\U(\F^{*}_{l-1}(T')))\right\}.
\end{cases}
\end{equation}
Notice that $\S^d$ is not empty since the skeleton that we used to build $E''$ is in it, and by induction it is easy to check that $\S^0$ is not empty. Without changing the above proof we can assume that we took $E''=(E\setminus\overline{W})\cup(\mathcal{U}(S''))$ where $S''\in\mathcal{S}^0$. For $0\leq l\leq d$, put
\begin{equation}
S^l=\mathcal{F}^*_{l'}(S''),
\end{equation}
and use the same notations as in the last point of the lemma. We can use the same argument as we used above to prove the quasiminimality of $E''$, to prove that each $E^l$ is quasiminimal for $\H^l$ over $W^l$.
\end{proof}

Before stating and proving our main theorem we provide the following lemma which will allow us to consider minimization problems with respect to the integral functional $J_h^d$ instead of $\H^d$ only, and consider the case of almost-minimal sets as well with a gauge function closely related to $h$.

\begin{lemma}[Lower semicontinuity of $J_h^d$ with respect to $\H^d$]\label{lemmameaseuresemicontinuity}
Suppose that $U$ is an open domain, that $h\colon U\rightarrow[1,M]$ is lower semicontinuous and that $(E_k)_{k\in\mathbb{N}}$ is a sequence of measurable subsets of $U$.

If there is a measurable set $E\subset U$ such that for all open subset $V\subsubset U$:
\begin{equation}
\H^d(E\cap V)\leq\liminf_{k\rightarrow\infty}\H^d(E_k\cap V),
\end{equation}
then the following holds:
\begin{equation}
J_h^d(E)\leq\liminf_{k\rightarrow\infty}J_h^d(E_k).
\end{equation}
\end{lemma}

\begin{proof}
Fix an integer $m>0$ and for $l\geq 0$, put
\begin{equation}
X_l=\{x\in U:\;h(x)>2^{-m}l\}
\end{equation}
Notice that $X_l$ is open because $h$ is lower semicontinuous, and for $x\in U$ set
\begin{equation}\label{equationmeaseuresemicontinuityA}
h_m(x)=2^{-m}\sum_{l\geq 0}\mathds{1}_{X_l}(x),
\end{equation}
where $\mathds{1}_{X_l}$ stands for the characteristic set function of $X_l$. Since $h$ is bounded, the sum in~\eqref{equationmeaseuresemicontinuityA} is finite and for all open subset $V\subsubset U$:
\begin{equation}
\begin{split}
J_{h_m}^d(E\cap V)&=2^{-m}\sum_{l\geq 0}l\H^d(E\cap X_l\cap V)\\
&\leq2^{-m}\sum_{l\geq 0}l\limsup_{k\rightarrow\infty}\H^d(E_k\cap X_l\cap V)\\
&=\limsup_{k\rightarrow\infty}J_{h_m}^d(E_k\cap V).
\end{split}
\end{equation}
Besides, notice that
\begin{equation}
h\leq h_m\leq h+2^{-m},
\end{equation}
which gives
\begin{equation}\label{equationmeaseuresemicontinuityB}
J_h^d(E)\leq J_{h_m}^d(E)\leq\liminf_{k\rightarrow\infty}J_{h_m}^d(E_k)\leq\liminf_{k\rightarrow\infty}J_h^d(E_k)+2^{-m}H,
\end{equation}
where $H=\sup_k\H^d(E_k)$. Consider the two possibles cases:
\begin{enumerate}
\item if $H<\infty$ then by taking limits in~\eqref{equationmeaseuresemicontinuityB} we have finished;
\item if $H=\infty$, there are two more possible cases:
\begin{itemize}
\item if $\liminf_{k\rightarrow\infty}\H^d(E_k\cap V)=\infty$ we have finished;
\item otherwise, we can extract a subsequence $(E'_k)_{k\in\mathbb{N}}$ of $(E_k)$ such that $\sup_k\H^d(E'_k)<\infty$,
\begin{equation}
\lim_{k\rightarrow\infty}\H^d(E'_k\cap V)=\liminf_{k\rightarrow\infty}\H^d(E_k\cap V),
\end{equation}
and go back to the above case 1 by replacing $E_k$ with $E'_k$.\qedhere
\end{itemize}
\end{enumerate}
\end{proof}

We now have all the required ingredients to proceed into proving the main result. A large part of the argument is based upon the proof of the second point of theorem \ref{theoremdavid} (see~\cite{david1}). Our polyhedral optimization theorem~\ref{theorempolyhedraloptimization} allows us to build a polyhedric minimizing sequence for which we have to make sure that the subfaces of dimension less than $d$ do not converge towards a set of positive $d$-dimensional measure. This will be done using the optimality of subdimensional cores we obtained before. Notice that we do not require our minimizing sequence to be made of reduced sets, which might prove to be convenient when trying to control the topological constraint when taking limits, since the subdimensional cores can play a topological role.

\begin{theorem}[Main result]\label{theoremexistence}
Let $U$ be an open, bounded domain of $\mathbb{R}^n$, $0\leq d<n$, $\mathfrak{F}$ a non-empty family of relatively closed subsets of $U$ stable under the $\diam(U)$-deformations over $U$ such that $\inf_{F\in\mathfrak{F}}\H^d(F)<\infty$ and $h\colon U\rightarrow[1,M]$ such that
\begin{equation}\label{equationexistence}
\forall(x,y)\in U^2\colon h(y)\leq(1+\tilde h(\Vert x-y\Vert))h(x)
\end{equation}
where $\tilde h\colon]0,\diam(U)[\rightarrow[0,\infty]$ verifies
\begin{equation}
\lim_{r\rightarrow 0}\tilde h(r)=0.
\end{equation}

Then one can build a sequence $(E_k)_{k\in\mathbb{N}}$ of elements of $\mathfrak{F}$ satisfying the following properties:
\begin{itemize}
\item$E_k\Hconverge{U} E'$;
\item$\displaystyle J_h^d(E')\leq\inf_{F\in\mathfrak{F}}J_h^d(F)$;
\item$E'$ is almost-minimal with gauge function $\tilde h$ over $U$.
\end{itemize}

More precisely, by putting
\begin{equation}
E^l=\ker^l_d(E)\qquad\text{and}\qquad U^l=U\setminus\bigsqcup_{d\geq l'>l}E^{l'}
\end{equation}
for $0\leq l\leq d$, the following holds:
\begin{itemize}
\item $\ker^l_d(E_k)\Hconverge{U^l}E^l$;
\item $\displaystyle J_h^l(E^l)\leq\inf_{\substack{F\in\mathfrak{F}\\F\setminus U^l=E\setminus U^l}}J_h^l(\ker^l_d(F))$;
\item $E^l$ is almost-minimal with gauge function $\widetilde{h}$ over $U^l$.
\end{itemize}
\end{theorem}

Notice that we did not require that all the elements of $\mathfrak{F}$ have finite measure. However, we can always consider the subclass of $d$-sets of $\mathfrak{F}$ (which is not empty since $\inf_{F\in\mathfrak{F}}\H^d(F)<\infty$), which is stable under deformations on $U$ due to the Lipschitz condition in definition~\ref{definitiondeformation}.

\begin{proof}
Suppose that $(U_k)_{k\in\mathbb{N}}$ is an increasing sequence of open and relatively compact subsets of $U$ such that
\begin{equation}\label{equationexistenceA}
\bigcup_{k\in\mathbb{N}}U_k=U.
\end{equation}
For instance, one can take
\begin{equation}
U_k=\left\{x\in\open{B}(0,2^k)\colon B(x,2^{-k})\subset U\right\}.
\end{equation}
For $k\geq0$, set $\epsilon_k=2^{-k}$ and choose $E_k\in\mathfrak{F}$ such that:
\begin{equation}\label{equationexistenceB}
J_h^d(E_k\cap U_k)\leq\epsilon_k+\inf_{F\in\mathfrak{F}}J_h^d(F\cap U_k).
\end{equation}
Set $\eta>0$. By applying lemma~\ref{theorempolyhedraloptimization} to $E_k\cap U_k$ inside $U_k$, one can build an open set $W_k$ such that $U_k\subsubset W_k\subsubset U$, a $n$-dimensional complex $S_k$ and an Almgren competitor
\begin{equation}\label{equationexistenceC}
E'_k=(E_k\setminus U_k)\sqcup\bigsqcup_{0\leq l\leq d}E^l_k\in\mathfrak{F}
\end{equation}
such that $J_h^d(E'_k\cap U_k)\leq(1+\eta)J_h^d(E_k\cap U_k)$, and $E^l_k=\U(S^l)$ (for $0\leq l\leq d$) where $S^l\subset\F_l(S_k)$, where each $E^l_k$ is $(MM',\diam(U_k)$-quasiminimal over $U_k^l$ for $\H^l$. Notice that $J_h^d(E_k\cap U_k)\leq M\mathcal{H}^d(E_k\cap U_k)<+\infty$, because $E_k$ is a $d$-set included in $U$, which is bounded. By taking $\eta=\frac{\epsilon_k}{J_h^d(E_k\cap U_k)}>0$ and by~\eqref{equationexistenceB}, we get:
\begin{equation}
J_h^d(E'_k\cap U_k)\leq\epsilon_k+J_h^d(E\cap U_k)\leq2\epsilon_k+\inf_{F\in\mathfrak{F}}J_h^d(F\cap U_k).
\end{equation}

We can extract from $(E'_k$) a convergent subsequence that converges towards a relatively closed subset $E^l$ of $U$ locally on every compact of $U$. By setting
\begin{equation}\label{equationexistenceD}
U^l=U\setminus\bigcup_{d\geq l'>l}E^{l'}
\end{equation}
and extracting multiple subsequences, we can even assume that $E^l_k\cap U^l$ converges towards $E^l$ locally on every compact of $U^l$. To summarize, once we have extracted all our convergent subsequences, we obtain:
\begin{equation}\label{equationexistenceE}
\forall l\leq d\colon E_k^l\Hconverge{U^l}E^l\quad,\quad E'_k=\bigsqcup_{0\leq l\leq d}E^l_k\quad,\quad E'_k\Hconverge{U}E\quad\text{and}\quad E=\bigsqcup_{0\leq l\leq d}E^l.
\end{equation}

Now fix $l\leq d$, suppose that $V\subsubset U^l$ and for $\epsilon>0$, put:
\begin{equation}\label{equationexistenceF}
W_\epsilon=\bigcup_{x\in E^d\cup E^{d-1}\ldots\cup E^{l+1}}B(x,\epsilon).
\end{equation}
Since $E^{l'}_k\Hconverge{U^{l'}}E^{l'}$ when $l'>l$, one can find $k_0$ such that
\begin{equation}\label{equationexistenceG}
\forall k\geq k_0\colon\bigcup_{d\geq l'>l}E^{l'}_k\subset W_\epsilon.
\end{equation}
Besides, $V\subsubset U^l$ and $U^l\cap E^{l'}=\emptyset$ when $l'>l$, so we can take $\epsilon>0$ small enough such that $V\cap W_\epsilon=\emptyset$, which in turn gives:
\begin{equation}\label{equationexistenceH}
\forall k\geq k_0,\forall l'>l\colon E^{l'}_k\cap V=\emptyset.
\end{equation}
Since $E^l_k$ is $(MM',\diam(W_k))$-quasiminimal over $W^l_k$, it is also $(MM',\diam(U))$-quasiminimal over $V\cap W^l_k$ when $k\geq k_0$. Furthermore, since $\overline{V}$ is compact and included in $U$, which is covered by $\bigcup_kU_k$ by~\eqref{equationexistenceA}, by taking a finite covering we can assume that $k_0$ is large enough such that
\begin{equation}\label{equationexistenceI}
\forall k\geq k_0\colon\overline{V}\subset U_k.
\end{equation}

We can also assume that, for instance, $\overline{\R}(S_k)\leq\epsilon_k$ --- by taking $R$ small enough in lemma's~\ref{theorempolyhedraloptimization} proof. That way, again by taking $k_0$ large enough, we can assume that we can extract a subset $S'_k$ from $S_k$ verifying
\begin{equation}\label{equationexistenceJ}
V\subsubset\U(S'_k)\subsubset W^l_k.
\end{equation}
By extracting another subsequence, we can even suppose that for instance $\underline{\mathcal{R}}(S_{k+1})\geq8\sqrt{n}\overline{\mathcal{R}}(S_k)$, and extract our complexes $S'_k$ such that:
\begin{equation}\label{equationexistenceK}
\forall k\geq k_0\colon\U(S'_k)\subset\U(S'_{k+1}).
\end{equation}

Suppose that we have done all the above setup, put
\begin{equation}\label{equationexistenceL}
D=\H^l(\U(\F_l(S_{k_0})))\qquad\text{and}\qquad D'=\min_{\alpha\in\F_l\left(S'_{k_0}\right)}\H^l(\alpha),
\end{equation}
and suppose that $k\geq k_0$. Our next goal is to prove the two following statements:
\begin{gather}\label{equationexistenceM}
\H^l(E^l_k\cap V)\leq MM'D,\\\label{equationexistenceN}
\H^l(E^l_k\cap\overline{V})\in\{0\}\cup[\frac{D'}{M'},+\infty[.
\end{gather}
Firstly, put $W'_k=\open{\U(S'_k)}$. By applying lemma~\ref{lemmapolyhedraldeformation} to the closed $l$-set $E_k^l$ in the complex $S'_{k_0}$, we get a deformation $(\psi_t)$ over $W'_{k_0}$ such that, by putting $E_k'^l=\psi_1(E_k^l)$ we have $E_k'^l\cap W'_l\subset\U(\F_l(S'_k))$,
\begin{equation}
\H^l(E'^l_k)\leq D\qquad\text{and}\qquad\H^l(E'^l_k)\leq M'\H^l(E^l_k\cap W'_{k_0}).
\end{equation}
Using the quasiminimality of $E_k^l$ over $W_k^l$ we get directly
\begin{equation}
\H^l(E^l_k\cap W'_{k_0})\leq MM'\H^l(E'^l_k)\leq MM'D,
\end{equation}
and since $V\subset W'_{k_0}$ by~\eqref{equationexistenceJ}, we obtain~\eqref{equationexistenceM}.

Now, if we suppose that $\H^l(E^l_k\cap W'_{k_0})<\frac{D'}{M'}$ then by~\eqref{equationexistenceN} we have
\begin{equation}\label{equationexistenceO}
\H^l(E'^l_k)<D'=\min_{\alpha\in\F_l(S'_{k_0})}\H^l(\alpha)
\end{equation}
and since $E_k'^l\subset W'_{k_0}\subset\U(\F_l(S'_{k_0}))$ this means that for all $\alpha\in\F_l(S'_{k_0}$, $\alpha\cap E_k'^l\neq\alpha$. By using lemma~\ref{lemmapolyhedralerosion} we can build a deformation $(\psi'_t)$ over $W'_{k_0}$ such that:
\begin{equation}\label{equationexistenceP}
E''^l_k=\psi'_1(E'^l_k)\subset\U(\F_{l-1}(S'_{k_0})).
\end{equation}
Using the quasiminimality of $E_k^l$ again and by~\eqref{equationexistenceP}, we get
\begin{equation}
\H^l(E^l_k\cap W'_{k_0})\leq MM'\H^l(E''^l_k\cap W'_{k_0})=0,
\end{equation}
and since $\overline{V}\subset\open{W'_{k_0}}$ by~\eqref{equationexistenceJ}, we get~\eqref{equationexistenceN}.

Applying theorem~\ref{theoremdavid} to the sequence $(E_k^l\cap U)_{k\geq k_0}$ --- which converges towards $E^l\cap V$ --- gives us the following points:
\begin{itemize}
\item $E^l\cap V$ is $(MM',\diam(U))$-quasiminimal over $U$ and $\ker^l(E^l\cap V)E^l\cap V$;
\item\begin{equation}
\H^l(E^l\cap V)\leq\liminf_{k\rightarrow\infty}\H^l(E^l_k\cap V)\leq MM'D<\infty
\end{equation}
so $E^l\cap V$ is a closed relative $l$-subset of $V$;
\item\begin{equation}
\H^l(E^l\cap\overline{V})\geq C\limsup_{k\rightarrow\infty}\H^l(E^l_k\cap\overline{V})\in\{0\}\cup[\frac{CD'}{M'},+\infty[,
\end{equation}
and consequently two cases are possible:
\begin{itemize}
\item if $\H^l(E^l\cap V)=0$ then $\ker^l(E^l\cap V)=\emptyset$ and
\begin{equation}
\limsup_{k\rightarrow\infty}\H^l(E^l_k\cap V)=0,
\end{equation}
which means that for $k$ large enough, $E^l_k\cap V=\emptyset$ and $E^l\cap V=\emptyset=\ker^l(E^l\cap V)$ (since $E^l_k\cap V\Hconverge{V}E^l\cap V$);
\item if $\H^l(E^l\cap V)>0$ then for $k$ large enough we have $\H^l(E^l_k\cap V)\geq\frac{D'}{M'}>0$, so $E^l_k\cap V\neq\emptyset$.
\end{itemize}
\end{itemize}
We get from the second point that $\ker^{l'}(E^l)=\emptyset$ and $\H^d(E^{l'})=0$ if $l'>l$. If we take for $V$ a ball centered on $E^l$ and relatively compact in $U^l$ with arbitrary small radius, the third point tells us that $E^l_k\Hconverge{U^l}\ker^l(E^l)$ and that $\ker^l(E^l)=E^l$, and as a consequence:
\begin{equation}
\ker^l_d(E'_k)\Hconverge{U^l}\ker^l_d(E).
\end{equation}
The first point implies that $E^l$ is $(MM',\diam(U))$-quasiminimal over $U^l$. And using lemma~\ref{lemmameaseuresemicontinuity} we also get:
\begin{equation}
\begin{split}
J_h^d(E)&=J_h^d(E^d)\\
&\leq\liminf_{k\rightarrow\infty}J_h^d(E^d_k\cap U_k)\\
&=\liminf_{k\rightarrow\infty}\inf_{F\in\mathfrak{F}}J_h^d(F\cap U_k)\\
&\leq\inf_{F\in\mathfrak{F}}J_h^d(F).
\end{split}
\end{equation}
Notice that we could start again all the above process with an initial sequence $(E'_k)$ such that, for instance:
\begin{equation}\label{equationexistenceQ}
J_h^l(E'_k\cap V^l_{\epsilon_k})\leq 2\epsilon_k+\inf_{F\in\mathfrak{F},F\setminus V^l_{\epsilon_k}=E'_k\setminus V^l_{\epsilon_k}}J_h^l(F\cap V^l_{\epsilon_k}),
\end{equation}
where
\begin{equation}
V^l_{\epsilon_k}=U_k\setminus\left(\bigcup_{\substack{l+1\leq l'\leq d\\x\in E^{l'}}}B(x,\epsilon_k)\right).
\end{equation}
That way we would ensure that
\begin{equation}
J_h^l(\ker^l_d(E))\leq\inf_{\substack{F\in\mathfrak{F}\\F\setminus U^l=E\setminus U^l}}J_h^l(\ker^l_d(F)),
\end{equation}
and by defining new limits $E^{l'}$ for $l'\leq l$ and a new $V_{\epsilon_k}^l$, by induction we could prove the last point of theorem~\ref{theoremexistence}.

At this point, all that is left to prove is the almost-minimality of the $E^l$: suppose that $\delta>0$ and that $(f_t)$ is a $\delta$-deformation over $U$. Let us apply the last point of theorem~\ref{theoremdavid} to the sequence $(E'_k)$: we get a Lipschitz map $g$ over $U$ and using proposal~\ref{proposalautomaticdeformation} we can build a $\delta$-deformation $(g_t)$ over $U$, with $g_1$ verifying equation~\eqref{equationdavidA}.

Suppose that $M=1$. For $k$ large enough we have $\xi_{f_1}\cup \xi_{g_1}=\xi_{f_1}\subset U_k$ and we can even suppose that $\support(g)\subsubset U_k$. By~\eqref{equationexistenceQ} and since $g_1(E_k)\in\mathfrak{F}$ we automatically have
\begin{equation}
\H^d(g_1(E'_k\cap\xi_{g_1}))\geq\H^d(E'_k\cap\xi_{g_1})-2\epsilon_k.
\end{equation}
By~\eqref{equationdavidA} and provided that $k$ is large enough we also get
\begin{equation}\label{equationexistenceS}
\H^d(f_1(E\cap\xi_{f_1}))\geq\H^d(E\cap\xi_{f_1})-4\epsilon_k,
\end{equation}
and since $\epsilon_k\rightarrow 0$: 
\begin{equation}
\H^d(f_1(E\cap\xi_{f_1}))\geq\H^d(E\cap\xi_{f_1}),
\end{equation}
which achieves proving that $E$ is minimal over $U$.

When $M>1$, we can find a ball $B$ with radius $\delta$ such that $\xi_{g_1}\subset\xi_{f_1}\subset B$. Again, by~\eqref{equationexistenceQ} we have
\begin{equation}
J_h^d(g_1(E'_k\cap\xi_{g_1}))\geq J_h^d(E'_k\cap\xi_{g_1})-2\epsilon_k,
\end{equation}
and by~\eqref{equationexistence}:
\begin{equation}
\H^d(g_1(E'_k\cap\xi_{g_1}))\geq(1+\widetilde{h}(\delta))\H^d(E'_k\cap\xi_{g_1})-2\epsilon_k.
\end{equation}
Using~\eqref{equationdavidA} again, we obtain
\begin{equation}
\H^d(f_1(E\cap\xi_{f_1}))\geq(1+\widetilde{h}(\delta))(\H^d(E\cap\xi_{f_1})-2\epsilon_k)-2\epsilon_k,
\end{equation}
which similarly gives in turn
\begin{equation}
\H^d(f_1(E\cap\xi_{f_1}))\geq(1+\widetilde{h}(\delta))\H^d(E\cap\xi_{f_1}).
\end{equation}

The above argument we used to prove the almost-minimality of $E$ could be done again in decreasing dimension for $E^l$ inside $U^l$, which achieves proving the last point of theorem~\ref{theoremexistence}.
\end{proof}

\subsection{Two examples of application}\label{subsectionapplication}

As we outlined before, we cannot ensure that the minimal candidate given by theorem~\ref{theoremexistence} is still in our topological class $\mathfrak{F}$. More precisely, it is easy to find cases for which there is not even a solution to our measure minimization problem in $\mathfrak{F}$ --- and even more since we supposed that $U$ is open. For instance, when $n=2$ and $d=1$, take $U=]-2,2[^2\setminus[-1,1]^2$ and consider the class $\mathfrak{F}$ of paths joining $x=(1,-2)$ to $y=(1,2)$ with open extremities, and included in $U$. Clearly, $\mathfrak{F}$ is stable under the $\diam(U)$-deformations over $U$ and it is easy to check that $\inf_{F\in\mathfrak{F}}\H^1(F)=4$ but every element of $\mathfrak{F}$ is of length greater than $4$ since the open line segment joining $x$ to $y$ is not in $\mathfrak{F}$. Notice that in that case, the minimal candidate given by theorem~\ref{theoremexistence} is in fact the union of the two open line segments joining $x$ to $(1,-1)$ and $y$ to $(1,1)$. The convergence notion ``over all compact set of $U$'' we had to use because $U$ is open is rather weak near the boundary of $U$ and for this reason we can cause gap to appear in $E$ when taking limit in our minimizing sequence, as in the previous example.

However, in the context of a more restrictive notion of minimal sets than Almgren-minimal sets (see below) our approach can give complete existence results. The definition of this other kind of minimal sets is borrowed from Guy David in~\cite{david3}, where the reader might find more details about how they can be useful to study the regularity of minimal segmentations for the Mumford-Shah functional.

Let $E$ be a closed set in $\mathbb{R}^n$. A Mumford-Shah competitor for $E$ (a ``MS-competitor'' in short) is a closed set $F$ such that we can find a closed ball $B$ verifying
\begin{equation}
F\setminus B=E\setminus B
\end{equation}
and for all $x,y\in\mathbb{R}^n\setminus (B\cup E)$, ``$F$ separates $x$ from $y$ whenever $E$ does'' (i.e. if $x$ and $y$ lie in different connex components of $\mathbb{R}^n\setminus E$ then they also lie in different connex components of $\mathbb{R}^n\setminus F$). We say that $E$ is MS-minimal if
\begin{equation}\label{equationMScompetitor}
\H^{n-1}(E\setminus F)\leq\H^{n-1}(F\setminus E)
\end{equation}
for all MS-competitor $F$ of $E$.

The following statement can be used to find MS-minimizers inside a localized class of MS-competitors. In fact, we have to give an upper bound on the size of the ball in which we allow our sets to be changed. Also, we have to give some way to ensure that our minimizing limit will not come too close to the boundary of the ball when taking limit in our minimizing sequence, to avoid gaps to appear as we explained above.

\begin{corollary}[MS-minimal competitor inside a ball]\label{corollaryapplication1}
Suppose that $n\geq 1$ and that $E$ is a closed set. For all ball $B$, $E$ has a MS-minimal competitor $E'$ inside $\overline{B}$ (i.e. $E'$ is minimal like in~\eqref{equationMScompetitor} amongst all MS-competitors $F$ of $E$ such that $F\setminus\overline{B}=E\setminus\overline{B}$).
\end{corollary}

The statement still holds when $B$ is any compact convex set, although we will prove it only in the case when $B$ is a ball. However, the proof may be adapted easily for this case.

\begin{proof}
For convenience, let us suppose that $B$ is open and centered at the origin, denote by $\pi$ the radial projection onto $\partial B$ centered at the origin and set $B'=2B$. For $F\subset\mathbb{R}^n$ we also define the two set functions
\begin{equation}
\begin{split}
H(F)&=(F\cap B)\cup\{tx\colon x\in F\cap\partial B\text{ and }1\leq t\leq 2\}\cup\{2x\colon x\in F\setminus B\}\\
I(F)&=(F\cap B)\cup\left\{\frac{x}{2}\colon x\in F\setminus B'\right\}.
\end{split}
\end{equation}
These functions will be used to turn $E$ into a cone inside $B'\setminus B$ and to easily build deformations on $B'\setminus H(E)$. Notice that $I\circ H(F)=F$ and that $F$ is a MS-competitor of $F'$ if and only if $H(F)$ is a MS-competitor of $H(F')$.

For $R>0$ small enough, one can build a dyadic complex $S$ inside $B'$ such that $\overline{B}\subset\U(S)\subset B'$. In fact, by fixing an orthonormal basis with origin at $(R/2,\ldots,R/2)$ and by taking all possible cubes inside $B'$ we can even assume that for all $x\in B'\setminus\U(S)$, the line segment $[0,x]$ intersects $\U(\F_\partial(S))=\partial\U(S)$ at an unique point $y$. In that case the map $f\colon x\mapsto y$ is Lipschitz, possibly with a very large constant depending on $R$.

Fix $A\geq 1$, set $U=B'\setminus(H(E)\setminus B)$, and define $h\colon U\rightarrow[1,A]$ by
\begin{equation}
h(x)=\begin{cases}1\text{ if $x\in\U(S)$,}\\A\text{ otherwise.}\end{cases}
\end{equation}
We also consider the class $\mathfrak{E}$ of relatively closed subsets $F$ of $\overline{B}\cap U$ such that $F\cup (H(E)\setminus\overline{B})$ are MS-competitors of $H(E)$. Notice that $\mathfrak{E}$ is not empty since $H(E)\cap U\in\mathfrak{E}$ and that $\inf_{F\in\mathfrak{E}}\H^{n-1}(F)<\infty$ since $\partial B\cap U\in\mathfrak{E}$. We also denote by $\mathfrak{F}$ the class of deformations over $U$ of the elements of $\mathfrak{E}$, which are also MS-competitors of $H(E)$ (see~\cite{dugundji,david2}).

Our function $h$ is only lower semicontinuous over $U$ although theorem~\ref{theoremexistence} requires $h$ to be continuous over $U$. However, if we consider how we proved theorem~\ref{theorempolyhedraloptimization} back then we can always suppose that we did a covering of $E\setminus\partial\U(S)$ by balls included in $U\setminus\partial\U(S)$, and assume that we built our global dyadic grid (the one we used to merge all the grids in the balls of our almost covering together) such that its faces cover $\partial\U(S)\cap W_k$ --- in fact, we have to consider a dyadic complex in the same basis as $S$ whose stride divides the stride of $S$. In that case, the upper semicontinuity of $h$ is not needed anymore, since it is only used when doing our magnetic projections to locally flatten $E$ onto its tangent planes, and $E\cap\partial\U(S)$ is already flattened onto the faces of our polyhedric grid.

With that minor modification we can therefore apply theorem~\ref{theoremexistence} to $\mathfrak{F}$, $h$ and $U$: we get a measure-minimizing sequence $(E_k)$ of elements of $\mathfrak{F}$ such that $E_k\cap\open{W_k}=\U(S'_k)\cap\open{W_k}$ (where $S'_k$ is an optimal $n-1$-dimensional skeleton of a complex $S_k$ with $W_k=\U(S_k)$ and $\U(S)\subset\overline{W_k}$),
\begin{equation}\label{equationapplication1A}
E_k\Hconverge{U} E'\qquad\text{and}\qquad J_h^{n-1}(E')\leq\inf_{F\in\mathfrak{F}}J_h^{n-1}(F).
\end{equation}
Fix $k>0$. Since $E_k\setminus\U(S)\subset W_k\subsubset U$, by using Kirszbraun theorem with $f$ and since $U\setminus B$ is a cone we can build a $\diam(W_k)$-deformation $(\phi_t)$ over $W_k$ such that $\phi_1\vert_{E_k\setminus\U(S)}=f\vert_{E_k\setminus\U(S)}$ and $\phi_1\vert_{\U(S)}=\identity_{\U(S)}$. Since $E_k\setminus\overline{B}\subsubset U$ we can even suppose (by taking $S_k$ large enough) that $E_k\setminus\U(S)\subset W_k$. Therefore, if we denote by $\alpha$ the Lipschitz constant of $f$ we get:
\begin{equation}
\H^{n-1}(\phi_1(E_k\setminus\U(S)))\leq\alpha^{n-1}\H^{n-1}(E_k\setminus\U(S)).
\end{equation}
Using our polyedric deformation lemmas~\ref{lemmapolyhedraldeformation} and \ref{lemmapolyhedralerosion} with $\phi_1(E_k)$ we can build a deformation $(\psi_t)$ over $U$ such that
\begin{equation}
\psi_1\circ\phi_1(E_k)\subset\U(S)
\end{equation}
and $\psi_1\circ\phi_1(E_k)$ is polyhedric inside $\U(S)$ (i.e. it is a finite union of subfaces of dimension at most $n-1$ of $S_k$). However, since $E_k$ was already polyhedric inside $\U(S)$ we also have
\begin{equation}
\psi_1\vert_{E_k\cap\U(S)}=\identity_{E_k\cap\U(S)}
\end{equation}
and
\begin{equation}
\begin{split}
\H^{n-1}(\psi_1\circ\phi_1(E_k\setminus\U(S)))&\leq C\H^{n-1}(\phi_1(E_k\setminus\U(S)))\\
&\leq C\alpha^{n-1}\H^{n-1}(E_k\setminus\U(S)),
\end{split}
\end{equation}
with $C$ depending only on $n$. Therefore, we get
\begin{equation}
J_h^{n-1}(\psi_1\circ\phi_1(E_k))\leq\frac{C\alpha^{n-1}}{A}J_h^{n-1}(E_k\setminus\U(S))+J_h^{n-1}(E_k\cap\U(S)).
\end{equation}

If we suppose that we took $A>C\alpha^{n-1}$ then necessarily $\H^{n-1}(E_k\setminus\U(S))=0$ since $E_k$ is optimal amongst all its polyhedric deformations. Notice that this argument also applies for the $n-2$-dimensional measure (since $\ker^{n-2}_{n-1}(E_k)$ is also optimal in theorem~\ref{theorempolyhedraloptimization}), and so on till dimension $0$. Therefore, this proves that $E_k\setminus\U(S)=\emptyset$, which means that $E_k$ never gets too close to $\partial B'$ and that $E'\subset\U(S)$.

We are now ready to show that $E'$ is a MS-competitor of $H(E)$. For that purpose, suppose that $x,y\in\mathbb{R}^n\setminus(B'\cup H(E))$ are separated by $H(E)$, pick a path $\gamma$ from $x$ to $y$ and let us show that $\gamma$ intersects $E'\cup(H(E)\setminus U)$. Since $E_k\in\mathfrak{F}$, $\gamma$ intersects $E_k\cup(H(E)\setminus U)$ at some point $x_k$ and by compacity of $\gamma$ we can find $x\in\gamma$ and extract a subsequence such that $\lim_{k\rightarrow\infty}x_k=x$. Also, notice that either $x_k\in H(E)\setminus U$ or $x_k\in E_k\subset\U(S)$ for all $k$ and therefore $x\in(H(E)\setminus U)\cup\U(S)$. If $x\in H(E)\setminus U$ we have finished. If $x\in\U(S)\setminus(H(E)\setminus U)$, then for $k_0$ large enough and $k\geq k_0$ we have $\overline{B}(x,\Vert x-x_k\Vert)\subset U$. Since $E_k\Hconverge{U} E'$ we can find a sequence $y_k$ of points of $E'\cap B(x,\Vert x-x_{k_0}\Vert)$ that converges towards $x$, and since $E'\cap\overline{B}(x,\Vert x-x_{k_0}\Vert)$ is closed this is enough to prove that $x\in E'$.

To conclude, let us denote by $\pi$ the radial projection onto $\partial B$ centered at the origin, and for $x\in\U(S)$ put
\begin{equation}
g(x)=\begin{cases}\pi(x)\text{ if $x\notin B$,}\\x\text{ otherwise.}\end{cases}
\end{equation}
Again by applying Kirszbraun theorem we can build a $\diam(U)$-deformation $(\phi_t)$ over $B'$ such that $\phi_1(H(E)\cap B)=H(E)\cap B$ and $\phi_1\vert_{E'}=g\vert_{E'}$. Notice that $g$ is $1$-Lipschitz and therefore $g(E')\in\mathfrak{F}$ with $\H^{n-1}(g(E'))\leq\H^{n-1}(E')$. Put 
\begin{equation}
E''=I(g(E')\cup (H(E)\setminus U))
\end{equation}
and notice that by~\eqref{equationapplication1A}, $E''$ is a MS-competitor of $E$ that meets the following requirements:
\begin{equation}
E''\setminus\overline{B}=E\setminus\overline{B}\qquad\text{and}\qquad\H^{n-1}(E''\cap\overline{B})\leq\inf_{\substack{F\text{ MS-competitor of }E\\F\setminus\overline{B}=E\setminus\overline{B}}}\H^{n-1}(F\cap\overline{B}).
\end{equation}
\end{proof}

Let us give another simple example of problem for which we do not need to control the topology near the boundary of the domain. In what follows we place ourself in the periodized cube $(\mathbb T^n,\dist)$ where 
\begin{equation}
\mathbb T^n=\mathbb R^n/\mathbb Z^n,
\end{equation}
$\dist$ is the natural induced distance
\begin{equation}
\dist(x,y)=\min_{z\in\mathbb{Z}^n}\Vert\tilde x-\tilde y+z\Vert
\end{equation}
and $\tilde x$ and $\tilde y$ denote equivalent points of $x$ and $y$ in $\mathbb R^n$. We say that a one-parameter family $(\phi_t)_{0\leq t\leq 1}$ of maps from $\mathbb T^n$ onto itself is a periodic deformation if $\phi_0=\identity_{\mathbb T^n}$, $(t,x)\mapsto\phi_t(x)$ is continuous over $[0,1]\times\mathbb T^n$ and $\phi_1$ is Hölder-regular, that is
\begin{equation}
\forall x,y\in\mathbb T^n\colon\dist(\phi_1(x)-\phi_1(y))<C\dist(x-y)^{1-\alpha}
\end{equation}
for some $C>0$ and $\alpha\in[0,1[$. For convenience, we will also denote by $\H^d$ the Hausdorff measure on $\mathbb T^n$ and keep the same definition for $J_h^d$ as in~\eqref{equationfunctionnal}.

\begin{corollary}[Periodic minimizer]\label{corollaryapplication2}
Suppose that $n>2$, $M>0$, $R>0$, that $h\colon\mathbb T^n\rightarrow[1,M]$ is a continuous function such that
\begin{equation}\label{equationapplication2A}
\forall x,y\in\mathbb T^n\colon\vert h(x)-h(y)\vert\leq\bar h(\dist(x,y))\quad\text{and}\quad\int_0^R\frac{\bar h(r)}{r}dr<\infty
\end{equation}
and that $\mathfrak{F}$ is a non-empty class of closed sets in $\mathbb{T}$ stable under periodic deformations. Then $\mathfrak{F}$ contains a set $E$ such that
\begin{equation}
J_h^2(E)=\inf_{F\in\mathfrak{F}}J_h^2(F).
\end{equation}
\end{corollary}

In fact, as we will see in the proof below, we could give a slightly more general result when $n=3$ by assuming that $\mathfrak F$ is stable by deformations $(\phi_t)$ over $\mathbb T^n$ such that $\phi_1$ is Lipschitz.

In what follows, for $F\subset\mathbb T^n$ we denote by $\tilde F$ its natural periodized equivalent set in $\mathbb R^n$ which verifies
\begin{equation}
\forall z\in\mathbb Z^n\colon\tilde F=z+\tilde F,
\end{equation}
and by $\tilde h$ the periodized equivalent of $h$ such that
\begin{equation}
\forall z\in\mathbb T^n\colon\tilde h(\tilde z)=h(z)
\end{equation}
for any equivalent point $\tilde z$ of $z$ in $\mathbb R^n$. It is easy to check that
\begin{equation}\label{equationapplication2B}
\forall F\subset\mathbb T^n,\forall z\in\mathbb R^n\colon J_h^d(F)=J_{\tilde h}^d\left(\tilde F\cap(z+[0,1[^n)\right).
\end{equation}
Similarly, for a given periodic deformation $(\phi_t)$ over $\mathbb T^n$ we denote by $(\tilde\phi_t)$ its periodized equivalent over $\mathbb R^n$ such that
\begin{equation}
\forall t\in[0,1],\forall x\in\mathbb R^n,\forall z\in\mathbb Z^n\colon\tilde\phi_t(x+z)=\tilde\phi_t(x)+z.
\end{equation}
These notations will be used to show that the optimization process described in theorems~\ref{theorempolyhedraloptimization} and~\ref{theoremexistence} can also be adapted to this periodic setup. However, the reader that is already convinced of that fact may skip the first part of the proof till~\eqref{equationapplication2F}.

\begin{proof}
Suppose that $\inf_{F\in\mathfrak F}\H^2(F)<\infty$ (otherwise our problem does not make sense and we have finished) and that $(E_k)_{k\geq 0}$ is a minimizing sequence of elements of $\mathfrak F$ with finite measure:
\begin{equation}
\lim_{k\rightarrow\infty}J_h^2(E_k)=\inf_{F\in\mathfrak{F}}J_h^2(F)<\infty.
\end{equation}

Fix $k\geq 0$. If we consider a global dyadic complex (used to merge the complexes of the almost-covering together) such that its $n-1$-dimensional subfaces cover $\partial[0,1]^n$, we can apply theorem~\ref{theorempolyhedralapproximation} and lemma~\ref{lemmapolyhedralerosion} to $\tilde E_k\cap[0,1]^n$. We get a deformation $(\phi_t)$ over $[0,1]^n$ and a complex $S_k$ such that $\U(S_k)=[0,1]^n$, $\phi_1(\tilde E_k\cap[0,1]^n)\subset\U(\F_2(S_k))$ and
\begin{equation}\label{equationapplication2C}
J_{\tilde h}^2(\phi_1(\tilde E_k\cap[0,1[^n))\leq(1+2^{-k})J_{\tilde h}^2(\tilde E_k\cap[0,1[^n).
\end{equation}
By using polyhedrons that are small enough, we can even suppose that
\begin{equation}\label{equationapplication2D}
\overline\R(S_k)\leq\frac1{100}.
\end{equation}
Since 
\begin{equation}
\forall t\in[0,1]\colon\phi_t\vert_{\partial[0,1]^n}=\identity_{\partial[0,1]^n},
\end{equation}
we can also extend $(\phi_t)$ as a periodized deformation $(\tilde\phi_t)$ over $\mathbb{R}^n$. Set
\begin{equation}
\tilde E'_k=\tilde\phi_1(\tilde E_k)
\end{equation}
and notice that the corresponding sequence $(E'_k)$ in $\mathbb T^n$ is a minimizing sequence of elements of $\mathfrak F$ for $J_h^2$. As we did before to prove theorem~\ref{theorempolyhedraloptimization}, we can do a finite minimization of $J_{\tilde h}^2(\psi_1(\tilde E'_k)\cap[0,1[)$ amongst the deformations over $\mathbb R^n$ such that $\psi_1(\tilde E'_k)$ is carried by a $2$-dimensional skeleton of $S_k$, $\support\psi\subset[0,1]^n$ and 
\begin{multline}
\forall t\in[0,1],\forall k\in\{1,\ldots,n\},\forall (z_1,\ldots,z_{n-1})\in[0,1]^{n-1}\colon\\
\psi_t(z_1,\ldots,z_l,0,z_{l+1},\ldots,z_n)=\psi_t(z_1,\ldots,z_l,1,z_{l+1},\ldots,z_n).
\end{multline}
Let us call $(\tilde\psi_t)$ an optimal deformation after having periodized it over $\mathbb R^n$ and put
\begin{equation}
\tilde E''_k=\tilde\psi_1(\tilde E_k).
\end{equation}
Notice that again, the corresponding set $E''_k\subset\mathbb T^n$ is in $\mathfrak F$. We also consider the infinite periodized complex $\tilde S_k$ defined by
\begin{equation}
\tilde S_k=\{z+\delta\colon z\in\mathbb Z\text{ and }\delta\in S_k\},
\end{equation}
and for $z\in\mathbb R^n$ and $r>0$, denote by $\Delta(z,r)$ the cube defined by
\begin{equation}
\Delta(z,r)=z+\left[-\frac r2,\frac r2\right[^n.
\end{equation}
By~\eqref{equationapplication2D}, for all $z\in\mathbb R^n$ one can find two finite subsets $T(z)$ and $T'(z)$ of $\tilde S_k$ such that
\begin{equation}\label{equationapplication2E}
B\left(z,\frac17\right)\subset\U(T(z))\subset B\left(z,\frac27\right)\subset\U(T'(z))\subset B\left(z,\frac37\right)\subsubset\Delta(z,1).
\end{equation}
Notice that $\tilde E''_k\cap\Delta(z,1)$ is also optimal amongst all its polyedric Almgren-competitors (ie. amongst all its images by a deformation with support in $\Delta(z,1)$ that are carried by a $2$-dimensional skeleton of $\tilde S$). Now suppose that $(f_t)_{0\leq t\leq 1}$ is a $\frac1{15}$-deformation over $\mathbb R^n$ and let $z\in\mathbb R^n$ such that $\support f\subset B\left(z,\frac17\right)$. By~\eqref{equationapplication2E}, the polyedral optimality of $\tilde E''_k$ and a similar argument as in theorem~\ref{theorempolyhedraloptimization} we get that
\begin{equation}
\H^2(\tilde E''_k\cap\xi_{f_1})\leq MM'\H^d(f_1(\tilde E''_k\cap\xi_{f_1}))
\end{equation}
where $M'$ depends only on $n$. Therefore, $(\tilde E''_k)$ is a sequence of quasiminimal sets with uniform constants and by~\eqref{equationapplication2C} and~\eqref{equationapplication2B}, $(E''_k)$ is a sequence of elements of $\mathfrak F$ minimizing $J_h^2$ for which the Hausdorff measure is lower semicontinuous. If we extract a convergent subsequence for the local Hausdorff convergence on every compact set of $\mathbb R^n$ --- the limit will also be periodized --- and by a similar argument as in theorem~\ref{theoremexistence}, we get that
\begin{equation}\label{equationapplication2F}
\tilde E''_k\Hconverge U\tilde E\qquad\text{and}\qquad J_h^d(E)\leq\liminf_{k\rightarrow\infty}J_h^d(E''_k)=\inf_{F\in\mathfrak F}J_h^2(F).
\end{equation}
Furthermore,
\begin{equation}
\tilde E=\tilde E^2\sqcup\tilde E^1\sqcup\tilde E^0
\end{equation}
where $\tilde E^l=\ker_2^l\tilde E$ is a reduced almost-minimal set over
\begin{equation}
\tilde U^l=\mathbb R^n\setminus\left(\bigcup_{2\geq k>l}\tilde E^k\right)
\end{equation}
with gauge function $\bar h$.

Let us now show that for $l=0,1,2$ we can build a deformation $(\phi^l_t)$ over $\mathbb T^n$ that sends an open neighborhood $W^l$ of $E^l$ onto $E^l$. For that purpose, we will use the biHölder equivalence of one- and two-dimensional reduced almost-minimal sets with one- and two-dimensional reduced minimal cones (see~\cite{taylor} for a biLipschitz version when $(l,n)=(2,3)$ or \cite{morgan} when $l=1$ with a slightly different requirement on $\bar h$, and see~\cite{david2,david3} for the biHölder regularity we will actually be using below). For $l=1,2$, denote by $\mathcal Z^l$ the set of reduced $l$-dimensional minimal cones over $\mathbb{R}^n$ (we will give a better description of $\mathcal Z^l$ later). Fix $\tau\in]0,1[$. By~\eqref{equationapplication2A}, Proposition~12.6 and Theorem~15.5 in \cite{david2}, for all $x\in\tilde E^l$ there is $r\in\left]0,\frac12-\tau\right[$, a cone $Z\in\mathcal Z^l$ centered at $x$ and a biHölder map $f\colon B(x,2r)\rightarrow\mathbb R^n$ such that:
\begin{multline}\label{equationapplication2G}
\forall y,z\in B(x,2r)\colon\\(1-\tau)\Vert z-y\Vert^{1+\tau}\leq\Vert f(z)-f(y)\Vert\leq(1+\tau)\Vert z-y\Vert^{1-\tau},
\end{multline}
\begin{gather}\label{equationapplication2H}
B(x,r(2+2\tau))\subset\tilde U^l,\\\label{equationapplication2I}
B(x,r(2-\tau))\subset f(B(x,2r))\quad\text{and}\quad\Vert f-\identity_{B(x,2r)}\Vert_\infty\leq r\tau\\
\intertext{and}\label{equationapplication2J}
\tilde E^l\cap B(x,r(2-\tau))\subset f(Z\cap B(x,2r))\subset E.
\end{gather}
Additionally, suppose that there is $C>0$, an open set $U_x$ and a map $g\colon U_x\rightarrow B(x,2r)$ such that
\begin{gather}\label{equationapplication2K}
B(x,2r)\cap Z\subset U_x\subset B(x,2r),\\
g(U_x)\subset Z\cap B(x,2r),\\\label{equationapplication2L}
\forall z\in Z\cap B(x,2r)\colon g(z)=z
\intertext{and}\label{equationapplication2M}
\forall y,z\in U_x\colon\Vert g(z)-g(y)\Vert\leq C\Vert z-y\Vert.
\end{gather}
We will explain in the last part of the proof how we can obtain such a Lipschitz map. Put
\begin{equation}\label{equationapplication2N}
V_x=f^{-1}(U_x)
\end{equation}
and for all $z\in V_x$, set
\begin{equation}
\pi_x(z)=f\circ g\circ f^{-1}(z).
\end{equation}
Notice that this definition is consistent because of~\eqref{equationapplication2I}, and that $V_x$ is an open set containing $\tilde E^l\cap B(x,2(r-\tau))$. Also, notice that
\begin{equation}
\forall z\in\tilde E^l\cap B(x,r(2-\tau))\colon\pi_x(z)=z
\end{equation}
by~\eqref{equationapplication2L}, that
\begin{equation}
\pi_x(V_x)\subset Z\cap B(x,r(2+\tau))
\end{equation}
by~\eqref{equationapplication2I}, \eqref{equationapplication2K} and~\eqref{equationapplication2J}, and that
\begin{multline}\label{equationapplication2O}
\forall y,z\in V_x\colon\\\Vert\pi_x(z)-\pi_x(y)\Vert\leq C\frac{1+\tau}{(1-\tau)^{\frac1{1+\tau}}}\Vert z-y\Vert^{\frac{1-\tau}{1+\tau}}=C'\Vert z-y\Vert^{1-\tau'}
\end{multline}
by~\eqref{equationapplication2G} and~\eqref{equationapplication2M}. Since we supposed that $r\in\left]0,\frac12-\tau\right[$ then $V_x\subsubset\Delta(x,1)$ and by using Mickle's extension theorem \cite{mickle}, we can extend $\pi_x$ over $\Delta(x,1)$ such that it stills verifies~\eqref{equationapplication2O} (possibly with a larger constant $C'$ and by taking a smaller set for $V_x$) and
\begin{equation}
\pi_x\vert_{\Delta(x,1)\setminus B(x,r(2-2\tau))}=\identity_{\Delta(x,1)\setminus B(x,r(2-2\tau))}.
\end{equation}
Therefore, we can consider the equivalent of $\pi_x$ inside $\mathbb T^n$ --- which we will also denote by $\pi_x$ for convenience, as well as $V_x$.

Since $E^2$ is compact in $\mathbb T^n$, and $\{V_x\colon x\in E^2\}$ is a covering of $E^2$ we can extract a finite covering $\{V_{x_1},\ldots,V_{x_p}\}$. Put
\begin{gather}
V^2=\bigcup_{1\leq i\leq p}V_{x_i}
\intertext{and}
\forall (t,z)\in[0,1]\times\mathbb T^N\colon\phi^2_t(z)=(1-t)z+t\pi_{x_1}\circ\pi_{x_2}...\circ\pi_{x_p}(z).
\end{gather}
Then, $(\phi_t^2)$ is a periodic deformation over $\mathbb T^n$ such that
\begin{equation}
E^2\subset V^2\qquad\text{and}\qquad\phi_1^2(V^2)=E^2.
\end{equation}

Similarly, the set $E^1\setminus V^2$ is compact, covered by $\{V_x\colon x\in E^1\setminus V^2\}$ and we can build an open set $V^1$ and a periodic deformation $(\phi_t^1)$ such that
\begin{equation}
E^1\setminus V^2\subset V^1\qquad\text{and}\qquad E^1\setminus V^2\subset\phi_1^1(V^1)\subset E^1.
\end{equation}
Additionally, by~\eqref{equationapplication2H} we have $B(x,r(2+2\tau))\cap E^2=\emptyset$ and we can also suppose that
\begin{equation}
\support\phi^1\cap E^2=\emptyset.
\end{equation}

Since $\H^0(\tilde E^0\cap[0,1[^n)<\infty$, $E^0$ is finite. Therefore, we can easily build an open set $V^0$ and a periodic deformation $(\phi^0_t)$ such that
\begin{gather}
\phi_1^0(V^0)=E^0\quad\text{and}\quad(\support\phi^0)\cap(E^2\cup E^1)=\emptyset.
\end{gather}

To conclude, put
\begin{equation}
V=V^1\cup V^2\cup V^2\quad\text{and}\quad\phi_t(z)=\phi_t^2\circ\phi_t^1\circ\phi_t^0(z)
\end{equation}
and notice that by construction, $\phi$ is a periodic deformation such that
\begin{equation}
\phi_1(V)=\phi_1(E)=E'\quad\text{and}\quad\phi_1(E^2)=E^2.
\end{equation}
Since $\phi_1$ is Hölder, then $\ker^2(E')=E^2$ and $J_h^2(E')=J_h^2(E)$. Recall that
\begin{equation}
E''_k\Hconverge{\mathbb T^n}E'\subset V
\end{equation}
so for $k$ large enough we have $E''_k\subset V$ and get that
\begin{equation}\label{equationapplication2P}
\phi_1(E''_k)\subset E'.
\end{equation}
To get the converse inclusion, notice that $J_h^2(E''_k)\geq J_h^2(E)=J_h^2(E')$ and since both $\ker^2(E''_k)$ and $\ker^2(E')$ are compact we get
\begin{equation}
\ker^2(E')\subset\phi_1^2(\ker^2(E''_k)).
\end{equation}
To be honest, the converse inclusion for the $1$-dimensional cores is a little more difficult to obtain if we consider $E^1$ as given by theorem~\ref{theoremexistence}. However, we can suppose that we minimized the measure of the $1$-dimensional core on the complementary of an open neighborhood of $E^2$ containing $V^1$ amongst the sets $F\in\mathfrak F$ such that $\ker^2(F)=E^2$, and use the same argument as above. Then again, since $E^0$ is finite the case of the $0$-dimensional core is easily treated, and we get as expected
\begin{equation}\label{equationapplication2Q}
E'\subset\phi_1(E''_k).
\end{equation}
Together with~\eqref{equationapplication2P} this achieves proving that
\begin{equation}
E'\in\mathfrak{F}.
\end{equation}

Notice that we did not prove that $E\in\mathfrak{F}$, because it was not needed in order to prove corollary~\ref{corollaryapplication2}. However, although the author feels quite inclined to believe that it is possible, it seems difficult to build a similar retraction that does not change anything to $E$ using the mere biHölder regularity of almost-minimal sets. Nonetheless, it seems easier when $n=3$ using Taylor and Morgan's versions which give a biLipschitz equivalence and thus, more control on the way $E^1$ meets $E^2$.

Remember that we still have to prove that we can build a local Hölder retraction on any $l$-dimensional reduced minimal cone for $l=1,2$ that meet the requirements~\eqref{equationapplication2K}, \eqref{equationapplication2L} and \eqref{equationapplication2M} as we announced before. 

Let us deal with the case $l=1$ first and suppose that $Z$ is a $1$-dimensional reduced minimal cone. For convenience, we also suppose that $Z$ is centered at the origin. According to~\cite{morgan} or~\cite{david2}, $Z$ can come in two flavors:
\begin{itemize}
\item a line, in that case we simply take the orthogonal projection onto it;
\item three half lines contained in a $2$-plane $P$ that meet at the origin and make $\frac{2\pi}3$ angles. In that case, denote by $p$ the orthogonal projection onto $P$. Notice that for all $z\in P\setminus Z$, the connected component of $P\setminus Z$ that contains $z$ is bounded by two of the three half lines in $Z$. Denote by $L$ the remaining half line without the origin and notice that the line through $z$ parallel to $L$ meets $Z\setminus L$ at an unique point (see figure~\ref{figureapplication2A}). Call it $q(z)$ and set $q(z)=z$ if $z\in Z$. It is easy to check that $q\circ p$ is Lipschitz, and meets all our requirements.
\end{itemize}

\begin{figure}\label{figureapplication2A}
\centering\framebox{\includegraphics[width=0.5\textwidth]{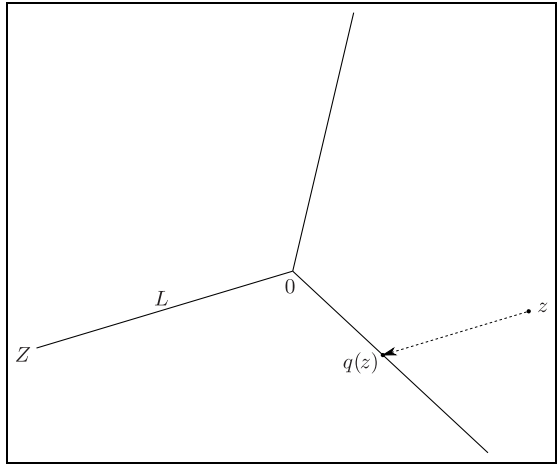}}
\caption{A simple Lipschitz retraction onto a $Y$-shaped minimal cone of dimension~$1$.}
\end{figure}

Now suppose that $Z$ is a $2$-dimensional reduced minimal cone centered at the origin. Then we know (again, see~\cite{david2}) that $Z$ is the cone over a set $K=K_1\sqcup K_2\subset\partial B(0,1)$ such that
\begin{itemize}
\item $K_1$ is a finite union of disjoint great circles,
\item $K_2$ is a finite union of closed arcs of great circles that only meet at their endpoints with $\frac{2\pi}3$ angles, and each endpoint is common to exactly three arcs.
\end{itemize}
Notice that according to this description, $K_1$ and $K_2$ are locally biLipschitz equivalent respectively to the two flavors of $1$-dimensional minimal cones we described previously. Since they are also compact and disjoint, we can build a Lipschitz map $q\colon\partial B(0,1)\rightarrow\partial B(0,1)$ that sends an open neighborhood of $K$ onto $K$. We can also extend $q$ to $\mathbb R^n$ by putting
\begin{equation}
q(z)=\begin{cases}
\Vert z\Vert q\left(\frac z{\Vert z\Vert}\right)\text{ if }z\neq 0\\
0\text{ otherwise}.
\end{cases}
\end{equation}
Finally (by using Kirszbraun theorem for instance), consider a Lipschitz map $p$ such that
\begin{equation}
p(z)=\begin{cases}
0\text{ if }\Vert z\Vert<\frac r2\\
z\text{ if }\Vert z\Vert>r.
\end{cases}
\end{equation}
Again, it is easy to check that $q\circ p$ meets our requirements.
\end{proof}

\ifCRM
  \vfill
  \section*{Acknowledgements}
  The author would like to express his thanks to Guy David for his many advices and suggestions. He also gratefully acknowledges partial support from the Centre de Recerca Matemàtica at the Universitat Autònoma de Barcelona.
\fi

\ifCRM
  \pagebreak
\else
  \newpage
  \addcontentsline{toc}{section}{\refname}
\fi
\bibliographystyle{alpha}
\bibliography{\jobname}

\end{document}